\documentclass[12pt]{amsart}%{elsarticle}
\usepackage[foot]{amsaddr}
\usepackage{amsmath,amsfonts}
\usepackage{verbatim}
\usepackage{lscape}
\usepackage{mathtools}
\usepackage{tikz}
\usepackage{algorithm}
\usepackage{algorithmic}
\usepackage{booktabs}
\usepackage{subcaption}
\newcommand{\norm}[1]{\left\lVert#1\right\rVert}

\newcommand{\paren}[1]{\left(#1\right)}
\newcommand{\cparen}[1]{\left[#1\right]}
\newcommand{\llparen}[1]{\left\{#1\right\}}

\newcommand*{\defeq}{\mathrel{\vcenter{\baselineskip0.5ex \lineskiplimit0pt
                     \hbox{\scriptsize.}\hbox{\scriptsize.}}}%
                     =}

\DeclareMathOperator{\skk}{sk}

\DeclareMathOperator{\SO}{\textsf{SO}}

\DeclareMathOperator{\bPi}{\mathbf{P}_{\boldsymbol{\tau}}}

\DeclareMathOperator{\btau}{\boldsymbol\tau}
\DeclareMathOperator{\bsig}{\boldsymbol\sigma}

\DeclareMathOperator{\bLamb}{\boldsymbol\Lambda}

\DeclareMathOperator{\bomega}{\boldsymbol\omega}

\DeclareMathOperator{\card}{card}

\DeclareMathOperator{\Rey}{Re}
\usetikzlibrary{backgrounds}
\graphicspath{{figures/}}

\begin{document}

%\begin{frontmatter}

\title{Simulating squirmers with volumetric solvers}%
%\title{Elsevier \LaTeX\ template\tnoteref{mytitlenote}}





\author[SP]{Stevens Paz}%\corref{cor1}}
\email{espeisan@usp.br}
%\ead{espeisan@usp.br}
%\cortext[cor1]{Corresponding author.}
\author[GCB]{Gustavo C. Buscaglia}
\email{gustavo.buscaglia@icmc.usp.br}
%\ead{gustavo.buscaglia@icmc.usp.br}
\address[SP,GCB]{Instituto de Ci\^encias Matem\'aticas e de Computa\c{c}\~ao – ICMC, Universidade de S\~ao Paulo, Campus de S\~ao Carlos, Caixa Postal 668, 13560-970 S\~ao Carlos, SP, Brazil}

\begin{abstract}
Squirmers are models of a class of microswimmers, such as ciliated org\-anisms and phoretic particles, that self-propel in fluids without significant deformation of their body shape. Available techniques for their simulation are based on the boundary-element method and do not contemplate nonlinearities such as those arising from the fluid's inertia or non-Newtonian rheology. This article describes a methodology to simulate squirmers that overcomes these limitations by using volumetric numerical methods, such as finite elements or finite volumes. It deals with  interface conditions at the squirmer's surface that generalize those in the published literature. The actual procedures to be performed on a fluid solver to implement the proposed methodology are provided, including the treatment of metachronal surface waves. Among the several numerical examples, a two-dimensional simulation is shown of the hydrodynamic interaction of two individuals of {\em Opalina ranarum}.
\end{abstract}

%\begin{keyword}
%squirmer model; numerical microfluidics; ciliated organisms; phoretic particles; fluid-solid interaction; finite ele\-ment/volume methods.
%\MSC[2010] 00-01\sep  99-00
%\end{keyword}

%\end{frontmatter}

%\linenumbers
\maketitle


%%%%%%%%%%%%%%%%%%%%%%%%%%%%%%%%%%%%%%%%%%%%%%%%%%
%%%%%%%%%%%%%%%%%%%%%%%%%%%%%%%%%%%%%%%%%%%%%%%%%%

\section{Introduction} 
Microswimmers are organisms or particles with self-driven capacity of locomotion \cite{lauga2014locomotion}. A large class of microswimmers is that of ciliated organisms, in which cilia act as oars that bend, stretch and rotate generating forces and displacements in the surrounding fluid \cite{childress2012natural,khanna2004biology}. A {\em squirmer}, initially introduced by Lighthill \cite{lighthill1952squirming}, is a model of a microswimmer consisting of a deformable body that swims via small shape oscillations \cite{debnath2008sir}. It was applied to ciliates by Blake \cite{blake1971spherical} using the concept of \textit{ciliary envelope}, in which the tips of the numerous cilia are treated as a deformable shell that covers the body. This model has been extensively used in the literature to study energy dissipation and swimming efficiency \cite{lighthill1975mathematical,michelin2010efficiency,kreissl2016efficiency}, nutrient uptake \cite{magar2003nutrient, magar2005average,michelin2011optimal,michelin2013unsteady,ishikawa2016nutrient} and the mechanical effect of the squirmer's geometry while swimming \cite{shapere1989geometry}.

Within the ciliary envelope model, the microswimmer has a smooth effective impermeable surface $\Gamma$ through which it interacts with the surrounding fluid. We restrict here to the important class of {\em tangential squirmers}, in which only the tangential motions of the envelope are considered \cite{giacche2010hydrodynamic}. Since normal-to-the-surface deformations are neglected, tangential squirmers move as rigid bodies that exhibit a tangential {\em slip velocity} ${\bf u}_s$ with respect to the adjacent fluid. To sustain the slippage between the body and the fluid a tangential force $\mathbf{f}_s$ (per unit area) develops at $\Gamma$ which consumes a power (per unit area) equal to $\mathbf{f}_s\cdot {\bf u}_s$. The organism must provide this power to the cilia at each point so that they can sustain their motion, which is the reason for being considered an {\em active particle}.

The mathematical treatment of squirmers has mainly dealt with those of spherical shape, to which analytical or semi-analytical (i.e., series expansion)  techniques can be applied \cite{lighthill1952squirming,blake1971spherical,wang2012inertial,evans2011orientational,khair2014expansions,pedley2016spherical}. Numerical approximations are needed to predict the motion of confined squirmers, of non-spherical squirmers, of squirmers interacting with other squirmers or other particles \cite{ishikawa2006hydrodynamic}, etc. The most frequent technique in the literature is the boundary element method \cite{kanevsky2010modeling,ishimoto2013squirmer,zhang2015paramecia,ishimoto2017boundary}, which expresses the velocity field in terms of Stokeslets (Green functions of the Stokes operator \cite{happel1983low}) \cite{lauga2009hydrodynamics,spagnolie2012hydrodynamics,zhu2013low,lauga2014locomotion,giuliani2018predicting}. Boundary element methods are attractive because only the squirmer's surface need to be meshed, and there is (essentially) no need of remeshing along the squirmers evolution, no matter how  large its displacements or rotations may be.

It should be noted that rigid bodies exhibiting active (power consuming) tangential velocity slippage with the adjacent fluid are not exclusive to ciliated organisms. Active phoretic particles \cite{anderson1989colloid,julicher2009generic,popescu2010phoretic}, % (Annual Review Fluid Mechanics 1989, incluir na pasta), 
such as Janus particles \cite{walther2013janus,de2015diffusiophoretic,zottl2016emergent}, 
are also modeled as rigid bodies with a tangential slip velocity ${\bf u}_s$, and the techniques described in this article apply to them as well \cite{shen2018hydrodynamic}.

There are some advantages in using finite element or finite volume methods to model squirmers. These methods are readily extended to non-New\-to\-nian rheology \cite{zhu2012self} and non-zero Reynolds numbers, whereas boundary elements rely on the problem's linearity. Further, finite elements/volumes provide a sparse representation of the volumetric velocity field for advection computations (e.g., of nutrients \cite{magar2003nutrient,magar2005average}), while boundary element results need to undergo a quite costly post-processing step. Finally, numerical analysis has over the years equipped these methods, especially finite element ones, with powerful a priori and a posteriori convergence assessment techniques, as well as with stabilization techniques \cite{codina2001stabilized}, that are less developed for their boundary elements counterparts. 

Notwithstanding, finite element/volume methods for squirmers are quite absent in the literature. The purpose of this contribution is to provide a fully detailed explanation of how to turn a finite element/volume Navier-Stokes solver into a squirmer simulator that contemplates squirmers of arbitrary shape and motion, and which allows general boundary conditions at the squirmer-fluid interface. To the authors' knowledge, all available volumetric formulations are restricted to the case in which ${\bf u}_s$ is imposed as a datum (``type-I case'' in what follows), though the ability of imposing the tangential force $\mathbf{f}_s$ (``type-II case'', with the force possibly depending on ${\bf u}_s$) is sometimes important \cite{short2006flows,kanevsky2010modeling,pedley2016spherical}. Restricting thus to precedents for the type-I case, Aguillon et al \cite{aguillon2012modelling} consider the squirmer problem discretized on a fixed mesh through the fictitious domain method \cite{glowinski2001fictitious}. Shen and Vernerey \cite{shen2017phoretic}, in their method for surface-active vesicles, also turn to a fixed-mesh technique by means of extended finite elements. Besides the generality of the boundary conditions, this contribution differs from these precedents in that the mesh conforms to the squirmers' boundaries in an Arbitrary-Lagrangian-Eulerian (ALE) manner \cite{sarrate2001arbitrary,hu2001direct,al2008motion,donea2017arbitrary}. This choice has well-known pros and cons. It carries with it a meshing difficulty, since the mesh needs to be deformed as the geometry changes in time and periodically rebuilt from scratch. This difficulty, however, is strictly a matter of computational geometry and the quality and availability of meshing software packages increases steadily. On the other hand, the ALE approach allows for {\em any} finite element/volume solver to be easily adapted to the squirming problem following some simple manipulations here described. Problems of particle sedimentation, of motion of very small bubbles and of swimming of articulated bodies can be addressed with variations of the proposed methodology.

The plan of this article is as follows: The mathematical formulation of the exact problem is developed in section 2, introducing useful notation for squirming kinematics and presenting the differential problems for both type-I and type-II squirmers, accompanied by the corresponding weak formulations. In section 3 the numerical method is given in full detail. The spatial discretization is worked out for the Galerkin finite element method (in 3.1-3.2) so as to put forward our specific implementation as an example. The rest of section 3 describes how to manipulate the matrices of a general fluid solver and perform the time marching so as to turn it into a squirmer simulator, and applies to essentially any nodal finite element solver or vertex-centered finite volume solver. In section 4 the verification of the method and code is reported by showing the results of convergence analyses and of comparisons with semi-analytical and numerical approximations in the literature. The verification is restricted to spherical steady squirmers, for which sufficient data are available. Steady squirmers, however, do not model the motion of ciliated bodies in their detailed dynamics since each cilium, being fixed to the surface, must describe an oscillatory motion. Section 5 explains the spatio-temporal organization of these oscillations that lead to self-propulsion, known as metachronal waves, and details their implementation as boundary conditions for squirmers of type I or II. This section is closed with a simulation of the interaction of two ciliated bodies inspired in the {\em Opalina ranarum}. The geometry is simplified to two dimensions for lack of a 3D remeshing algorithm in our implementation, but the techniques are described for the 3D case. Conclusions and suggestions for future work are compiled in section 6.

%%%%%%%%%%%%%%%%%%%%%%%%%%%%%%%%%%%%%%%%%%%%%%%%%%
%%%%%%%%%%%%%%%%%%%%%%%%%%%%%%%%%%%%%%%%%%%%%%%%%%

\section{Problem formulation}
The squirmers considered in this article move as rigid bodies, in the sense that for each one there exists a (closed) reference domain $\mathcal{B}^{*}\subset\mathbb{R}^{d}$ and all possible configurations of the squirmer are translations and rotations of $\mathcal{B}^*$. Taking an arbitrary point $\mathbf{X}_c$ as center of rotation, at all times $t$ there exists a point $\mathbf{x}_c(t)$ and a rotation matrix $\mathbf{Q}(t)$ such that the position $\mathbf{x}(\mathbf{X},t)$ of the material point $\mathbf{X}$ is given by
\begin{equation}\label{eqn:rig_motion_iso}
\mathbf{x}\left(\mathbf{X},t\right) = \mathbf{x}_{c}\paren{t} + \mathbf{Q}\paren{t}\left(\mathbf{X}-\mathbf{X}_{c}\right)~.
\end{equation}
Notice that vectors in $\mathbb{R}^d$ ($d=2$ or 3) operate as column matrices in the algebraic equations. The region occupied by the squirmer at time $t$ is, thus, 
\[%\begin{equation}
\mathcal{B}\paren{t}=\{\mathbf{x}\left(\mathbf{X},t\right)~, \mathbf{X}\,\in\,\mathcal{B}^* \}~.
\]%\end{equation}
We assume that the squirmer moves inside a fixed domain $\Omega$ filled with fluid, so that $\Omega_{f}\paren{t} = \Omega\setminus{\mathcal{B}}\paren{t}$ is the fluid domain at time $t$. For simplicity, in the exposition we take zero-velocity boundary conditions for the fluid at $\partial\Omega$. Other boundary conditions are dealt with in the usual way. In a nutshell, the squirming problem consists of finding a continuous function $[0,T]\to\mathbf{q}\paren{t} = \left(\mathbf{x}_{c}\paren{t},\mathbf{Q}\paren{t}\right)$ of generalized coordinates satisfying a given initial condition and suitable interface conditions with the surrounding fluid at $\partial \mathcal{B}(t)$.
%%%%%%%%%%%%%%%%%%%%%%%%%%%%%%%%%%%%%%%%%%%%%%%%%%
\subsection{Some useful notation for squirmer kinematics}

Because $\mathbf{Q}(t)$ belongs to the \textit{special rotation group}
\[
\SO\paren{d} = \llparen{\mathbf{Q}\in\mathbb{R}^{d\times d}\,:\, \mathbf{Q}^{-1} = \mathbf{Q}^{T}  ,\, \det\cparen{\mathbf{Q}} = 1},
\]
the manifold of possible configurations of the squirmer is $Q=\mathbb{R}^d\times \SO\paren{d}$ \cite{betsch2016structure,lew2016brief,shabana2013dynamics} and
the body's Eulerian velocity $\mathbf{u}_{\mathcal{B}}$ is given by
\begin{equation}\label{eqn:kin_rest}
\mathbf{u}_{\mathcal{B}}\paren{\mathbf{x},t} = \mathbf{v}_c\paren{t} + \boldsymbol{\omega}\paren{t}\times\left(\mathbf{x}-\mathbf{x}_{c}\paren{t}\right)~,
\end{equation}
where $\mathbf{v}_{c}\paren{t} = \mathbf{\dot{x}}_{c}\paren{t}$ is the translational velocity and $\boldsymbol{\omega}\paren{t}$ is the pseudovector of angular velocities in the spatial frame. It relates to $\mathbf{Q}(t)$ and $\mathbf{\dot{Q}}(t)$ by 
\begin{equation}\label{eqn:Qode}
\mathbf{\dot{Q}}\mathbf{Q}^{T} = \skk\cparen{\boldsymbol{\omega}} \defeq
\begin{pmatrix}
0 & -\omega_{z} & \omega_{y}\\
\omega_{z} & 0 & -\omega_{x}\\
-\omega_{y} & \omega_{x} & 0
\end{pmatrix}~,
\end{equation}
where the isomorphism $\skk[\cdot ]$ between vectors and skew-symmetric matrices has been introduced.
Notice that $\skk\cparen{\bomega}\mathbf{y} = \boldsymbol\omega\times\mathbf{y}$ for all $\mathbf{y}\in\mathbb{R}^{3}$. For $d = 2$ the generalized coordinates can be changed to $\mathbf{q}(t)=(\mathbf{x}_c(t),\theta(t))$, replacing the rotation matrix $\mathbf{Q}$ by the rotation angle $\theta\in\mathbb{R}$ of which the time derivative $\omega = \dot{\theta}$ is the rotational velocity of the body. The expression for $\mathbf{u}_{\mathcal{B}}$ in such a case  simplifies to
\begin{equation}\label{eqn:kin_rest_2d}
\mathbf{u}_{\mathcal{B}}\paren{\mathbf{x},t} = \mathbf{v}_c\paren{t} + \omega\paren{t}\bLamb\left(\mathbf{x}-\mathbf{x}_{c}\paren{t}\right)~,
\end{equation}
with 
\[
\mathbf{Q}\paren{t} = 
\begin{pmatrix}
\cos\theta(t) & -\sin\theta(t)\\
\sin\theta(t) & \cos\theta(t)
\end{pmatrix}\quad\mbox{ and } \quad
\bLamb = 
\begin{pmatrix}
0 & -1\\
1 & 0
\end{pmatrix}.
\]
Let us define now the {\em velocity array}
\begin{equation}\label{eqn:defs}
\mathbf{s} = \left ( \begin{array}{c}\mathbf{v}_{c}\\
\boldsymbol\omega \end{array} \right )\in\mathbb{R}^{n_{c}}, \quad n_{c} = d+\frac{d\paren{d-1}}{2}~.
\end{equation}
If $d=3$, $\mathbf{s}$ relates to $\mathbf{\dot{q}}=\paren{ \mathbf{v}_{c}, \mathbf{\dot{Q}} }$ through \eqref{eqn:Qode}. If $d=2$ we simply have $\mathbf{s}=\mathbf{\dot{q}}$.
Given a trajectory $\mathbf{q}\paren{t}$, equations \eqref{eqn:kin_rest}, \eqref{eqn:kin_rest_2d} and \eqref{eqn:defs} allow us to make explicit the linear dependence of the body's velocity field with the velocity vector, i.e.,
\[%\begin{equation}
\mathbf{u}_{\mathcal{B}}(\mathbf{x},t) = \mathbf{H}(\mathbf{q}\paren{t},\mathbf{x})\,\mathbf{s}\paren{t},
\]%\end{equation}
where we have introduced the matrix $\mathbf{H}(\mathbf{q}\paren{t},\mathbf{x})\,\in\,\mathbb{R}^{d\times n_c}$, of which the first $d$ columns (corresponding to pure translations) are the identity matrix $\mathbf{I}_d$ and the next $n_c-d$ columns (corresponding to pure rotations) are $-\skk\cparen{\mathbf{x} - \mathbf{x}_{c}\paren{t}}$ if $d=3$, or $\bLamb\cparen{\mathbf{x} - \mathbf{x}_{c}\paren{t}}$ if $d=2$, this is, 
\[%\begin{equation}
\mathbf{H}(\mathbf{q}\paren{t},\mathbf{x}) = 
\begin{cases}
\paren{\mathbf{I}_{d}~~|~~ -\skk\cparen{\mathbf{x} - \mathbf{x}_{c}\paren{t}}} &\text{if $d = 3$}\\
& \\%\cparec{\paren{\theta\paren{t}}\paren{\mathbf{x}^{*} - \mathbf{x}^{*}_{c}}}}
\paren{\mathbf{I}_{d}~~ |~~ \phantom{-}\bLamb\cparen{\mathbf{x} - \mathbf{x}_{c}\paren{t}}} &\text{if $d = 2$}.
\end{cases} 
\]%\end{equation}
The previous notation readily extends to the case of $N>1$ squirmers by defining $N$ sets of generalized coordinates $\mathbf{q}$, generalized velocities $\mathbf{\dot{q}}$, velocity vectors $\mathbf{s}$, etc., since each body will follow an independent rigid motion.
%%%%%%%%%%%%%%%%%%%%%%%%%%%%%%%%%%%%%%%%%%%%%%%%%%
\subsection{The fluid problem}

The ambient fluid is assumed incompressible, so that its governing equations are given by
\[%\begin{equation}\label{eqn:ns_eqn}
\rho\frac{D \mathbf{u}}{D t} - \nabla\cdot\bsig = \mathbf{0}, \quad 
\hfil\nabla\cdot\mathbf{u} = 0,
\quad\text{in }\Omega_{{f}}\left(t\right),\,t\in\left(0,T\right),
\]%\end{equation}
where $\rho$ is the density, $\mathbf{u}$ the Eulerian velocity field, $D/Dt$ the material derivative and $\bsig$ the Cauchy stress tensor. The fluid will be assumed Newtonian for simplicity, i.e.,
\[%\begin{equation}
    \bsig = - p\, \mathbf{I}_d + 2\mu\nabla^S\mathbf{u},
\]%\end{equation}
where $p$ is the pressure, $\mu$ is the viscosity and $\nabla^S\mathbf{u}=\frac12 \paren{\nabla \mathbf{u} +\nabla \mathbf{u}^T}$,
but other rheological models can be considered. The presence of the squirmer (or squirmers) intervenes through the flow domain $\Omega_f$, since $\Omega_f(t)=\Omega \setminus \mathcal{B}\paren{t}$, and through the kinematical and dynamical compatibility conditions at $\partial\mathcal{B}\paren{t}$. These are:
\begin{itemize}
    \item {\em Kinematical condition:} There exists a tangential slip velocity $\mathbf{u}_s$ between the body and the adjacent fluid, i.e.,
    \begin{equation}\label{eqn:slip_bc1}
\mathbf{u}\paren{\mathbf{x},t} = \mathbf{u}_{\mathcal{B}}\paren{\mathbf{x},t} + \mathbf{u}_{s}\paren{\mathbf{x},t},\quad \forall\mathbf{x}\in\partial\mathcal{B}\paren{t},~\forall t.
    \end{equation}
    Notice that the normal velocity is continuous, since the squirmer's surface is impermeable. For ciliated organisms, the slip velocity represents the velocity difference between the real surface of the body and that of the surrounding fluid near the cilia's tips. For electrophoretic particles, on the other hand, $\mathbf{u}_s$ represents the jump in velocity across the (nanometric) electric double layer.
    \item {\em Tangential force equilibrium:} The force $\mathbf{f}_s$ is exerted by the cilia on the adjacent fluid, i.e.,
    \[%\begin{equation} \label{eqn:tanforce}
        \bPi \,\bsig\paren{\mathbf{x},t}\, \mathbf{n}= \mathbf{f}_s\paren{\mathbf{x},t},\quad \forall\mathbf{x}\in\partial\mathcal{B}\paren{t},~\forall t,
    \]%\end{equation}
    where $\bPi=\mathbf{I}_d-\mathbf{n}\,\mathbf{n}^T$ is the projection matrix onto the tangent plane to $\partial\mathcal{B}(t)$ at $\mathbf{x}$. The normal unit vector $\mathbf{n}$ points into the body.
    \item {\em Global force and torque balance:} Neglecting the inertia of the squirmer, the total force and torque on it must be zero. That is, for all $t$,
    \begin{eqnarray}
    \int_{\partial\mathcal{B}(t)} \bsig\mathbf{n}\ dS &=& \mathbf{0}, \label{eqn:forcefree}\\
    \int_{\partial\mathcal{B}(t)} \paren{\mathbf{x}-\mathbf{x}_{c}}\times\bsig\mathbf{n}\ dS &=& \mathbf{0}.\label{eqn:torquefree}
    \end{eqnarray}
    If the inertia of the squirmer is considered, the right-hand sides above must change appropriately (e.g., if $\mathbf{X}_c$ is chosen as the center of mass the right-hand side of \eqref{eqn:forcefree} changes to $M\,\mathbf{\dot{v}}_c$, $M$ being the body's mass). We concentrate here in cases with negligible inertia not just because they are physically realistic, but also because the usual algorithms for fluid-structure interaction (weak coupling, iterative coupling) cannot be applied at all. The methods proposed here can be readily extended to consider the inertia in an implicit, strongly coupled way.
\end{itemize}

Though conditions \eqref{eqn:slip_bc1}-\eqref{eqn:torquefree} must necessarily hold for the solution to be physically meaningful, the two quantities ${\bf u}_s$ and ${\bf f}_s$ cannot be simultaneously imposed as data of the problem. Just as what happens in a heat conduction problem, in which one can impose the boundary temperature or the heat flux, but not both, in the squirming problem one can impose the slip velocity or the tangential force, but not both. This gives rise to two kinds of squirmers. Those in which $\mathbf{u}_s$ is given will be denoted here as {\em type-I} squirmers. This is the case considered in practically all previous studies. The squirmers in which $\mathbf{f}_s$ is given, as a known quantity or as a known function of $\mathbf{u}_s$, will be referred to as {\em type-II} squirmers. As discussed by Short et al \cite{short2006flows}, type-II squirmers can help model organisms in which data of the effective slip velocity are unavailable. The mathematical formulations for type-I and type-II squirmers are somewhat different and are thus presented separately below.
%%%%%%%%%%%%%%%%%%%%%%%%%%%%%%%%%%%%%%%%%%%%%%%%%%
\subsection{Type-I squirmer}
In this case one imposes the slip velocity, which is given by a known {\em tangent} vector $\mathbf{u}_{s}^{*}(\mathbf{X},t)$ in the material frame. Specifically, if  $\mathbf{x}\in\partial\mathcal{B}(t)$, the slip velocity is a datum calculated from
\[%\begin{equation}\label{eqn:typeIa}
\mathbf{u}_{s}\paren{\mathbf{x},t} = \mathbf{u}_{s}\paren{\mathbf{x}\paren{\mathbf{X},t},t} = \mathbf{Q}\paren{t}\mathbf{u}_{s}^{*}\paren{\mathbf{X},t}, \quad \mathbf{X}\in\partial\mathcal{B}^{*}.
\]%\end{equation}
It can also be given as a scalar field $u_s$ which is multiplied by the unit tangent vector at each instant to obtain $\mathbf{u}_s$.
The mathematical problem (for $d=3$, the case $d=2$ is an easy exercise) reads as follows: Given $\mathbf{q}(t=0)$ and $\mathbf{u}(\mathbf{x},t=0)$ (this latter datum is only needed if $\rho>0$), determine $\mathbf{q}(t)=\paren{\mathbf{x}_c(t),\mathbf{Q}(t)}$, $\mathbf{s}(t)=\paren{\mathbf{v}_c(t),\bomega(t)}$, $\mathbf{u}(\mathbf{x},t)$ and $p(\mathbf{x},t)$ for $0<t\leq T$ and $\mathbf{x}\,\in\,\Omega_f(t)$ satisfying
\begin{eqnarray}
\frac{d\mathbf{x}_c}{dt}&=&\mathbf{v}_c, \label{eqtypeIa}\\
\frac{d\mathbf{Q}}{dt}&=&\skk\cparen{\bomega}\,\mathbf{Q},\label{eqtypeIb}\\
\mathbf{u}\paren{\mathbf{x},t} - \mathbf{H}\paren{\mathbf{q}(t),\mathbf{x}}\,\mathbf{s}(t) &=& \mathbf{u}_{s}\paren{\mathbf{x},t},\qquad \mbox{ on }\partial\mathcal{B}\paren{t},\label{eqtypeIc}\\
\rho\frac{D \mathbf{u}}{D t} - \mu\,\nabla^2\mathbf{u} + \nabla p &=& \mathbf{0}, \qquad\qquad \mbox{ in } \Omega_f(t),\label{eqtypeId}\\
\nabla\cdot\mathbf{u} &=& 0,\qquad\qquad\mbox{ in } \Omega_f(t),\label{eqtypeIe}\\
    \int_{\partial\mathcal{B}(t)} \bsig\mathbf{n}\ dS &=& \mathbf{0},\label{eqtypeIf}\\
    \int_{\partial\mathcal{B}(t)} \paren{\mathbf{x}-\mathbf{x}_{c}}\times\bsig\mathbf{n}\ dS &=& \mathbf{0}.\label{eqtypeIg}
\end{eqnarray}
Considering for interpretation purposes $\rho=0$, we see that the main equations to be solved are \eqref{eqtypeIa}-\eqref{eqtypeIb}, of which the right-hand side contains the {\em unique} values of $\mathbf{s}=\paren{\mathbf{v}_c,\bomega}$ that introduced in \eqref{eqtypeIc} impose velocity boundary conditions for the Navier-Stokes equations \eqref{eqtypeId}-\eqref{eqtypeIe} that produce a force-free and torque-free solution. Notice that $\mathbf{q}=\paren{\mathbf{x}_c,\mathbf{Q}}$ intervenes in \eqref{eqtypeIc}-\eqref{eqtypeIg} not just explicitly (in \eqref{eqtypeIc}) but also through the geometry (i.e., $\Omega_f$, $\partial\mathcal{B}$). The problem clearly belongs to the class of fluid-solid-interaction ones, with negligible inertia in the solid.

The weak form of the problem above can be derived as usual by multiplication by a test function and integration by parts. The variational problem is formulated on the space \cite{glowinski2001fictitious,walker2013analysis}
\[%\begin{equation}%\label{eqn:mech_rest_D}
%\begin{split}
W\paren{\mathbf{q}} = \Big\{ \mathbf{w}\in H^{1}\paren{\Omega_{f}(\mathbf{q})}^{d}\, : \mathbf{w} = \mathbf{0} \text{ on } \partial\Omega,~\mathbf{w} = \mathbf{H}({\mathbf{q}})\mathbf{d} \text{ on } \partial\mathcal{B}, ~\mathbf{d}\,\in\,\mathbb{R}^{n_c} \Big\},
%\end{split}
\]%\end{equation}
where the dependence on time has been replaced by a dependence on $\mathbf{q}$, which depends on time. Let us now introduce an extension (or lifting) linear operator $\mathcal{E}$ that, given a (regular enough) function $f$ defined on $\partial \mathcal{B}$, assigns to it $\mathcal{E}f\in H^1(\Omega_f)$ that coincides with $f$ on $\partial\mathcal{B}$ and is zero on $\partial\Omega$. The action of this operator on vector or matrix fields defined on $\partial\mathcal{B}$ is defined by applying $\mathcal{E}$ componentwise.

The space $W(\mathbf{q})$ then decomposes as
\[%\begin{equation}
    W(\mathbf{q})=W_0(\mathbf{q})\oplus V(\mathbf{q})  
\]%\end{equation}
where $V(\mathbf{q})$ is the finite-dimensional space (of dimension $n_c$) of extensions of rigid-body motions, i.e.,
\begin{equation} \label{eqvq}
    V(\mathbf{q})=\{ \mathbf{w}=\mathcal{E}\,\mathbf{H}(\mathbf{q})\,\mathbf{d}, ~\mathbf{d}\in\mathbb{R}^{n_c}\}
\end{equation}
and $W_0(\mathbf{q})=H^1_0(\Omega_f)^d$ consists of vector fields that vanish at all the boundaries. Above, $\mathbf{H}(\mathbf{q})$ is shorthand for the matrix field $\mathbf{H}(\mathbf{q},\cdot)$ and $\mathbf{H}(\mathbf{q})\mathbf{d}$ is the vector field defined on $\partial\mathcal{B}$ by $\cparen{\mathbf{H}(\mathbf{q})\mathbf{d}}\paren{\mathbf{x}}=\mathbf{H}(\mathbf{q},\mathbf{x})\,\mathbf{d}$. In \eqref{eqvq} the operator $\mathcal{E}$ acts on this field, or, equivalently, the matrix field $\mathcal{E}\mathbf{H}(\mathbf{q})$ acts on the vector $\mathbf{d}$, since
$$
\mathcal{E} \left ( \mathbf{H}(\mathbf{q})\mathbf{d} \right ) = 
\left ( \mathcal{E}\mathbf{H}(\mathbf{q}) \right )\,\mathbf{d}.
$$

The matrix field
\[%\begin{equation}
    \widetilde{\mathbf{H}}(\mathbf{q})=\mathcal{E}\mathbf{H}(\mathbf{q})
\]%\end{equation}
plays an important role in the picture. 
A {\bf basis for $V(\mathbf{q})$}, that we will denote by $\{\mathbf{h}^i\in V(\mathbf{q}),i=1,\ldots,n_c\}$, is provided by the {\bf columns of $\widetilde{\mathbf{H}}$} considered as vector fields on $\Omega_{f}$.

The weak form of \eqref{eqtypeIc}-\eqref{eqtypeIg} has a remarkably compact form: {\em Find $(\mathbf{s}(t),\mathbf{u},p)$, where 
$\mathbf{u}$ must belong to $\mathcal{E}\mathbf{u}_s+\widetilde{\mathbf{H}}(\mathbf{q}(t))\,\mathbf{s}(t)+W_0(\mathbf{q}(t))$ and $p$ must belong to $L^2_0(\Omega_f(t))$, such that
\begin{eqnarray*}
    \int_{\Omega_{f}\paren{t}} \rho\,\frac{D\mathbf{u}}{Dt}\paren{\mathbf{x},t}\cdot\mathbf{w}\paren{\mathbf{x}}\, d\mathbf{x} + \int_{\Omega_{f}\paren{t}} \bsig\paren{\mathbf{x},t}:\nabla^{S}\mathbf{w}\paren{\mathbf{x}}\, d\mathbf{x}&=&0,\\%\label{eqmom1}\\
\int_{\Omega_f\paren{t}} z\paren{\mathbf{x}}\nabla \cdot \mathbf{u}\paren{\mathbf{x},t}\, d\mathbf{x} &=& 0,%\label{eqmass1}
\end{eqnarray*}
for all $\paren{\mathbf{w},z}\in W\paren{\mathbf{q}(t)}\times L^2_0(\Omega_f(t))$.} By taking, with $t$ fixed, $\mathbf{w}=\mathbf{u}(\cdot,t)-\mathcal{E}\,\mathbf{u}_s$ and using \eqref{eqtypeId}-\eqref{eqtypeIe} together with the Reynolds transport theorem, one gets the energy identity
\begin{equation}\label{eqenergy}
    \frac{d}{dt}\int_{\Omega_f(t)} \rho \frac{\|\mathbf{u}\|^2}{2}~d\mathbf{x}
    + \int_{\Omega_f(t)} 2\mu \|\nabla^S\mathbf{u}\|^2~d\mathbf{x}
    = \int_{\partial \mathcal{B}(t)} \mathbf{f}_s\cdot \mathbf{u}_s ,
\end{equation}
which shows that the power spent by the cilia at $\partial\mathcal{B}$ sustains the motion of the squirmer against viscous dissipation (second term above).
%%%%%%%%%%%%%%%%%%%%%%%%%%%%%%%%%%%%%%%%%%%%%%%%%%
\subsection{Type-II squirmer}
In this case one imposes the {\em tangential} force $\mathbf{f}_s$ exerted by the cilia on the fluid, while $\mathbf{u}_s=\mathbf{u}-\mathbf{u}_{\mathcal{B}}$ is an unknown of the problem.

The mathematical problem (for $d=3$, the case $d=2$ is an easy exercise) reads as follows: Given $\mathbf{q}(t=0)$ and $\mathbf{u}(\mathbf{x},t=0)$ (this latter datum is only needed if $\rho>0$), determine $\mathbf{q}(t)=\paren{\mathbf{x}_c(t),\mathbf{Q}(t)}$, $\mathbf{s}(t)=\paren{\mathbf{v}_c(t),\bomega(t)}$, $\mathbf{u}(\mathbf{x},t)$ and $p(\mathbf{x},t)$ for $0<t\leq T$ and $\mathbf{x}\,\in\,\Omega_f(t)$ satisfying
\begin{eqnarray}
\frac{d\mathbf{x}_c}{dt}&=&\mathbf{v}_c, \label{eqtypeIIa}\\
\frac{d\mathbf{Q}}{dt}&=&\skk\cparen{\bomega}\,\mathbf{Q},\label{eqtypeIIb}\\
\mathbf{n}\cdot \left [\mathbf{u}\paren{\mathbf{x},t} - \mathbf{H}\paren{\mathbf{q}(t),\mathbf{x}}\,\mathbf{s}(t) \right ]&=& 0,\qquad \qquad \mbox{ on }\partial\mathcal{B}\paren{t},\label{eqtypeIIc}\\
\rho\frac{D \mathbf{u}}{D t} - \mu\,\nabla^2\mathbf{u} + \nabla p &=& \mathbf{0}, \qquad\qquad \mbox{ in } \Omega_f(t),\label{eqtypeIId}\\
\nabla\cdot\mathbf{u} &=& 0,\qquad\qquad\mbox{ in } \Omega_f(t),\label{eqtypeIIe}\\
\bPi \bsig \mathbf{n} &=& \mathbf{f}_s ,\qquad \qquad \mbox{ on }\partial\mathcal{B}\paren{t},\label{eqtypeIIc2}\\
    \int_{\partial\mathcal{B}(t)} \bsig\mathbf{n}\ dS &=& \mathbf{0},\label{eqtypeIIf}\\
    \int_{\partial\mathcal{B}(t)} \paren{\mathbf{x}-\mathbf{x}_{c}}\times\bsig\mathbf{n}\ dS &=& \mathbf{0}.\label{eqtypeIIg}
\end{eqnarray}
The difference with the equations of a type-I squirmer is that now the boundary conditions for the Navier-Stokes equations \eqref{eqtypeIId}-\eqref{eqtypeIIe} have been split into two: The normal component of the velocity is constrained by \eqref{eqtypeIIc}, while the tangential force is imposed by \eqref{eqtypeIIc2}.

The variational problem is formulated on the space
\[%\begin{equation}%\label{eqn:mech_rest_D}
%\begin{split}
\widetilde{W}\paren{\mathbf{q}} = \Big\{ \mathbf{w}\in H^{1}\paren{\Omega_{f}(\mathbf{q})}^{d}\, : \mathbf{w} = \mathbf{0} \text{ on } \partial\Omega,~\mathbf{n}\cdot \left [ \mathbf{w} - \mathbf{H}({\mathbf{q}})\mathbf{d} \right ] = 0 \text{ on } \partial\mathcal{B}, ~\mathbf{d}\,\in\,\mathbb{R}^{n_c} \Big\}.
%\end{split}
\]%\end{equation}
The space $\widetilde{W}(\mathbf{q})$ decomposes as
\[%\begin{equation}
    \widetilde{W}(\mathbf{q})= W_{\parallel}(\mathbf{q})\oplus V(\mathbf{q})
\]%\end{equation}
where $V(\mathbf{q})$ is as before (Eq. \eqref{eqvq}) and
\[%\begin{equation}
    W_{\parallel}(\mathbf{q})=\{ \mathbf{w}\,\in\,H^1(\Omega_f)^d,~\mathbf{w}=\mathbf{0}\text{ on }\partial\Omega,~\mathbf{w}\cdot \mathbf{n}=0\text{ on }\partial\mathcal{B}\}~.
\]%\end{equation}
As a consequence, each $\mathbf{w}\in\widetilde{W}(\mathbf{q})$ can be uniquely decomposed as 
\[%\begin{equation}
    \mathbf{w}=\mathbf{w}_{\parallel}+\mathbf{w}_V,
\]%\end{equation}
with $\mathbf{w}_{\parallel}\in W_\parallel(\mathbf{q})$ and $\mathbf{w}_V\in V(\mathbf{q})$.

The weak form of \eqref{eqtypeIIc}-\eqref{eqtypeIIg} is: {\em Find $(\mathbf{s}(t),\mathbf{u},p)$, where $\mathbf{u}$ must belong to $\widetilde{\mathbf{H}}(\mathbf{q}(t))\,\mathbf{s}(t)+W_\parallel(\mathbf{q}(t))$, $p$ must belong to $L^2_0(\Omega_f(t))$, and 
\begin{eqnarray}
    \int_{\Omega_{f}\paren{t}} \rho\,\frac{D\mathbf{u}}{Dt}\cdot\mathbf{w}\, d\mathbf{x} + \int_{\Omega_{f}\paren{t}} \bsig:\nabla^{S}\mathbf{w}\, d\mathbf{x}&=&\int_{\partial\mathcal{B}(t)}\mathbf{f}_s\cdot\mathbf{w}_\parallel,\label{eqmom2}\\
\int_{\Omega_f\paren{t}} z\,\nabla \cdot \mathbf{u}\, d\mathbf{x} &=& 0,\label{eqmass2}
\end{eqnarray}
for all $\paren{\mathbf{w},z}\in \widetilde{W}\paren{\mathbf{q}(t)}\times L^2_0(\Omega_f(t))$.} By taking $\mathbf{w}=\mathbf{u}(\cdot,t)$ in \eqref{eqmom2}
and noting that $\mathbf{u}_\parallel=\mathbf{u}_s=\mathbf{u}-\mathbf{u}_{\mathcal{B}}$
one gets that the energy equation \eqref{eqenergy} also holds for type-II squirmers. In this case, of course, $\mathbf{u}_s$ is not a datum but instead is computed from $\mathbf{u}_s(\mathbf{x},t)=\mathbf{u}(\mathbf{x},t)-\mathbf{H}(\mathbf{q}\paren{t},\mathbf{x})\mathbf{s}(t)$.

%%%%%%%%%%%%%%%%%%%%%%%%%%%%%%%%%%%%%%%%%%%%%%%%%%
%%%%%%%%%%%%%%%%%%%%%%%%%%%%%%%%%%%%%%%%%%%%%%%%%%

\section{Numerical method} 
%%%%%%%%%%%%%%%%%%%%%%%%%%%%%%%%%%%%%%%%%%%%%%%%%%
\subsection{Discretization in space}

The spatial discretization is only discussed here for the case of conforming finite elements, but it can be translated quite straightforwardly to other techniques. The exposition recovers its generality once the matrix formulation each type of squirmer is established, and is largely independent of the discretization method that led to it. 

Let us thus proceed to discretize the proposed problem in space. For this purpose, an ALE moving mesh is adopted below. In general, this strategy requires periodic remeshing and subsequent interpolation of the variables. Though implemented in our code and used in the examples, this issue will not be addressed here.

For each $t\in[0,T]$, $T>0$, let $\mathcal{T}_{h}\paren{t}$ be an approximate triangulation of the fluid region $\Omega_{f}\paren{t}$, this is, a regular partition of the physical domain into non-empty compact subdomains, or elements, $\Omega^{e}\paren{t}$ of characteristic size $h$, which define a discrete domain $\overline{\Omega}_{fh}\paren{t}\subset\overline{\Omega}\paren{t}$ as
\[
\overline{\Omega}_{fh}\paren{t} = \bigcup_{e} \Omega^{e}\paren{t}.
\]
We assume for simplicity that $\Omega$ is polygonal and thus $\partial\Omega$ is exactly approximated.
The interpolated boundary of the squirmer is denoted by $\partial\mathcal{B}_h(t)$, so that
\[
\partial\Omega_{fh}(t)=\partial\Omega\,\cup\,\partial\mathcal{B}_h(t)~.
\]

The fluid velocity $\mathbf{u}$ and pressure $p$ are approximated as
\[%\begin{subequations}\label{eqn:app_vp}
\begin{aligned}%\begin{align}
\mathbf{u}_{h}\left(\mathbf{x},t\right) &= \sum_{j\in\eta^{U}} \mathcal{N}^{j}\left(\mathbf{x},t\right)\mathbf{u}^{j}\left(t\right),\\%\label{eqn:app_vp_u}\\ %+ \sum_{K=1}^{N} \sum_{J\in\eta_{{K}}^{D}} \mathcal{N}^{{J}}\left(\mathbf{x},t\right) \left[ \mathbb{H}_{\mathbf{q}^{K}}\left(\mathbf{x}^{{J}},t\right)\mathbf{z}^{{K}}\left(t\right) + \mathbf{u}_{s,\btau}^{K}\left(\mathbf{x}^{J},t\right) \right],\\
p_{h}\left(\mathbf{x},t\right) &= \sum_{k\in\eta^{P}} \mathcal{M}^{{k}}\left(\mathbf{x},t\right)p^{{k}}\left(t\right),%\label{eqn:app_vp_p}
\end{aligned}%\end{align}
%\quad \mathbf{x}\in\overline{\Omega}_{fh}\paren{t}
\]%\end{subequations}
for $\mathbf{x}\in\overline{\Omega}_{fh}\paren{t}$, in finite dimensional subspaces $U_{h}\paren{t}\subset H^{1}\paren{\Omega_f\paren{t}}^{d}$
and $M_{h}\paren{t}\subset L^2_0(\Omega_{f}(t))$. The shape functions $\mathcal{N}^{j}\paren{\cdot,t}$, $\mathcal{M}^{k}\paren{\cdot,t}$ satisfy the nodal value property, namely, 
\[
\mathcal{N}^{j}\paren{\mathbf{x}^{i}(t),t} =
\begin{cases}
1,\quad \text{if $i = j$}\\
0,\quad \text{if $i \neq j$}
\end{cases},
\]
where $\mathbf{x}^{i}(t)$ is the position of node $i$ of the mesh $\mathcal{T}_{h}\paren{t}$, for $i$ belonging to the velocity global index set $\eta^{U}$. In particular, $\mathbf{u}_{h}\paren{\mathbf{x}^{i},t} = \mathbf{u}^{i}\paren{t}$, for all $i\in\eta^{U}$. Similarly, $\mathcal{M}^{k}\paren{\mathbf{x}^{l}(t),t}=\delta_{kl}$, so that
$p_{h}\paren{\mathbf{x}^{l}(t),t} = p^{l}(t)$ for pressure nodes $\mathbf{x}^{l}$ indexed by the set $\eta^{P}$.

Assuming the mesh to have no hanging nodes, the interpolation space for the velocity is
\[%\begin{equation}\label{eqn:vel_vh_rn}
%\begin{split}
U_{h}\paren{\mathbf{q}} = \Big\{ \mathbf{w}\in H^{1}\paren{\Omega_{fh}\paren{\mathbf{q}}}^{d}\, : \, \mathbf{w}\big|_{\Omega^{e}}\in P_{m}\paren{\Omega^{e}}^{d}, \text{ for all } e,~\mathbf{w}\big|_{\partial\Omega} = \mathbf{0} \Big\},
%\end{split}
\]%\end{equation}
and for the pressure
\[
M_{h}\paren{\mathbf{q}} = \left\{q\in L^{2}_{0}\paren{\Omega_{fh}\paren{\mathbf{q}}}\, : \, q\big|_{\Omega^{e}}\in P_{m}\paren{\Omega^{e}}, \text{ for all } e \right\},
\]
where $P_{m}\paren{\Omega^{e}}$ is the space of polynomials in $\Omega^{e}$ of degree less than or equal to $m$. In particular, we consider a stabilized $P_{1}/P_{1}$ element \cite{hughes1989new} and $P_{2}/P_{1}$ Taylor-Hood element \cite{taylor1973numerical}.
%%%%%%%%%%%%%%%%%%%%%%%%%%%%%%%%%%%%%%%%%%%%%%%%%%
\subsection{Semidiscrete Galerkin formulation for type-I squirmers}

Let us define
\[%\begin{equation}
%   \begin{split}
W_h\paren{\mathbf{q}} = \Big\{ \mathbf{w}_h\in U_h\paren{\mathbf{q}}\, : \mathbf{w}_h = \mathbf{0} \text{ on } \partial\Omega,~\mathbf{w}_h = \mathbf{H}({\mathbf{q}})\mathbf{d} \text{ on } \partial\mathcal{B}_h, ~\mathbf{d}\,\in\,\mathbb{R}^{n_c} \Big\}
%\end{split}
\]%\end{equation}
and
\[%\begin{equation}
%   \begin{split}
W_{0h}\paren{\mathbf{q}} = \Big\{ \mathbf{w}_h\in U_h\paren{\mathbf{q}}\, : \mathbf{w}_h = \mathbf{0} \text{ on } \partial\Omega, %\\ &
~\mathbf{w}_h = \mathbf{0}\text{ on } \partial\mathcal{B}_h \Big\}.
%\end{split}
\]%\end{equation}
Further, for each time $t$, let the extension operator $\mathcal{E}$ be the simplest and most popular one: If $\eta^U_{\partial}$ is the subset of $\eta^U$ containing the indices of velocity nodes in $\partial\mathcal{B}_h$ and $f$ is a continuous (piecewise $P_m$) function defined on $\partial\mathcal{B}_h$, 
\[%\begin{equation}
    \mathcal{E}f=\sum_{i\in\eta^U_\partial} f(\mathbf{x}^i)\,\mathcal{N}^i(\mathbf{x}).
\]%\end{equation}
In other words, $f$ is extended to $\Omega_{fh}$ by setting all nodal values not belonging to $\partial\mathcal{B}_h$ to zero and interpolating according to the adopted finite element space.

Since $U_h(\mathbf{q})$ restricted to $\partial\mathcal{B}_h$ contains at least $P_1$ polynomials, $\mathbf{w}_h = \mathbf{H}({\mathbf{q}})\mathbf{d}$ is satisfied exactly for all $\mathbf{d}$. Up to the geometrical difference between $\partial\mathcal{B}$ and $\partial\mathcal{B}_h$, which is out of the scope of this contribution, it thus holds that
\begin{eqnarray}
    W_h\paren{\mathbf{q}} & \subset & W\paren{\mathbf{q}},\\
    W_{0h}\paren{\mathbf{q}} & \subset & W_0\paren{\mathbf{q}},\\
    W_h\paren{\mathbf{q}} & = & W_{0h}\paren{\mathbf{q}}\oplus V\paren{\mathbf{q}}. \label{eqwhc}
\end{eqnarray}

Let $\mathbf{u}_{sh}$ be the interpolant of $\mathbf{u}_{s}$ in $U_h$ (restricted to $\partial\mathcal{B}_h$). The configuration manifold $Q=\mathbb{R}^d\times \SO\paren{d}$ is kept exact, but of course in the semidiscrete problem one computes approximations of the exact functions $\mathbf{q}\paren{t}=\paren{\mathbf{x}_c(t),\mathbf{Q}(t)}:[0,T]\to Q$ and $\mathbf{s}(t)=\paren{\mathbf{v}_c(t),\bomega(t)}:[0,T]\to \mathbb{R}^{n_c}$. We  add the subscript $h$ to these functions to make this fact explicit.

The approximate velocity $\mathbf{u}_h(\cdot,t)$ is sought belonging to $W_h\paren{\mathbf{q}_h(t)}$ and satisfying \eqref{eqtypeIc}. Thus, from \eqref{eqwhc}, it can be decomposed as
\[%\begin{equation}
    \mathbf{u}_h = \widetilde{\mathbf{H}}\,\mathbf{s}_h+\mathbf{u}_{0h}+\mathcal{E}\mathbf{u}_{sh},
\]%\end{equation}
where $\mathbf{u}_{0h}\in W_{0h}(\mathbf{q}_h(t))$. Let
$\mathbf{H}^j(\mathbf{q})=\widetilde{\mathbf{H}}(\mathbf{q},\mathbf{x}^j)$ (i.e.,
$\mathbf{H}^j(\mathbf{q})=\mathbf{H}(\mathbf{q},\mathbf{x}^j)$ if $j\in\eta^U_\partial$, and $=\mathbf{0}$, the null $d\times n_c$ matrix, otherwise). Then the nodal values $\mathbf{u}^j(t)$ of $\mathbf{u}_h(\cdot,t)$ are unconstrained unknowns if $j\in \eta^U_0$ (interior nodes, i.e., $\eta^U_0=\eta^U\setminus\eta^U_\partial$) and, for $j\in\eta^U_\partial$, they obey
\[%\begin{equation}
    \mathbf{u}^j(t)=\mathbf{H}^j(\mathbf{q}_h(t)) \mathbf{s}_h(t) + \mathbf{u}^j_{s}(t),
\]%\end{equation}
where $\mathbf{u}^j_s=\mathbf{u}_s(\mathbf{x}^j(t),t)$.

The semidiscrete Galerkin formulation for a type-I squirmer in a Newtonian fluid reads: 
{\em 
Determine functions $\mathbf{q}_h\paren{t}=\paren{\mathbf{x}_{ch}(t),\mathbf{Q}_h(t)}:[0,T]\to Q$,
$\mathbf{s}_h\paren{t}=\paren{\mathbf{v}_{ch}(t),\bomega_h(t)}:[0,T]\to \mathbb{R}^{n_c}$, $\mathbf{u}_h(\cdot,t)\,\in\,W_h(\mathbf{q}_h(t))$ and $p_h(\cdot,t)\,\in\,M_h(\mathbf{q}_h(t))$ such that
\begin{equation}
\mathbf{u}_h(\cdot,t)-\widetilde{\mathbf{H}}(\cdot,t)\mathbf{s}_h(t)-\mathcal{E}\mathbf{u}_{sh}(\cdot,t)\,\in\,W_{0h}(\mathbf{q}_h(t)) \label{eqtypeIz2}
\end{equation}
and
\begin{eqnarray}
\frac{d\mathbf{x}_{ch}}{dt}-\mathbf{v}_{ch}&=&\mathbf{0}, \label{eqtypeIa2}\\
\frac{d\mathbf{Q}_h}{dt}-\skk\cparen{\bomega_h}\,\mathbf{Q}_h&=&\mathbf{0},\label{eqtypeIb2}\\
    \int_{\Omega_{fh}} \rho\,\frac{D\mathbf{u}_h}{Dt}\cdot\mathbf{w}_h\, d\mathbf{x} + \int_{\Omega_{fh}}2\mu\nabla^S\mathbf{u}_{h}:\nabla^s\mathbf{w}_h\, d\mathbf{x} & -&\nonumber\\ -\int_{\Omega{fh}} p_h\,\nabla\cdot\mathbf{w}_h\,d\mathbf{x}&=&0,\label{eqmom12}\\
\int_{\Omega_{fh}} z_h\nabla \cdot \mathbf{u}_h\, d\mathbf{x} &=& 0,\label{eqmass12}
\end{eqnarray}
for all $\mathbf{w}_h\in W_h(\mathbf{q}_h(t))$, for all $z_h\in M_h(\mathbf{q}_h(t))$, for all $t$.}

Recalling that, if $\{\mathbf{e}^i\}$, $i=1,\ldots,d$ is the canonical basis of $\mathbb{R}^d$, a basis for $U_h(\mathbf{q})$ is provided by
\[
\Big \{ \mathcal{N}^j\mathbf{e}^i,~j\in\eta^U, i=1,\ldots,d\Big\},
\]
and so a basis for $W_h(\mathbf{q})$ is provided by
\begin{equation}
   \Big \{ \mathcal{N}^j\mathbf{e}^i,~j\in\eta^U_0, i=1,\ldots,d\Big\}\cup \Big\{\mathbf{h}^k(\mathbf{q},\cdot),~k=1,\ldots,n_c\Big\}.
   \label{eqbasiswh}
\end{equation}
Because of the extension $\mathcal{E}$ chosen, the nodal values of $\mathbf{h}^k(\mathbf{q},\cdot)$ are equal to the $k$-th column of $\mathbf{H}^j(\mathbf{q})$ if $j\in\eta^U_\partial$, and $\mathbf{0}$ otherwise.
%%%%%%%%%%%%%%%%%%%%%%%%%%%%%%%%%%%%%%%%%%%%%%%%%%
\subsection{Matrix formulation for type-I squirmers}
For the sake of simplicity, we will present the matrix problem of the Galerkin formulation for the linear case ($\rho=0$). The extension to the case $\rho>0$, or to strain-rate dependent material viscosity, should be straightforward for FEM practitioners.

Let us collect the nodal velocity unknowns, vertically, into the time dependent column vector
\begin{eqnarray}
\underline{U}(t)& =& \cparen{\mathbf{u}^{i}(t)} \qquad\forall i\in\eta^{U}, \nonumber
\end{eqnarray}
and, similarly, the pressure unknowns as $\underline{P}(t) = \cparen{p^{i}(t)}$, $\forall i\in\eta^{P}$.

{\bf Assume for the moment that}, for a given configuration $\mathbf{q}$ of the system and a given instant $t$, {\bf the boundary $\partial\mathcal{B}_h$ is simply a force-free boundary}. Then there is no doubt as to how to proceed: The velocity space is the whole of $U_h(\mathbf{q})$ and standard finite element treatment of equations \eqref{eqmom12}-\eqref{eqmass12} lead to the algebraic system
\begin{eqnarray}
\mathbb{A}\paren{\mathbf{q}}\underline{U} + \mathbb{G}\paren{\mathbf{q}}\underline{P} 
&=& \underline{F},\label{eqmatrix1a}\\
\mathbb{D}\paren{\mathbf{q}}\underline{U}+\mathbb{E}\paren{\mathbf{q}}\underline{P} &=& \underline{G},\label{eqmatrix1b}
\end{eqnarray}
where $\mathbb{A} = \cparen{\mathbf{A}_{ij}}$, $\mathbb{G} = \cparen{\mathbf{G}_{ik}}$ and $\mathbb{D} = \cparen{\mathbf{D}_{ki}}$, with $i,j\in\eta^U$ and $k\in\eta^P$ 
are composed of the block matrices
\[%\begin{equation}
\mathbf{A}_{ij} = \int_{\Omega_{fh}(\mathbf{q})} \mu\paren{\nabla\mathcal{N}^{i}\cdot\nabla\mathcal{N}^{j}\mathbf{I}_{d} + \nabla\mathcal{N}^{j}\otimes\nabla\mathcal{N}^{i}}\, d\mathbf{x}, %\label{eqn:bl_mat_rn_A}
\]%\end{equation}
\[%\begin{equation}
\mathbf{G}_{ik} = -\int_{\Omega_{fh}(\mathbf{q})} \mathcal{M}^{k}\nabla\mathcal{N}^{i}\, d\mathbf{x}, \quad \mathbf{D}_{ki} = -\mathbf{G}_{ik}, %\label{eqn:bl_mat_rn_G}
\]%\end{equation}
\[%\begin{equation}
    \mathbb{E}=\mathbf{0},\quad \underline{F}=\mathbf{0},\text{  and } \underline{G}=\mathbf{0}.
\]%\end{equation}
Equations \eqref{eqmatrix1a}-\eqref{eqmatrix1b} form the classical Stokes matrix system that arises from the Galerkin formulation. They are algebraic materializations of the momentum equation \eqref{eqmom12} and the incompressibility equation \eqref{eqmass12}.

\bigskip

\noindent{\bf Important remark:} In what follows, the specific steps that lead to the algebraic system \eqref{eqmatrix1a}-\eqref{eqmatrix1b} are largely irrelevant. Pressure-stabili\-za\-tion schemes, for example, lead to a different (non-zero) matrix $\mathbb{E}$. The presence of volumetric forces in the liquid modifies $\underline{F}$ and possibly $\underline{G}$. {\bf The reader may take \eqref{eqmatrix1a}-\eqref{eqmatrix1b} as the algebraic system arising from her/his favorite finite element or finite volume solver}. Up to now, the case is rather dull, just a fluid domain which has some force-free holes in it. The procedures below show how to manipulate the system so as to turn those holes into interesting type-I squirmers.

\bigskip
Let $n_0=\card (\eta^U_0)$ be the number of interior velocity nodes and $n_\partial=\card (\eta^U_\partial )$ be the number of boundary nodes on $\partial\mathcal{B}_h$, so that $n_U=\card (\eta^U) = n_0+n_\partial$, also let $n_P = \card(\eta^P)$ the number of pressure nodes. The velocity unknowns are partitioned into the column arrays $\underline{U_0}$ (of dimension $n_0d$) and $\underline{U_\partial}$ (of dimension $n_\partial d$) as
\begin{eqnarray*}
\underline{U_0}(t)& =& \cparen{\mathbf{u}^{i}(t)} \qquad\forall i\in\eta^{U}_0, \\
\underline{U_\partial}(t)& =& \cparen{\mathbf{u}^{i}(t)} \qquad\forall i\in\eta^{U}_\partial, 
\end{eqnarray*}
 so that, though the index sets do not need to be consecutive, we can without loss of generality write
 \begin{equation}
\underline{U}(t)=\left [ \begin{array}{c} \underline{U_0}(t)\\\underline{U_\partial}(t)\end{array}\right ].
\label{eqdecompu}
 \end{equation}

A crucial role is played by the $(n_Ud)\times n_c$ block matrix $\mathbb{H}(\mathbf{q})$, obtained by arranging the matrices $\mathbf{H}^j(\mathbf{q})=\widetilde{\mathbf{H}}(\mathbf{q},\mathbf{x}^j)$, $j\in \eta^U$, in a column:
\[%\begin{equation}
    \mathbb{H}(\mathbf{q})=\cparen{\mathbf{H}^j(\mathbf{q})},\quad j\in \eta^U.
\]%\end{equation}
Following the decomposition \eqref{eqdecompu} of the velocity unknowns, the matrix $\mathbb{H}$  decomposes into a submatrix $\mathbb{H}^0$ consisting of the blocks $\mathbf{H}^j$ with $j\in\eta^U_0$, {\em which is identically zero}, and the submatrix $\mathbb{H}^\partial$, corresponding to $j\in\eta^U_\partial$. Acting analogously on $\mathbb{A}$, $\mathbb{G}$ and $\underline{F}$ (vertical partitioning)
we have
\[%\begin{equation}
    \mathbb{H}=\left [ \begin{array}{c} \mathbf{0} \\
    \mathbb{H}^{\partial}
    \end{array} \right ]~, \quad
    \mathbb{A}=\left [ \begin{array}{c} \mathbb{A}^{0} \\
    \mathbb{A}^{\partial}
    \end{array} \right ],\quad 
    \mathbb{G}=\left [ \begin{array}{c} \mathbb{G}^{0} \\
    \mathbb{G}^{\partial}
    \end{array} \right ]~,\quad 
    \underline{F}=\left [ \begin{array}{c} \underline{F^0} \\
    \underline{F^\partial}
    \end{array} \right ]~. 
\]%\end{equation}
Up to now the velocity lines of the block matrices of the Stokes system have been classified according to whether they correspond to nodes in the interior or on the squirmer's boundary. No operation has been performed. The list of required operations is:
\begin{enumerate}
    \item Take the lines related to nodes in $\eta^U_\partial$ out of $\mathbb{A}$ and replace them with lines of the identity matrix to obtain $\widehat{\mathbb{A}}$:
    \[%\begin{equation}
    \widehat{\mathbb{A}}=\left [ \begin{array}{c} \mathbb{A}^{0} \\
    \mathbb{I}^\partial
    \end{array} \right ].
    \]%\end{equation}
    The block $\mathbb{I}^\partial$ has all elements in each line equal to zero, except for the diagonal, which is equal to one.
    \item Take the lines related to nodes in $\eta^U_\partial$ out of $\mathbb{G}$ and replace them with the null matrix to obtain $\widehat{\mathbb{G}}$:
    \[%\begin{equation}
    \widehat{\mathbb{G}}=\left [ \begin{array}{c} \mathbb{G}^{0} \\ \mathbf{0}
    \end{array} \right ].
    \]%\end{equation}
    \item Pre-multiply the lines taken out in the two previous actions by the transpose of $\mathbb{H}^\partial$ to obtain the matrices $\mathbb{S}^\partial$ and $\mathbb{T}^\partial$:
    \[%\begin{equation}
    \mathbb{S}^{\partial}=\left ( \mathbb{H}^\partial \right )^T \mathbb{A}^{\partial},\quad
    \mathbb{T}^{\partial}=\left ( \mathbb{H}^\partial \right )^T \mathbb{G}^{\partial}.
    \]%\end{equation}
    \item Denoting by $\underline{U_s}=\cparen{\mathbf{u}_s^i}$, $i\in\eta^U_\partial$, the $(n^U_\partial d)\times 1$ column array with nodal values of the slip velocity $\mathbf{u}_{sh}$, 
    build $\underline{\widehat{F}}(\mathbf{u}_{sh})$ and $\underline{B^\partial}$ as
    \[%\begin{equation}
        \underline{\widehat{F}}(\mathbf{u}_{sh})=\left [ \begin{array}{c}
    \underline{F^0} \\ \underline{U_s} \end{array} \right ],\quad
    \underline{B^\partial}=\left ( \mathbb{H}^\partial \right )^T \underline{F^\partial}.
    \]%\end{equation}
\end{enumerate}

Once these matrices are built, a task that the experienced finite element coder easily figures out how to do {\em at the element level} (before the assembly operation), the semi-discrete finite element formulation becomes, in matrix form: 

{\em Determine functions $\mathbf{q}_h\paren{t}=\paren{\mathbf{x}_{ch}(t),\mathbf{Q}_h(t)}:[0,T]\to Q$,
$\mathbf{s}_h\paren{t}=\paren{\mathbf{v}_{ch}(t),\bomega_h(t)}:[0,T]\to \mathbb{R}^{n_c}$, $\underline{U}(t):[0,T]\to\mathbb{R}^{n_Ud}$ and $\underline{P}(t):[0,T]\to\mathbb{R}^{n_P}$ such that, for each $t$,}
\begin{eqnarray}
\frac{d\mathbf{x}_{ch}}{dt}-\mathbf{v}_{ch}&=&\mathbf{0}, \label{eqtypeIa3}\\
\frac{d\mathbf{Q}_h}{dt}-\skk\cparen{\bomega_h}\,\mathbf{Q}_h&=&\mathbf{0},\label{eqtypeIb3}\\
\widehat{\mathbb{A}}~\underline{U}+\widehat{\mathbb{G}}~\underline{P}-\mathbb{H}~\mathbf{s}_h&=&\underline{\widehat{F}}(\mathbf{u}_{sh}),\label{eqtypeIc3}\\
\mathbb{S}^\partial~\underline{U}+\mathbb{T}^\partial~ \underline{P} &=& \underline{B^\partial}, \label{eqtypeId3}\\
\mathbb{D}~ \underline{U} + \mathbb{E}~ \underline{P} &=& \underline{G}. \label{eqtypeIe3}
\end{eqnarray}

Several comments are in order:
\begin{itemize}
    \item In the Stokes case above, equations \eqref{eqtypeIc3}-\eqref{eqtypeIe3} can be solved to produce a function $Q\times U_h|_{\partial\mathcal{B}_h}\to \mathbb{R}^{n_c}$ giving (abusing the notation)
    \[%\begin{equation}
        \mathbf{s}_h = \mathbf{s}_h\paren{\mathbf{q}_h,\mathbf{u}_{sh}}.
    \]%\end{equation}
    The dependence on $\mathbf{q}_h$ arises because, though not made explicit, all matrices depend on the geometry of $\Omega_{fh}$ and thus on $\mathbf{q}_h$. The dependence on $\mathbf{u}_{sh}$ is {\em linear}, in the linear case ($\rho=0$, constant $\mu$). To prove that for a given $\mathbf{q}_h$ the mapping $\mathbf{u}_{sh}\mapsto \mathbf{s}_h$ is well defined, it then suffices to show that if $\mathbf{u}_{sh}=0$ then $\mathbf{s}_h$ is necessarily zero (along with $\underline{U}$ and $\underline{P}$). In the continuous case this is immediate from \eqref{eqenergy} and Korn's inequality. In the discrete case the argument is analogous, but since $\nabla\cdot {\bf u}_h$ is not automatically zero one has to rely on the stability of the discrete pressure-velocity coupling (possibly stabilized). Since \eqref{eqtypeIa3}-\eqref{eqtypeIb3} can be rewritten as $\frac{d\mathbf{q}_h}{dt}=\mathbf{g}(\mathbf{q}_h,\mathbf{s}_h)$, the whole problem turns into the ODE
    \begin{equation}
        \frac{d\mathbf{q}_h}{dt}=\mathbf{g}\paren{\mathbf{q}_h,\mathbf{s}_h(\mathbf{q}_h,\mathbf{u}_{sh})}
        \label{eqode}
    \end{equation}
    to be solved on the manifold $Q=\mathbb{R}^d\times\SO(d)$, to which the vector field $\mathbf{g}$ is tangent. 
    \item A unique local (for $t<T$ small enough) solution of \eqref{eqode} starting at some $\mathbf{q}_h(0)$ for which the mesh is good enough can be shown to exist, since $\mathbf{g}$ is indeed Lipschitz. When trying to make $T\to +\infty$ to prove a global result in time, two kinds of difficulties appear. The most immediate one is purely numerical. The mesh may become distorted turning $\mathbf{g}$ singular. This can be overcome by a suitable remeshing algorithm. A more profound problem however persists, which is also present in the exact problem. As the squirmers evolve over the domain they may head towards the walls, or one towards the other. This may make the tangent force $\mathbf{f}_s$ to grow without bound, taking the dissipation to infinity ($\mathbf{u}\not\in H^1(\Omega_f)^d$) and thus not just making $\mathbf{g}$ singular but also making the model unrealistic (no squirmer can spend infinite power).
    \item The previous comments made use of the correspondence between the variational problem \eqref{eqtypeIz2}-\eqref{eqmass12} and its matrix formulation \eqref{eqtypeIa3}-\eqref{eqtypeIe3}. Let us make it explicit. The first $n_0d$ lines of \eqref{eqtypeIc3} express $\mathbb{A}^{0}\underline{U}+\mathbb{G}^0\underline{P}=0$, which enforces \eqref{eqmom12} for all $\mathbf{w}_h\in W_{0h}(\mathbf{q}_h)$. The last $n_\partial d$ lines express $\underline{U_\partial}-\mathbb{H}^\partial \mathbf{s}_h=\underline{U_s}$, which enforces \eqref{eqtypeIz2}. Because of \eqref{eqbasiswh}, equation \eqref{eqtypeId3} enforces \eqref{eqmom12} for all $\mathbf{w}_h\in V(\mathbf{q}_h)$. Finally, \eqref{eqtypeIe3}, which was left untouched, enforces \eqref{eqmass12}.
    \item In the case of vertex-centered finite volumes the correspondence must be made with integral versions of the differential problem \eqref{eqtypeIa}-\eqref{eqtypeIg}. The first $n_0d$ lines of \eqref{eqtypeIc3} enforce the momentum conservation equation \eqref{eqtypeId} at interior volumes and are left untouched. The last $n_\partial d$ lines enforce \eqref{eqtypeIc}, and \eqref{eqtypeId3} enforces the force-free and torque-free constraints \eqref{eqtypeIf}-\eqref{eqtypeIg}. Finally, \eqref{eqmass12} enforces the incompressibility equation \eqref{eqtypeIe} and was left untouched.
    \item If the fluid's inertia is considered, besides the matrix $\mathbb{A}$ being modified (and possibly some others too), a term $\mathbb{M}\frac{d\underline{U}}{dt}$ will appear in \eqref{eqmatrix1a}. One has then to operate in this matrix as follows: Decomposing $\mathbb{M}$ as
    \[%\begin{equation}
        \mathbb{M} = \left [ \begin{array}{c} \mathbb{M}^0 \\
        \mathbb{M}^\partial \end{array} \right ],
    \]%\end{equation}
    one builds $\widehat{\mathbb{M}}$ and $\mathbb{L}^\partial$ following
    \[%\begin{equation}
        \widehat{\mathbb{M}}=\left [ \begin{array}{c} \mathbb{M}^0 \\
        \mathbf{0} \end{array} \right ],\quad
        \mathbb{L}^\partial = \left ( \mathbb{H}^\partial \right )^T\mathbb{M}^\partial.
    \]%\end{equation}
    Finally, one adds the term $\widehat{\mathbb{M}}\frac{d\underline{U}}{dt}$ to the left-hand side of \eqref{eqtypeIc3} and the term $\mathbb{L}^\partial\frac{d\underline{U}}{dt}$ to the left-hand side of \eqref{eqtypeId3} and the matrix formulation now considers the fluid's inertia.
    \item The incorporation of variable viscosity, strain-rate-dependent for example, does not require any change in the formulation or in the manipulation of the matrices. The only consequence is that $\widetilde{\mathbb{A}}$ and $\mathbb{S}^\partial$ will depend on $\underline{U}$. 
    \item The extension to many ($N$) squirmers is straightforward. There will be one set of equations \eqref{eqtypeIa3}-\eqref{eqtypeIb3} per squirmer, of course, and the vector $\mathbf{s}_h=\paren{\mathbf{v}_c^{(1)},\bomega^{(1)},\ldots,\mathbf{v}_c^{(N)},\bomega^{(N)}}$ will have $Nn_c$ unknowns. Nevertheless, the operations (1)-(4) above can be performed sequentially squirmer by squirmer because no two squirmers share the same boundary node. Each squirmer adds $n_c$ columns to the matrix $\mathbb{H}$, modifies blocks of lines of $\mathbb{A}$, $\mathbb{G}$, $\underline{F}$ (and $\mathbb{M}$) and adds $n_c$ lines to $\mathbb{S}^\partial$, $\mathbb{T}^\partial$, $\underline{B^\partial}$ (and $\mathbb{L}^\partial$). 
\end{itemize}
%%%%%%%%%%%%%%%%%%%%%%%%%%%%%%%%%%%%%%%%%%%%%%%%%%
\subsection{Semidiscrete Galerkin formulation for type-II squirmers}

Most of the notation introduced in the discretization of type-I squirmers is also useful for type-II squirmers and will be used in what follows. Care was taken in not having the same symbol denoting something for type-I squirmers and something different for type-II ones. If the symbol is the same, it {\em is} the same entity, with the same definition.

The discrete version of the spaces, in this case, is obtained enforcing the no-penetration condition pointwise at each node $j\in\eta^U_\partial$, for which a unit normal vector $\mathbf{n}^j$ is assumed given. The spaces are
\[%\begin{equation}
   %\begin{split}
\widetilde{W}_h\paren{\mathbf{q}} = \Big\{ \mathbf{w}_h\in U_h\paren{\mathbf{q}}\, : \mathbf{w}_h = \mathbf{0} \text{ on } \partial\Omega,~\mathbf{n}^j\cdot\mathbf{w}_h(\mathbf{x}^j) = \mathbf{n}^j\cdot \mathbf{H}({\mathbf{q}})\mathbf{d} \text{ for } j\in\eta^U_\partial, ~\mathbf{d}\,\in\,\mathbb{R}^{n_c} \Big\}
%\end{split}
\]%\end{equation}
and
\[%\begin{equation}
%   \begin{split}
W_{\parallel h}\paren{\mathbf{q}} = \Big\{ \mathbf{w}_h\in U_h\paren{\mathbf{q}}\, : \mathbf{w}_h = \mathbf{0} \text{ on } \partial\Omega, %\\ &
~\mathbf{n}^j\cdot\mathbf{w}_h(\mathbf{x}^j) = 0 \text{ for } j\in\eta^U_\partial \Big\}.
%\end{split}
\]%\end{equation}

Though the inclusion $\widetilde{W}_h\paren{\mathbf{q}}\subset \widetilde{W}\paren{\mathbf{q}}$ is not valid (it only holds approximately),  it holds that
\begin{eqnarray}
    \widetilde{W}_h\paren{\mathbf{q}} & = & {W}_{\parallel h} \paren{\mathbf{q}}\oplus V\paren{\mathbf{q}}. \label{eqwhcII}
\end{eqnarray}

The approximate velocity $\mathbf{u}_h(\cdot,t)$ is sought belonging to $\widetilde{W}_h\paren{\mathbf{q}_h(t)}$ and satisfying \eqref{eqtypeIIc} {\em pointwise} at all nodes of $\partial\mathcal{B}_h$. Thus, from \eqref{eqwhcII}, it can be decomposed as
\[%\begin{equation}
    \mathbf{u}_h = \widetilde{\mathbf{H}}\,\mathbf{s}_h+\mathbf{u}_{\parallel h},
\]%\end{equation}
where $\mathbf{u}_{\parallel h}\in W_{\parallel h}(\mathbf{q}_h(t))$. In other words, the unconstrained unknowns are the nodal values $\mathbf{u}^j(t)$ of $\mathbf{u}_h(\cdot,t)$ if $j\in \eta^U_0$ {\em and their tangential component} $\bPi \mathbf{u}^j(t)$ if $j\in\eta^U_\partial$ (recall that $\bPi=\mathbf{I}_d-\mathbf{n}\mathbf{n}^T$),
whereas the normal component of $\mathbf{u}^j(t)$, $j\in\eta^U_\partial$, obeys
\begin{equation}
    \mathbf{n}^j\cdot\mathbf{u}^j(t)=\mathbf{n}^j\cdot\mathbf{H}^j(\mathbf{q}_h(t)) \mathbf{s}_h(t).
    \label{eqtypeIIkin}
\end{equation}

The semidiscrete Galerkin formulation for a type-II squirmer in a Newtonian fluid reads: 
{\em 
Determine functions $\mathbf{q}_h\paren{t}=\paren{\mathbf{x}_{ch}(t),\mathbf{Q}_h(t)}:[0,T]\to Q$,
$\mathbf{s}_h\paren{t}=\paren{\mathbf{v}_{ch}(t),\bomega_h(t)}:[0,T]\to \mathbb{R}^{n_c}$, $\mathbf{u}_h(\cdot,t)\,\in\,\widetilde{W}_h(\mathbf{q}_h(t))$ and $p_h(\cdot,t)\,\in\,M_h(\mathbf{q}_h(t))$ such that
\begin{equation}
\mathbf{u}_h(\cdot,t)-\widetilde{\mathbf{H}}(\cdot,t)\mathbf{s}_h(t)\,\in\,W_{\parallel h}(\mathbf{q}_h(t)) \label{eqtypeIIz2}
\end{equation}
and
\begin{eqnarray}
\frac{d\mathbf{x}_{ch}}{dt}-\mathbf{v}_{ch}&=&\mathbf{0}, \label{eqtypeIIa2}\\
\frac{d\mathbf{Q}_h}{dt}-\skk\cparen{\bomega_h}\,\mathbf{Q}_h&=&\mathbf{0},\label{eqtypeIIb2}\\
    \int_{\Omega_{fh}} \rho\,\frac{D\mathbf{u}_h}{Dt}\cdot\mathbf{w}_h\, d\mathbf{x} + \int_{\Omega_{fh}}2\mu\nabla^S\mathbf{u}_{h}:\nabla^s\mathbf{w}_h\, d\mathbf{x} & -&\nonumber\\ -\int_{\Omega{fh}} p_h\,\nabla\cdot\mathbf{w}_h\,d\mathbf{x}&=&\int_{\partial\mathcal{B}_h}\mathbf{f}_s\cdot \mathbf{w}_{\parallel h},\label{eqmom22}\\
\int_{\Omega_{fh}} z_h\nabla \cdot \mathbf{u}_h\, d\mathbf{x} &=& 0,\label{eqmass22}
\end{eqnarray}
for all $\mathbf{w}_h\in \widetilde{W}_h(\mathbf{q}_h(t))$, for all $z_h\in M_h(\mathbf{q}_h(t))$, for all $t$.}

Type-II squirmers have more velocity degrees of freedom than type-I squirmers because 
the tangential components of the velocity are unknowns, i.e.,
\begin{equation}
    W_{\parallel h}(\mathbf{q})= W_{0h}(\mathbf{q})\oplus T_h(\mathbf{q}), \label{eqw0plusth}
\end{equation}
where
\[%\begin{equation}
    T_h(\mathbf{q}) = \Big \{ \mathbf{w}_h\in W_{\parallel h}(\mathbf{q}),~\mathbf{w}_h(\mathbf{x}^j)=\mathbf{0}~\forall j\in\eta^U_0 \Big \}.
\]%\end{equation}

Denoting by $\{\boldsymbol{\tau}^i\}$, $i=1,\ldots,d-1$ unit tangent vectors that, added to $\mathbf{n}$, form a basis of $\mathbb{R}^d$, a basis for $T_h(\mathbf{q})$ is provided by
\[
\Big \{ \mathcal{N}^j\mathbf{\boldsymbol{\tau}}^i,~j\in\eta^U_\partial, i=1,\ldots,d-1\Big\}.
\]
These vector fields must be added to the basis of $W_h(\mathbf{q})$ given in \eqref{eqbasiswh} to obtain a basis of $\widetilde{W}_h(\mathbf{q})$.
%%%%%%%%%%%%%%%%%%%%%%%%%%%%%%%%%%%%%%%%%%%%%%%%%%
\subsection{Matrix formulation for type-II squirmers}
We again present the matrix problem of the Galerkin formulation for the linear case ($\rho=0$), and the point of departure is the system \eqref{eqmatrix1a}-\eqref{eqmatrix1b} that arises when the holes are considered as imposed-force boundaries, only that this time the force imposed is not zero but $\mathbf{f}_s$. For this reason, to the array $\underline{F}^\partial$ (arising from volumetric forces, for example) will be added another array $\underline{R^\partial} \paren{\mathbf{f}_s}$ that is the contribution of $\mathbf{f}_s$ on the boundary nodes. This is still totally standard. We now show how to turn the holes into type-II squirmers.

One needs to build, for each node $j\in\eta^U_\partial$, the $d\times d$ projection matrices
\[%\begin{equation}
    \mathbf{P}_{\boldsymbol{\tau}}^j=\mathbf{I}_d-\mathbf{n}^j\left ( \mathbf{n}^j \right )^T,\quad
    \mathbf{P}_{\mathbf{n}}^j=\mathbf{n}^j\left ( \mathbf{n}^j \right )^T,
\]%\end{equation}
and collect them into two {\em block-diagonal} matrices $\mathbb{P}_{\boldsymbol{\tau}}$ and $\mathbb{P}_{\mathbf{n}\alpha}$, with entries given by
\[%\begin{equation}
    \left (\mathbb{P}_{\boldsymbol{\tau}} \right )_{jj} = \mathbf{P}_{\boldsymbol{\tau}}^j,\quad
    \left (\mathbb{P}_{\mathbf{n}\alpha} \right )_{jj}=\alpha_j\mathbf{P}_{\mathbf{n}}^j.
\]%\end{equation}
Above, no summation in $j$ is implied, $j$ runs over the index set $\eta^U_\partial$, and a set of positive numbers $\alpha_j$ has been incorporated that is useful to avoid ill-conditioning (though all cases shown here have $\alpha_j=1$). The dimensions of these block-diagonal matrices is $n_\partial d\times n_\partial d$.

The list of required operations is:
\begin{enumerate}
    \item Modify the lines related to nodes in $\eta^U_\partial$ of $\mathbb{H}$, $\mathbb{A}$ and $\mathbb{G}$ to obtain $\widetilde{\mathbb{H}}$, $\widetilde{\mathbb{A}}$ and $\widetilde{\mathbb{G}}$:
    \[%\begin{equation}
    \widetilde{\mathbb{H}}=\left [ \begin{array}{c} \mathbf{0}\\
    \mathbb{P}_{\mathbf{n}\alpha} \mathbb{H}^\partial 
    \end{array} \right ],\quad
    \widetilde{\mathbb{A}}=\left [ \begin{array}{c} \mathbb{A}^{0}\\
    \mathbb{P}_{\boldsymbol{\tau}} \mathbb{A}^\partial + \mathbb{P}_{\mathbf{n}\alpha}
    \end{array} \right ],\quad
    \widetilde{\mathbb{G}}=\left [ \begin{array}{c} \mathbb{G}^{0}\\
    \mathbb{P}_{\boldsymbol{\tau}} \mathbb{G}^\partial 
    \end{array} \right ].
    \]%\end{equation}
    \item Compute the matrices ${\mathbb{S}}^{\partial}$ and ${\mathbb{T}}^\partial$:
    \[%\begin{equation}
    {\mathbb{S}}^{\partial}=\left ( \mathbb{H}^\partial \right )^T  \mathbb{A}^{\partial},\quad
    {\mathbb{T}}^{\partial}=\left ( \mathbb{H}^\partial \right )^T  \mathbb{G}^{\partial}.
    \]%\end{equation}
    \item Build $\underline{\widetilde{F}}(\mathbf{f}_{s})$ and $\underline{{B}^\partial}$ as
    \[%\begin{equation}
        \underline{\widetilde{F}}(\mathbf{f}_{s})=\left [ \begin{array}{c}
    \underline{F^0} \\  
    \mathbb{P}_{\boldsymbol{\tau}} \left ( \underline{F^\partial}+ \underline{R^\partial}(\mathbf{f}_s) \right )
    \end{array} \right ],\quad
    \underline{{B}^\partial}=\left ( \mathbb{H}^\partial \right )^T    \underline{F^\partial}.
    \]%\end{equation}
    Notice that $\underline{B^\partial}$ is the same as in the case of type-I squirmers. The contribution of the tangential force $\mathbf{f}_s$ does not intervene in its calculation.
\end{enumerate}

As before, all these operations can be performed {\em at the element level} (before the assembly operation). The semi-discrete formulation in matrix form reads: 

{\em Determine functions $\mathbf{q}_h\paren{t}=\paren{\mathbf{x}_{ch}(t),\mathbf{Q}_h(t)}:[0,T]\to Q$,
$\mathbf{s}_h\paren{t}=\paren{\mathbf{v}_{ch}(t),\bomega_h(t)}:[0,T]\to \mathbb{R}^{n_c}$, $\underline{U}(t):[0,T]\to\mathbb{R}^{n_Ud}$ and $\underline{P}(t):[0,T]\to\mathbb{R}^{n_P}$ such that, for each $t$,}
\begin{eqnarray}
\frac{d\mathbf{x}_{ch}}{dt}-\mathbf{v}_{ch}&=&\mathbf{0}, \label{eqtypeIIa3}\\
\frac{d\mathbf{Q}_h}{dt}-\skk\cparen{\bomega_h}\,\mathbf{Q}_h&=&\mathbf{0},\label{eqtypeIIb3}\\
\widetilde{\mathbb{A}}~\underline{U}+\widetilde{\mathbb{G}}~\underline{P}-\widetilde{\mathbb{H}}~\mathbf{s}_h&=&\underline{\widetilde{F}}(\mathbf{f}_{s}),\label{eqtypeIIc3}\\
{\mathbb{S}}^\partial~\underline{U}+{\mathbb{T}}^\partial~ \underline{P} &=& \underline{{B}}^\partial, \label{eqtypeIId3}\\
\mathbb{D}~ \underline{U} + \mathbb{E}~ \underline{P} &=& \underline{G}. \label{eqtypeIIe3}
\end{eqnarray}

Most of the remarks made for type-I squirmers also hold for type-II ones, but some differences deserve additional comments:
\begin{itemize}
    \item Following \eqref{eqw0plusth}, the weak momentum equation needs to hold $\forall\mathbf{w}_h$ in $W_{0h}(\mathbf{q}_h)$, $\forall\mathbf{w}_h$ in $T_h(\mathbf{q}_h)$ and $\forall\mathbf{w}_h$ in $V(\mathbf{q}_h)$. The first $n_0 d$ lines of \eqref{eqtypeIIc3}, as before, enforce the first of these conditions. For type-II squirmers, the $n_\partial d$ lines of \eqref{eqtypeIIc3} enforce, {\em simultaneously}, the kinematical constraint \eqref{eqtypeIIkin} (equivalent to \eqref{eqtypeIIz2}) in the {\em normal} velocity components, together with \eqref{eqmom22} for the {\em tangential} component, more specifically for all $\mathbf{w}_h\in T_h(\mathbf{q}_h)$. The introduction of the projection matrices serves this purpose, borrowing from previous works on slip boundary conditions for fluid flow \cite{engelman1982implementation,Montefuscolo:2012ts}. In a nutshell, the equations expressing conservation of momentum are rotated to the normal-tangential frame, so as to then keep just the tangential components and replace the normal ones by the kinematical constraint \eqref{eqtypeIIz2}. The modified equations are then rotated back to the canonical frame. In this last rotation the two equations (momentum conservation and kinematical constraint) become linearly combined, which is why it is cautionary to select numbers $\alpha_j$ of the same order of magnitude as that of the diagonal entries of matrix $\mathbb{A}$. The momentum equation at each node still needs to be enforced for all $\mathbf{w}_h\in V(\mathbf{q}_h)$, which is accomplished through \eqref{eqtypeIId3}. There, the contribution of $\mathbf{f}_s$ is zero because if $\mathbf{w}_h\in V(\mathbf{q}_h)$ then $\mathbf{w}_{\parallel h}=0$.
    \item In the case of finite volumes the correspondence is analogous, since the basis of the rotation strategy is that each of the $d$ momentum equations corresponding to a node in $\eta^U_\partial$ enforces momentum conservation along one cartesian direction.
    \item It is possible to define $\mathbf{n}^j$ for $j\in\eta^U_\partial$ such that all fields $\mathbf{w}_h\in W_{\parallel h}$ satisfy $\int_{\partial\mathcal{B}_h}\mathbf{w}_h\cdot \mathbf{n}=0$ exactly, where $\mathbf{n}$ is the exact normal to $\partial\mathcal{B}_h$ \cite{Montefuscolo:2012ts}. This choice guarantees that $\mathbf{u}_h$, whatever it is, does not ``create fluid mass''.
    \item It is important for finite element practitioners to notice that, since $\mathbf{n}^j$ is the same for all elements sharing node $j\in\eta^U_\partial$, all matrix manipulations above {\em can still be performed at the element level}, followed by standard assembly to build $\widetilde{\mathbb{A}}$, $\widetilde{\mathbb{G}}$ and $\underline{\widetilde{F}}(\mathbf{f}_s)$, and an additional assembly operation to build ${\mathbb{S}}^\partial$, ${\mathbb{T}}^\partial$ and ${\widetilde{B}^\partial}$.
    \item The treatment of fluid's inertia of strain-rate-dependent viscosity is exactly as for type-I squirmers. So is the extension to multiple squirmers.
\end{itemize}
%%%%%%%%%%%%%%%%%%%%%%%%%%%%%%%%%%%%%%%%%%%%%%%%%%
\subsection{Time marching}
Given a scalar or vector function of time $f$, we denote by $f^{n}$ its approximation at time level $t^{n} = n\Delta t$, with $n\in\mathbb{N}_{0}$ and time step $\Delta t > 0$. For simplicity, let us omit the suffix $h$ and restrict to the linear Stokes case. The inertia of the fluid in the examples is treated with the ALE formulation, as described by Montefuscolo et al \cite{montefuscolo2014high}.

Both types of squirmers lead to the same differential-algebraic equation (DAE), which can be written as
\begin{eqnarray}
\frac{d\mathbf{q}}{dt}&=&\mathbf{g}\paren{\mathbf{q},\mathbf{s}} \label{eqdae1}\\
\mathbb{C}(\mathbf{q})\,\left ( \begin{array}{c}
      \mathbf{s} \\ \underline{U} \\ \underline{P}
\end{array}\right ) & = & \underline{Z}(\mathbf{q}).\label{eqdae2}
\end{eqnarray}
with initial condition $\mathbf{q}(t=0)=\mathbf{q}_0$. It is important to notice that the dependence of $\mathbb{C}$ and $\underline{Z}$ on $\mathbf{q}$ is quite involved. Every time $\mathbf{q}$ is updated, the coordinates of the body and thus of the nodes on its surface change according to \eqref{eqn:rig_motion_iso}. This change is then extended to the interior nodes by some smoothing algorithm, which in our case invokes an elastic solver \cite{montefuscolo2014high}. This updated mesh is then passed to the Stokes solver to build the matrix and right-hand side following the steps described in the previous sections. This being said, any convergent scheme for DAEs could be used for \eqref{eqdae1}-\eqref{eqdae2}, in particular, in our implementation we adopted the second-order scheme described in Algorithm \ref{alg:stk_solver}.

\begin{algorithm}
\caption{Time marching for squirmer in Stokes fluid}\label{alg:stk_solver}
\begin{algorithmic}[1]
\REQUIRE $\mathbf{q}^{n-1}$, $\mathbf{s}^{n-1}$, $\mathbf{s}^{n-2}$ 
\STATE $\mathbf{q}^{n} = \mathcal{P}\cparen{\mathbf{q}^{n-1} + \Delta t\paren{\frac{3}{2}\mathbf{g}\paren{\mathbf{q}^{n-1},\mathbf{s}^{n-1}}-\frac{1}{2}\mathbf{g}\paren{\mathbf{q}^{n-2},\mathbf{s}^{n-2}}}}$ 
\STATE Update mesh and arrays $\mathbb{C}$ and $\underline{Z}$.
\STATE Find $\underline{Y}^n=\cparen{\mathbf{s}^n,\underline{U}^{n},\underline{P}^{n}}$ by solving $\mathbb{C}(\mathbf{q}^n) \underline{Y}^n=\underline{Z}(\mathbf{q}^n)$.
\ENSURE $\mathbf{q}^{n}$, $\mathbf{s}^{n}$, $\mathbf{U}^{n}$, $\mathbf{P}^{n}$
\end{algorithmic}
\end{algorithm}

The operator $\mathcal{P}$ in step 1 of the algorithm is a projection onto the configuration manifold $Q$. It is not needed for the translational degrees of freedom, that is, $\mathbf{x}_c$ is updated by
\[%\begin{equation}\label{eqn:trans_upd}
\mathbf{x}_{c}^{n} = \mathbf{x}_{c}^{n-1} + \Delta t \paren{\frac{3}{2}\mathbf{v}_{c}^{n-1} - \frac{1}{2}\mathbf{v}_{c}^{n-2}}.
\]%\end{equation}
It is also not needed for the rotational degree of freedom if $d = 2$, that is, the orientation angle $\theta$ of the body is updated by
\[%\begin{equation}\label{eqn:theta_upd}
\theta^{n} = \theta^{n-1} + \Delta t \paren{\frac{3}{2}\omega^{n-1} - \frac{1}{2}\omega^{n-2}}.
\]%\end{equation}
On the other hand, if $d = 3$, the matrix
\[%\begin{equation}\label{eqn:rotat_upd}
\mathbf{Q}^{n[0]} = \mathbf{Q}^{n-1} + \Delta t \paren{\frac32 \skk [\bomega^{n-1}] \mathbf{Q}^{n-1} - \frac12 \skk [\bomega^{n-2}] \mathbf{Q}^{n-2}}
\]%\end{equation}
will in general be only approximately orthogonal. If it is not projected back onto $\SO(d)$ one observes the bodies to loose their original shape along the simulation. We have implemented two projection algorithms that are both efficient and lead to essentially the same accuracy. The first one invokes the singular value decomposition $\mathbf{Q}^{n[0]}=\mathbf{U}\boldsymbol{\Sigma}\mathbf{V}^T$ and sets $\mathbf{Q}^n=\mathcal{P}\paren{\mathbf{Q}^{n[0]}}=\mathbf{U}\mathbf{V}^T$. The second is an iterative scheme, described in Algorithm \ref{alg:ort_proj}, taken from Sofroniou and Spaletta \cite{sofroniou2002solving} which converges very rapidly.
\begin{algorithm}
\caption{Iterative projection onto $\SO(d)$}\label{alg:ort_proj}
\begin{algorithmic}[1]
\REQUIRE $\mathbf{Q}^{n[0]}$, $\epsilon$, $k = 0$
\WHILE{$\|\mathbf{E}\|>\epsilon$}
  \STATE $\mathbf{E} = \mathbf{I}_{d} - \mathbf{Q}^{n[k]T}\mathbf{Q}^{n[k]}$
  \STATE $\mathbf{Q}^{n[k+1]} = \mathbf{Q}^{n[k]} + \frac{1}{2}\mathbf{Q}^{n[k]}\mathbf{E}$
  \STATE $k = k+1$
\ENDWHILE
\RETURN $\mathbf{Q}^{n} = \mathbf{Q}^{n[k+1]}$
\ENSURE $\mathbf{Q}^{n}$
\end{algorithmic}
\end{algorithm}

%%%%%%%%%%%%%%%%%%%%%%%%%%%%%%%%%%%%%%%%%%%%%%%%%%

%%%%%%%%%%%%%%%%%%%%%%%%%%%%%%%%%%%%%%%%%%%%%%%%%%
%%%%%%%%%%%%%%%%%%%%%%%%%%%%%%%%%%%%%%%%%%%%%%%%%%

\section{Verification experiments}

The verification of the method and code is carried out on steady squirmers of spherical shape, for which analytical, asymptotic and/or numerical solutions are available. Steady squirmers are squirmers in which $\mathbf{u}_s$ and $\mathbf{f}_s$ do not depend on time. They are adequate models for phoretic particles, while for ciliated organisms they can only capture time-averaged quantities. The modeling of the oscillatory boundary condition imposed by roaming cilia is discussed later on.

\subsection{Convergence assessment}

Convergence is assessed in Stokes flow ($\rho=0$) for a spherical squirmer of radius $R>0$, since for this case an analytical solution exists. The $3$D domain is obtained by rotating a $2$D domain $\Omega$ about the axis of symmetry. 
Let $r \geq R$ be the distance from the squirmer's centroid to any point in the fluid domain $\Omega_{f}\subset\Omega$, with $\vartheta\in\cparen{0,\pi}$ the polar coordinate measured from the direction of locomotion, $\mathbf{n}_{b} = \paren{\sin\vartheta,\cos\vartheta}$ 
the exterior normal vector and $\btau_{b} = \paren{\cos\vartheta,-\sin\vartheta}$
the polar tangent vector. We consider a type-I squirmer with imposed slip velocity $\mathbf{u}_{s} = u_{s}\btau_{b}$, being $u_{s} = B_{1}\sin\vartheta+B_{2}\sin\vartheta\cos\vartheta $ where $B_{1}, B_{2}\in\mathbb{R}$. The analytical velocity and pressure fields are given by \cite{lighthill1952squirming,blake1971spherical}

\begin{subequations}%\label{eqn:sphsqr}
\[%\begin{equation}
\begin{split}
\mathbf{u} 
=& \paren{\frac{2}{3}B_{1}\frac{R^{3}}{r^{3}}\cos\vartheta + \frac{1}{2}B_{2}\paren{\frac{R^{4}}{r^{4}} - \frac{R^{2}}{r^{2}}}\paren{3\cos^{2}\vartheta-1}}\mathbf{n}_{b}\\
&+ \paren{\frac{1}{3}B_{1}\frac{R^{3}}{r^{3}}\sin\vartheta + B_{2}\frac{R^{4}}{r^{4}}\sin\vartheta\cos\vartheta}\btau_{b},\label{eqn:sphsqru}
\end{split}
\]%\end{equation}
\[%\begin{equation}\label{eqn:sphsqrp}
p = -\mu B_{2}\frac{R^{2}}{r^{3}}\paren{3\cos^{2}\vartheta - 1},
\]%\end{equation}
\end{subequations}
and the exact swimming speed is $v_{c} = \frac{2}{3}B_{1}$ \cite{pedley2016spherical}. 

We also consider the type-II squirmer with imposed tangential force $\mathbf{f}_{s} = \frac{\mu}{R}\paren{2B_{1}\sin\vartheta+5B_{2}\sin\vartheta\cos\vartheta}\btau_{b}$, which produces the same exact solution.

Depending on the sign of the parameter $\beta = \frac{B_{2}}{B_{1}}$, the squirmer can be classified as \textit{neutral} if $\beta = 0$ (see the numerical streamlines in Figure \ref{fig:sqrmodelneutral}, in all cases $B_1>0$), \textit{pusher} if $\beta < 0$ (Figure \ref{fig:sqrmodelpusher}) and \textit{puller} if $\beta > 0$ (Figure \ref{fig:sqrmodelpuller}) \cite{pedley2016spherical}. When ${\beta} \neq 0$ a region of recirculation is created in the front (if pusher) or the back (if puller) of the squirmer, due to the change of sign of $u_{s}$ for some $\vartheta\in\paren{0,\pi}$.

%%%%%%%%%%%%%%%%%%%%%%%%%%%%%%%%%%%%%%%%
\newcommand\cw{1}
\newcommand\bw{.192}
\newcommand\dw{.24}
\newcommand\fw{.75}
\setlength{\fboxsep}{0pt}
\setlength{\fboxrule}{0pt}
%%%%%%%%%%%%%%%%%%%%%%%%%%%%%%%%%%%%%%%%
\begin{figure}[pt]
\begin{center}
    \begin{subfigure}{\dw\textwidth}%%%%%%%%%%%%%%%%%%%%
    \begin{center}
    $\beta = 0$\par
    \begin{tikzpicture}
    \node[anchor=south west,inner sep=0] (image) at (0,0) {\fbox{\includegraphics[width=\textwidth,trim=0cm 10cm 0cm 10cm,clip]{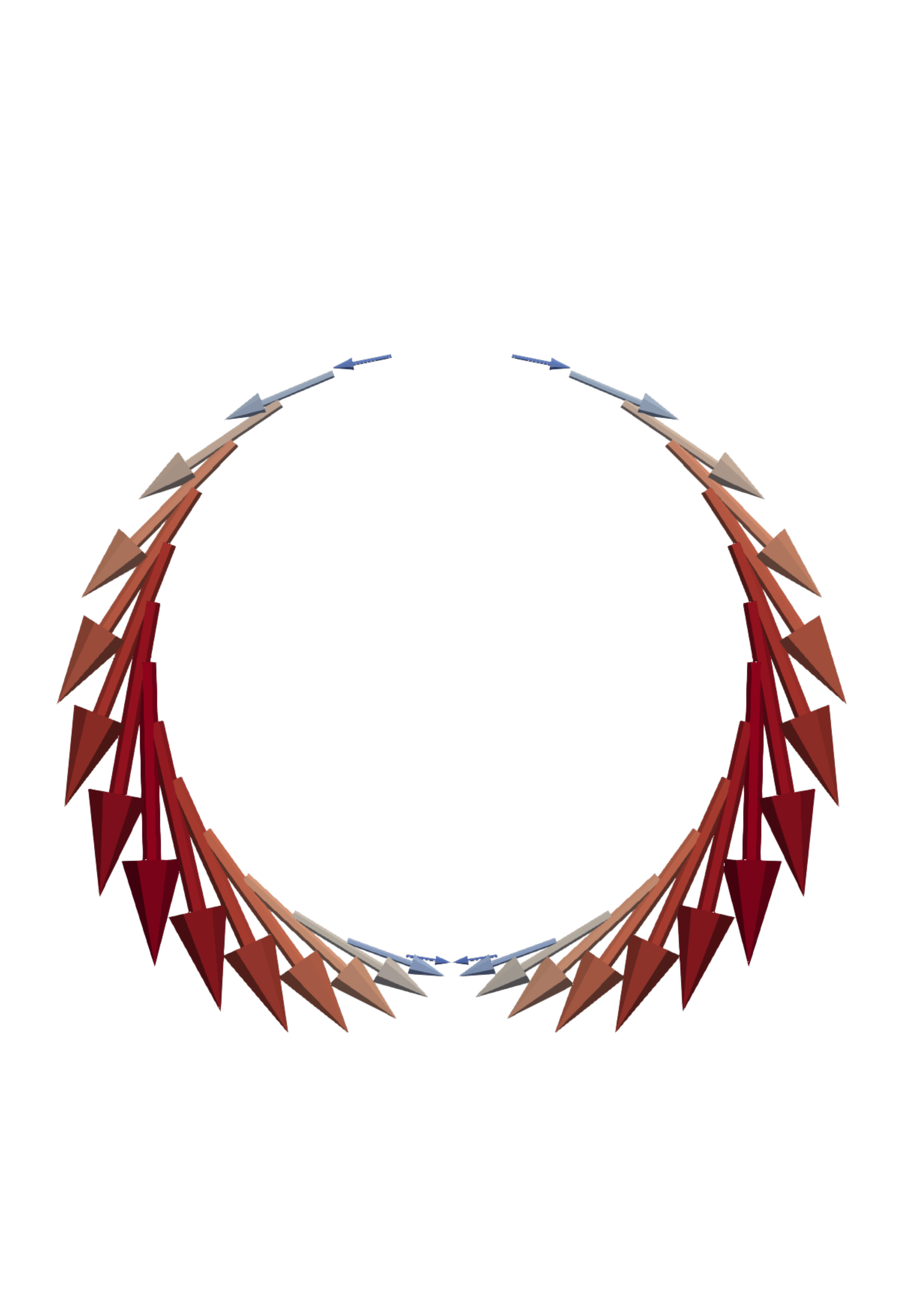}}};
    \begin{scope}[x={(image.south east)},y={(image.north west)}]
    \draw[->,>=stealth,thick] (0.5,0.5) -- (0.5,0.8);
    \node at (0.40,0.60) {\normalsize{$\mathbf{v}_{c}$}};
    \node at (0.15,0.85) {\normalsize{$\mathbf{u}_{s}$}};
    \node at (0.56,0.67) {\normalsize{$\vartheta$}};
    \draw [->] (0.5,0.6) arc (90:65:15pt);
    \node at (0.85,0.85) {\normalsize{$r$}};
    \draw[->,>=stealth,thick] (0.5,0.5) -- (0.80,0.80);
    \end{scope}
    \end{tikzpicture}\par
    \begin{tikzpicture}
    \node[anchor=south west,inner sep=0] (image) at (0,0) {\fbox{\includegraphics[width=\textwidth,trim=0cm 10cm 0cm 10cm,clip]{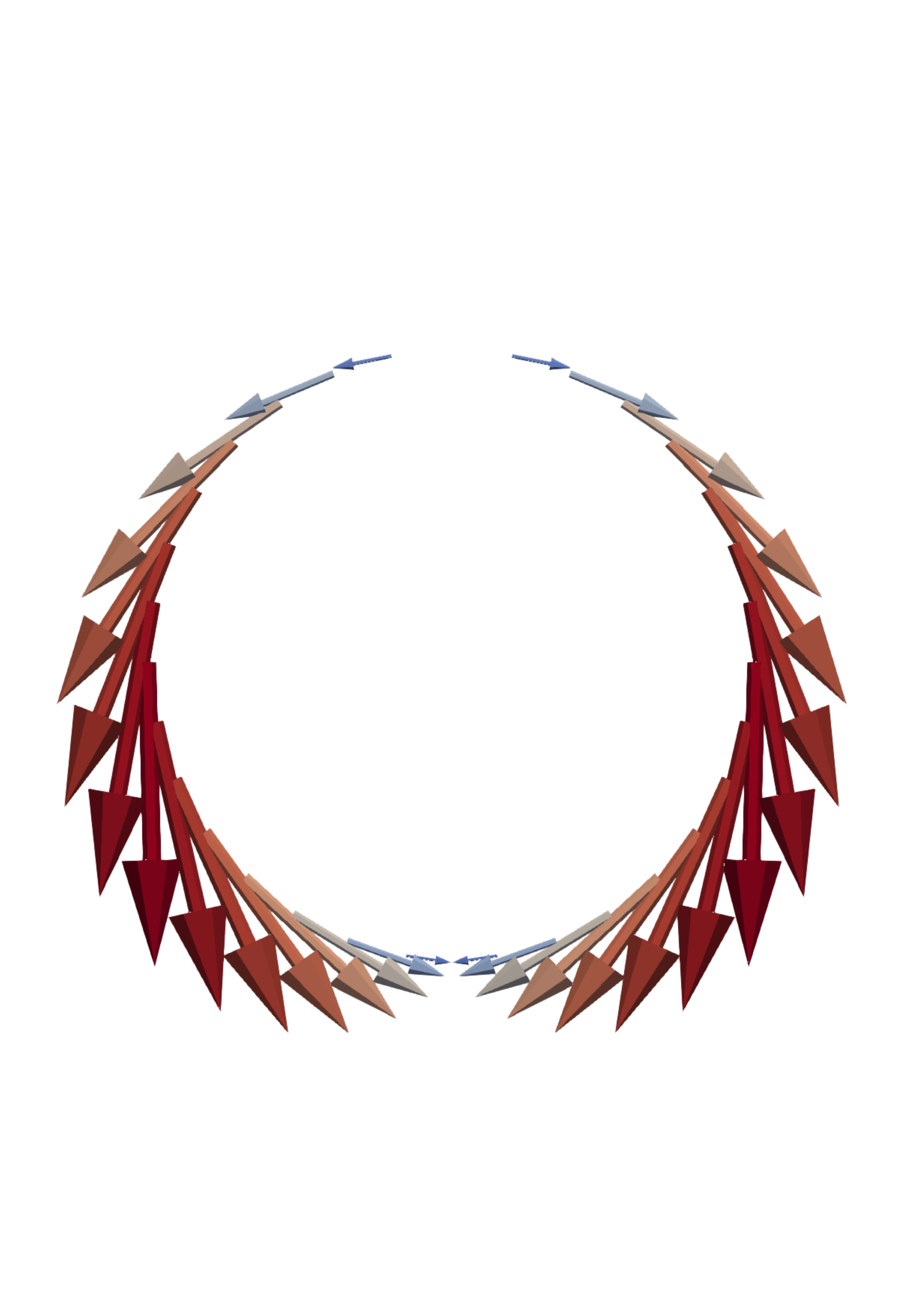}}};
    \begin{scope}[x={(image.south east)},y={(image.north west)}]
    \draw[->,>=stealth,thick] (0.5,0.5) -- (0.5,0.8);
    \node at (0.40,0.60) {\normalsize{$\mathbf{v}_{c}$}};
    \node at (0.15,0.85) {\normalsize{$\mathbf{f}_{s}$}};
    \end{scope}
    \end{tikzpicture}
    \end{center}
    \end{subfigure}
    \begin{subfigure}{\fw\textwidth}%%%%%%%%%%%%%%%%%%%%
    \begin{center}
    \foreach \y in {eul,lag}{
    \begin{tikzpicture}
    \node[anchor=south west,inner sep=0] (image) at (0,0) {\fbox{\includegraphics[width=.485\textwidth,trim=0cm 0cm 0cm 0cm,clip]{beta0\y.pdf}}};
    \begin{scope}[x={(image.south east)},y={(image.north west)}]
    \draw[->,>=stealth,thick] (0.5,0.5) -- (0.5,0.6);
    \end{scope}
    \end{tikzpicture}
    }
    \end{center}
    \end{subfigure}
\end{center}
\caption{{\em Neutral squirmer ($\beta=0$): } Slip velocity $\mathbf{u}_{s}$ and tangential force $\mathbf{f}_{s}$, indicating the radial and polar directions, $r$ and $\vartheta$ respectively, and the direction of movement $\mathbf{v}_{c}$ (left). Numerically-obtained streamlines (light blue) and some velocity vectors (dark blue) in the laboratory, fixed, frame (center) and in a frame moving with the particle (right).} \label{fig:sqrmodelneutral}
\end{figure}
%%%%%%%%%%%%%%%%%%%%%%%%%%%%%%%%%%%%%%%%
\begin{figure}[!ht]
\begin{center}
    \begin{subfigure}{\dw\textwidth}%%%%%%%%%%%%%%%%%%%%
    \begin{center}
    $\beta = -0.5$\par
    \begin{tikzpicture}
    \node[anchor=south west,inner sep=0] (image) at (0,0) {\fbox{\includegraphics[width=\textwidth,trim=0cm 10cm 0cm 10cm,clip]{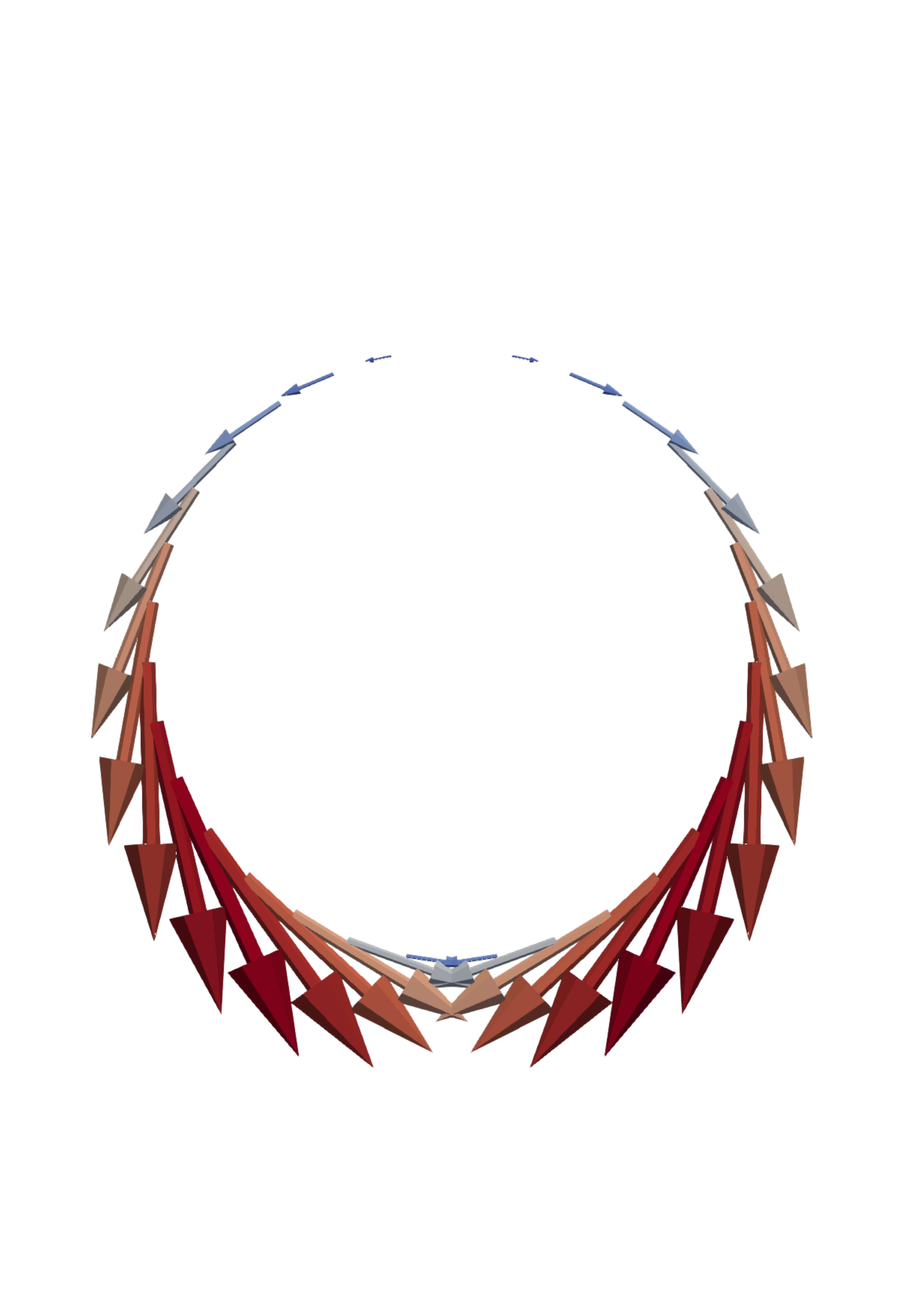}}};
    \begin{scope}[x={(image.south east)},y={(image.north west)}]
    \draw[->,>=stealth,thick] (0.5,0.5) -- (0.5,0.8);
    \node at (0.40,0.60) {\normalsize{$\mathbf{v}_{c}$}};
    \node at (0.15,0.85) {\normalsize{$\mathbf{u}_{s}$}};
    \node at (0.56,0.67) {\normalsize{$\vartheta$}};
    \draw [->] (0.5,0.6) arc (90:65:15pt);
    \node at (0.85,0.85) {\normalsize{$r$}};
    \draw[->,>=stealth,thick] (0.5,0.5) -- (0.80,0.80);
    \end{scope}
    \end{tikzpicture}\par
    \begin{tikzpicture}
    \node[anchor=south west,inner sep=0] (image) at (0,0) {\fbox{\includegraphics[width=\textwidth,trim=0cm 10cm 0cm 10cm,clip]{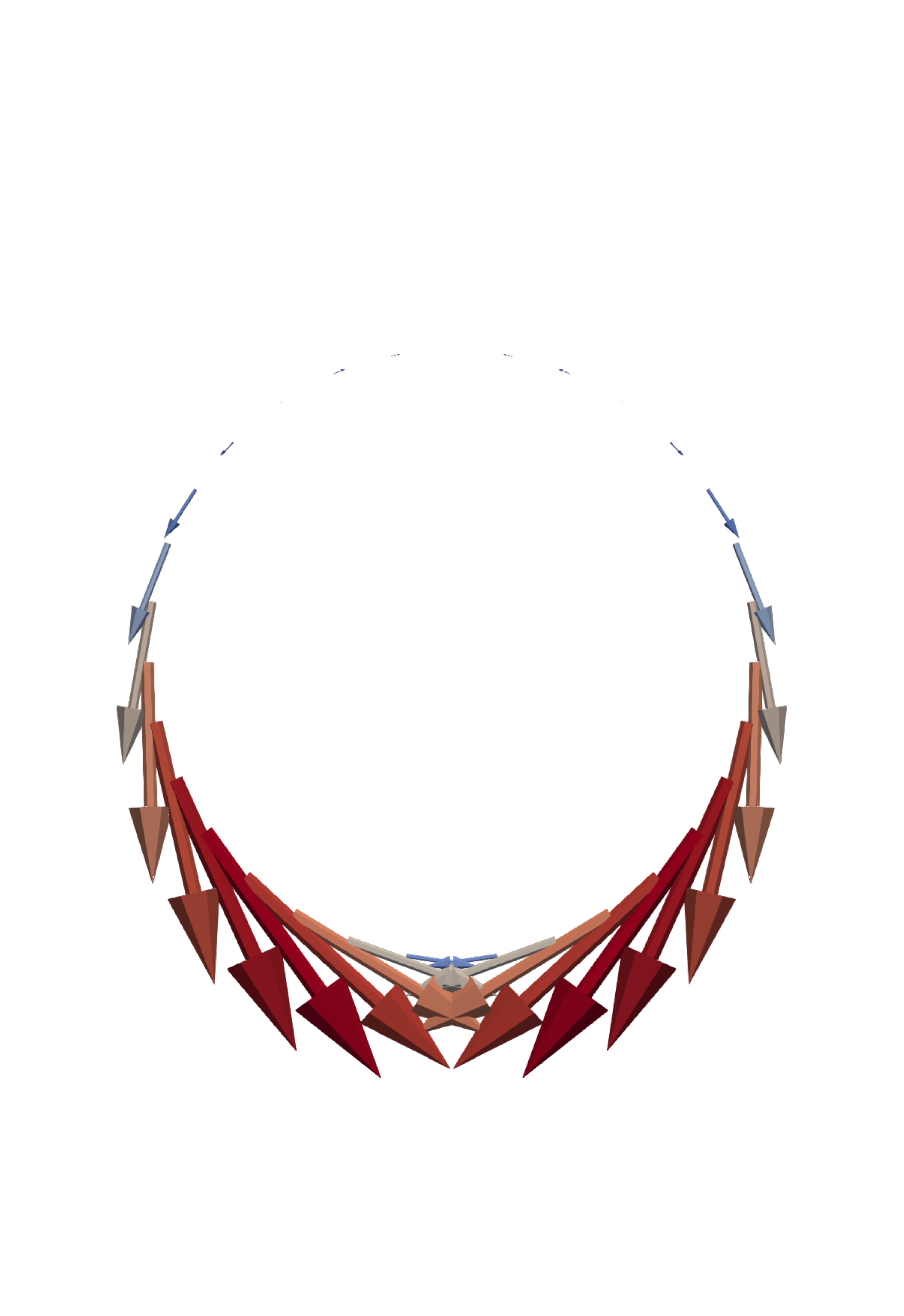}}};
    \begin{scope}[x={(image.south east)},y={(image.north west)}]
    \draw[->,>=stealth,thick] (0.5,0.5) -- (0.5,0.8);
    \node at (0.40,0.60) {\normalsize{$\mathbf{v}_{c}$}};
    \node at (0.15,0.85) {\normalsize{$\mathbf{f}_{s}$}};
    \end{scope}
    \end{tikzpicture}
    \end{center}
    \end{subfigure}
    \begin{subfigure}{\fw\textwidth}%%%%%%%%%%%%%%%%%%%%
    \begin{center}
    \foreach \y in {eul,lag}{
    \begin{tikzpicture}
    \node[anchor=south west,inner sep=0] (image) at (0,0) {\fbox{\includegraphics[width=.485\textwidth,trim=0cm 0cm 0cm 0cm,clip]{betan05\y.pdf}}};
    \begin{scope}[x={(image.south east)},y={(image.north west)}]
    \draw[->,>=stealth,thick] (0.5,0.5) -- (0.5,0.6);
    \end{scope}
    \end{tikzpicture}
    }
    \end{center}
    \end{subfigure}\par
    \vspace{.5\baselineskip}
    \begin{subfigure}{\dw\textwidth}%%%%%%%%%%%%%%%%%%%%
    \begin{center}
    $\beta = -5$\par
    \begin{tikzpicture}
    \node[anchor=south west,inner sep=0] (image) at (0,0) {\fbox{\includegraphics[width=\textwidth,trim=0cm 10cm 0cm 10cm,clip]{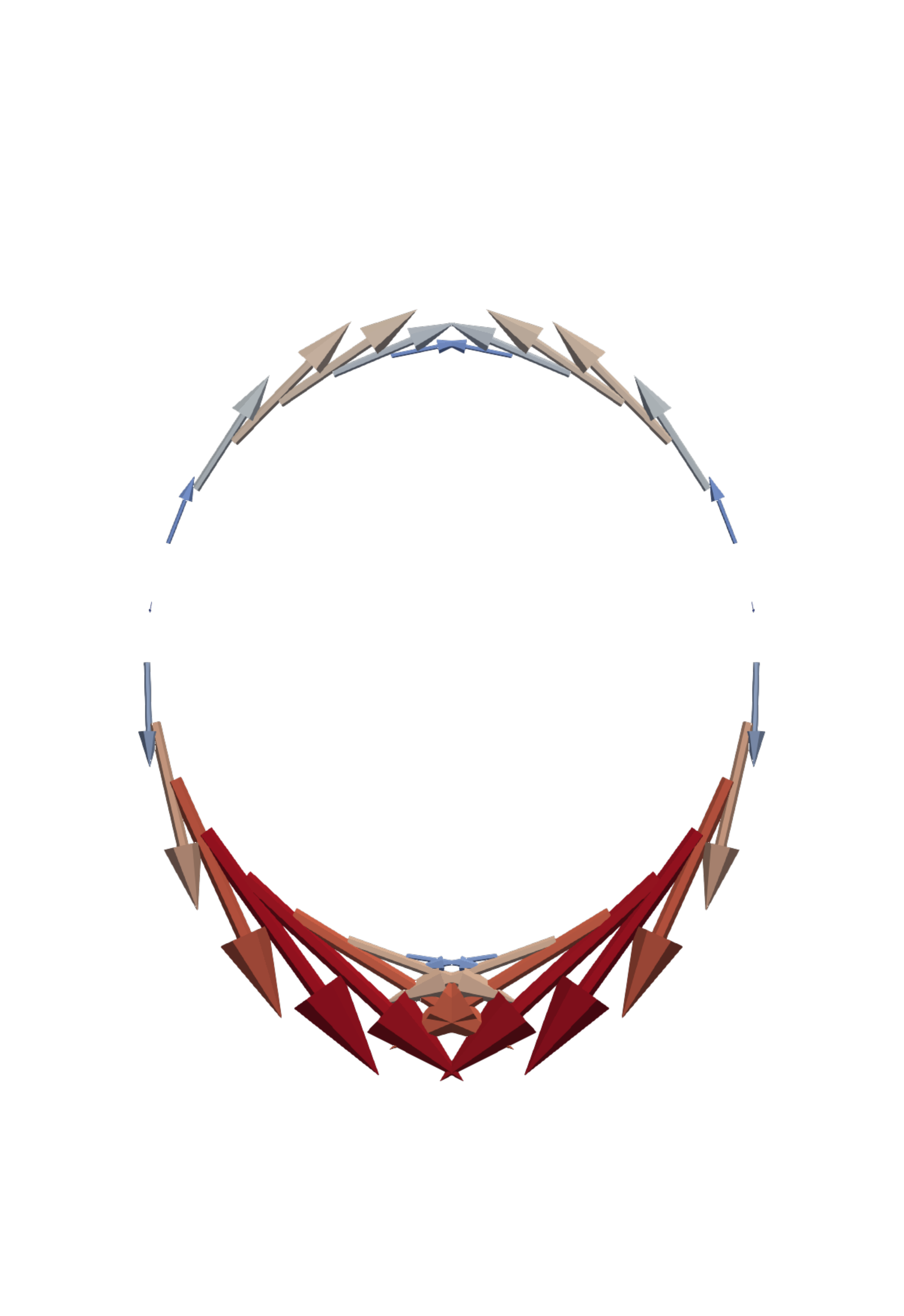}}};
    \begin{scope}[x={(image.south east)},y={(image.north west)}]
    \draw[->,>=stealth,thick] (0.5,0.5) -- (0.5,0.8);
    \node at (0.40,0.60) {\normalsize{$\mathbf{v}_{c}$}};
    \node at (0.15,0.85) {\normalsize{$\mathbf{u}_{s}$}};
    \node at (0.56,0.67) {\normalsize{$\vartheta$}};
    \draw [->] (0.5,0.6) arc (90:65:15pt);
    \node at (0.85,0.85) {\normalsize{$r$}};
    \draw[->,>=stealth,thick] (0.5,0.5) -- (0.80,0.80);
    \end{scope}
    \end{tikzpicture}\par
    \begin{tikzpicture}
    \node[anchor=south west,inner sep=0] (image) at (0,0) {\fbox{\includegraphics[width=\textwidth,trim=0cm 10cm 0cm 10cm,clip]{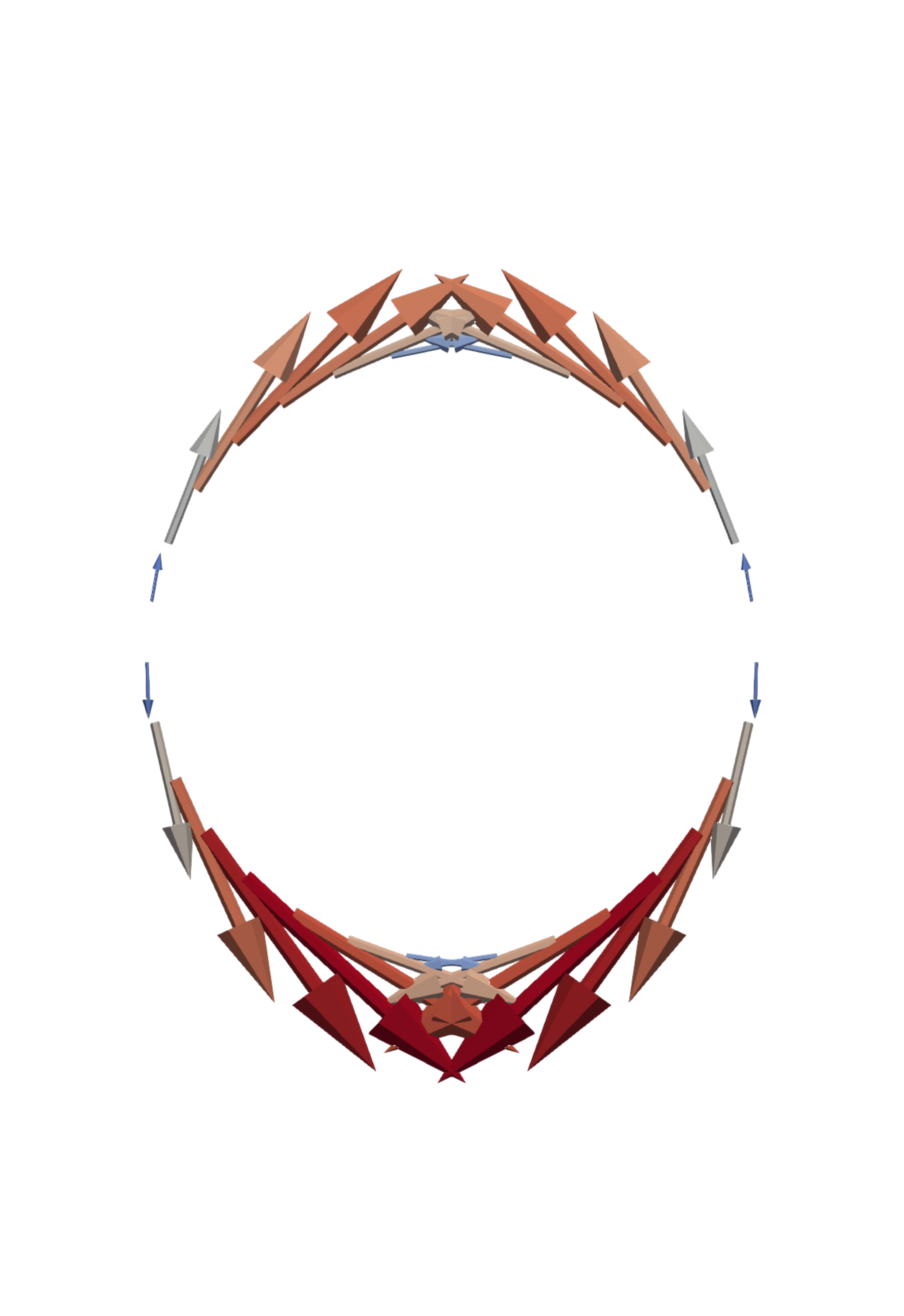}}};
    \begin{scope}[x={(image.south east)},y={(image.north west)}]
    \draw[->,>=stealth,thick] (0.5,0.5) -- (0.5,0.8);
    \node at (0.40,0.60) {\normalsize{$\mathbf{v}_{c}$}};
    \node at (0.15,0.85) {\normalsize{$\mathbf{f}_{s}$}};
    \end{scope}
    \end{tikzpicture}
    \end{center}
    \end{subfigure}
    \begin{subfigure}{\fw\textwidth}%%%%%%%%%%%%%%%%%%%%
    \begin{center}
    \foreach \y in {eul,lag}{
    \begin{tikzpicture}
    \node[anchor=south west,inner sep=0] (image) at (0,0) {\fbox{\includegraphics[width=.485\textwidth,trim=0cm 0cm 0cm 0cm,clip]{betan5\y.pdf}}};
    \begin{scope}[x={(image.south east)},y={(image.north west)}]
    \draw[->,>=stealth,thick] (0.5,0.5) -- (0.5,0.6);
    \end{scope}
    \end{tikzpicture}
    }
    \end{center}
    \end{subfigure}
\end{center}
\caption{{\em Pushers:} Same as Figure \ref{fig:sqrmodelneutral} for $\beta=-0.5$ (top) and $\beta=-5$ (bottom).} \label{fig:sqrmodelpusher}
\end{figure}
%%%%%%%%%%%%%%%%%%%%%%%%%%%%%%%%%%%%%%%%
\begin{figure}[!ht]
\begin{center}
    \begin{subfigure}{\dw\textwidth}%%%%%%%%%%%%%%%%%%%%
    \begin{center}
    $\beta = +0.5$\par
    \begin{tikzpicture}
    \node[anchor=south west,inner sep=0] (image) at (0,0) {\fbox{\includegraphics[width=\textwidth,trim=0cm 10cm 0cm 10cm,clip]{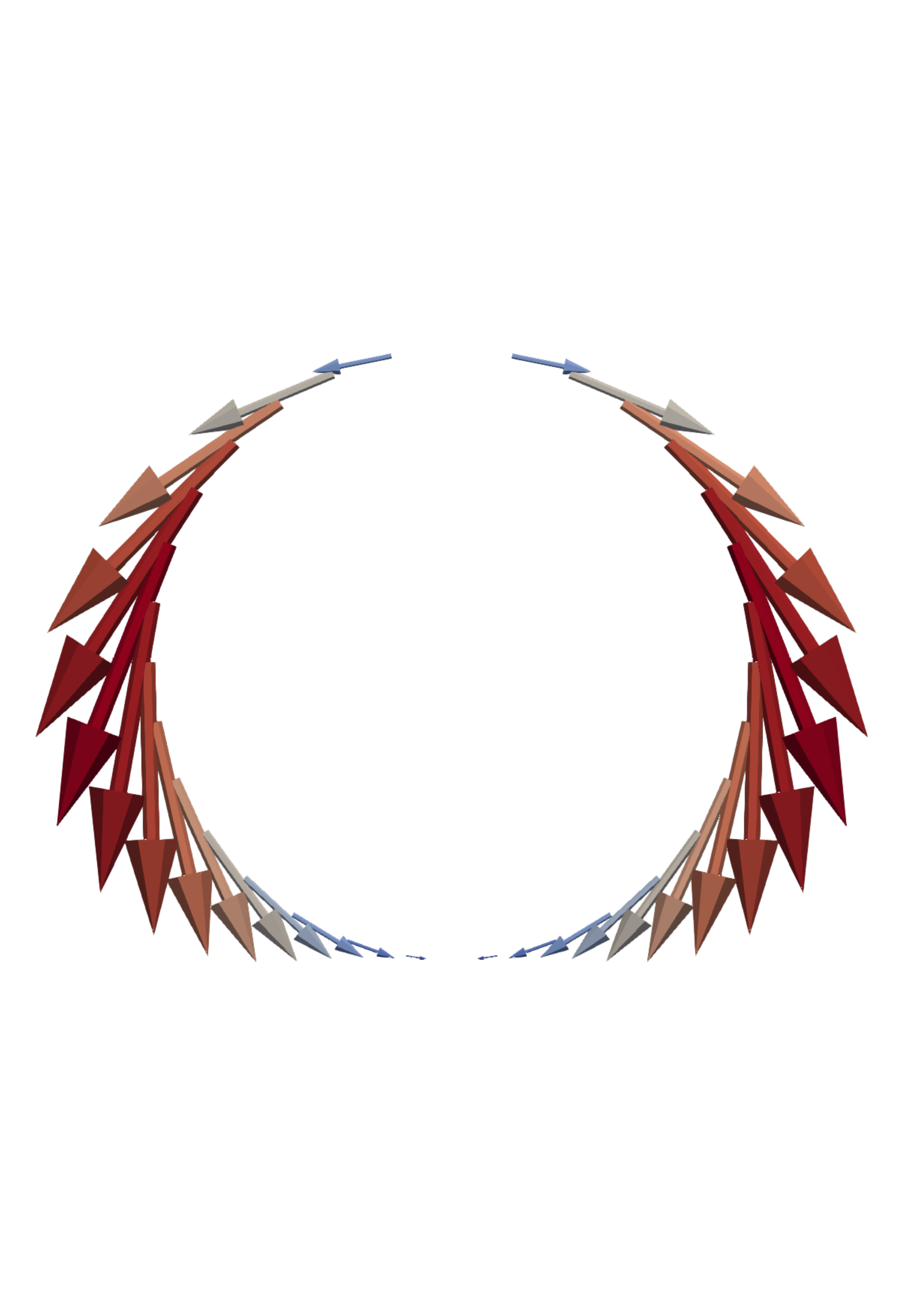}}};
    \begin{scope}[x={(image.south east)},y={(image.north west)}]
    \draw[->,>=stealth,thick] (0.5,0.5) -- (0.5,0.8);
    \node at (0.40,0.60) {\normalsize{$\mathbf{v}_{c}$}};
    \node at (0.15,0.85) {\normalsize{$\mathbf{u}_{s}$}};
    \node at (0.56,0.67) {\normalsize{$\vartheta$}};
    \draw [->] (0.5,0.6) arc (90:65:15pt);
    \node at (0.85,0.85) {\normalsize{$r$}};
    \draw[->,>=stealth,thick] (0.5,0.5) -- (0.80,0.80);
    \end{scope}
    \end{tikzpicture}\par
    \begin{tikzpicture}
    \node[anchor=south west,inner sep=0] (image) at (0,0) {\fbox{\includegraphics[width=\textwidth,trim=0cm 10cm 0cm 10cm,clip]{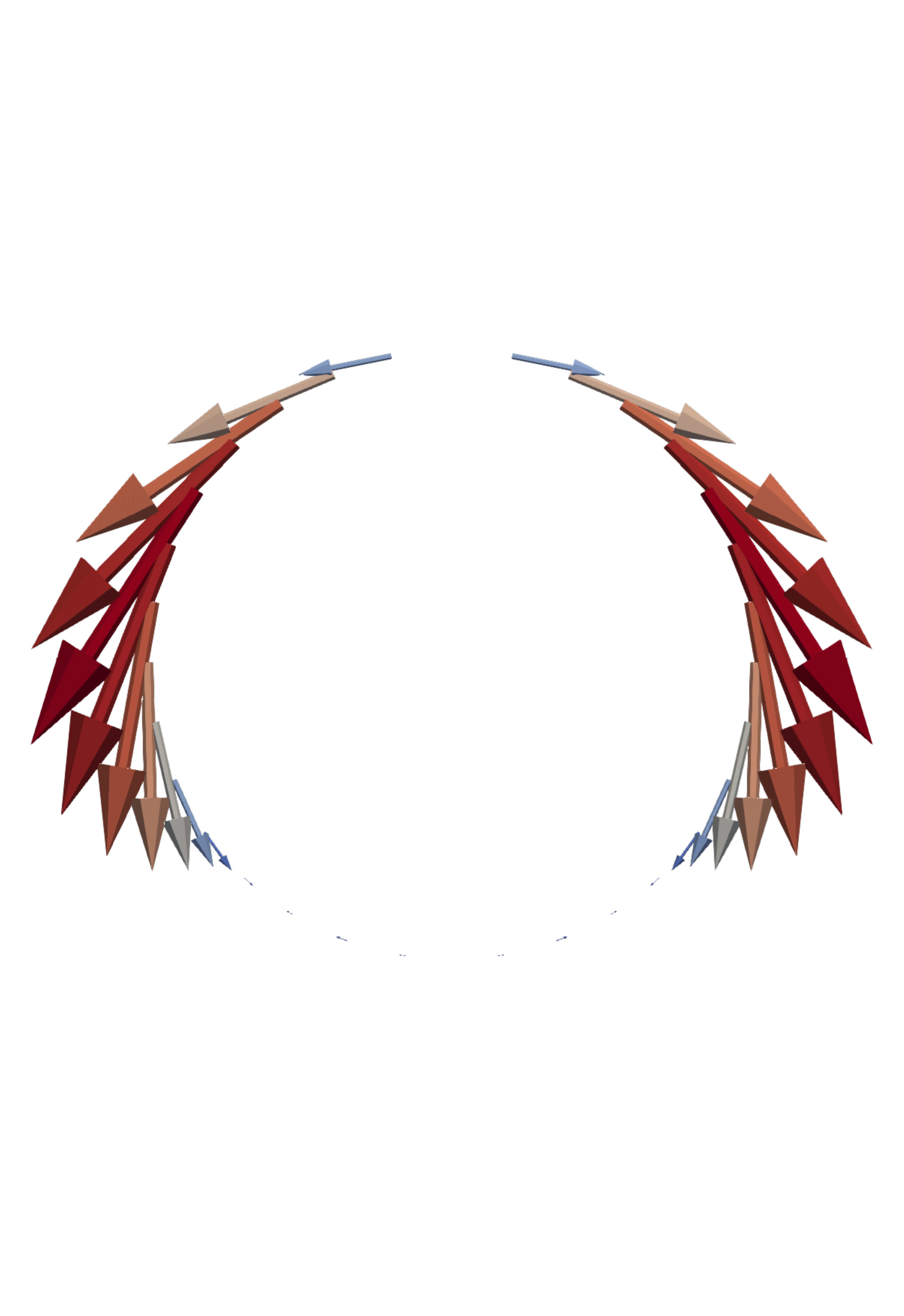}}};
    \begin{scope}[x={(image.south east)},y={(image.north west)}]
    \draw[->,>=stealth,thick] (0.5,0.5) -- (0.5,0.8);
    \node at (0.40,0.60) {\normalsize{$\mathbf{v}_{c}$}};
    \node at (0.15,0.85) {\normalsize{$\mathbf{f}_{s}$}};
    \end{scope}
    \end{tikzpicture}
    \end{center}
    \end{subfigure}
    \begin{subfigure}{\fw\textwidth}%%%%%%%%%%%%%%%%%%%%
    \begin{center}
    \foreach \y in {eul,lag}{
    \begin{tikzpicture}
    \node[anchor=south west,inner sep=0] (image) at (0,0) {\fbox{\includegraphics[width=.485\textwidth,trim=0cm 0cm 0cm 0cm,clip]{betap05\y.pdf}}};
    \begin{scope}[x={(image.south east)},y={(image.north west)}]
    \draw[->,>=stealth,thick] (0.5,0.5) -- (0.5,0.6);
    \end{scope}
    \end{tikzpicture}
    }
    \end{center}
    \end{subfigure}\par
    \vspace{.5\baselineskip}
    \begin{subfigure}{\dw\textwidth}%%%%%%%%%%%%%%%%%%%%
    \begin{center}
    $\beta = +5$\par
    \begin{tikzpicture}
    \node[anchor=south west,inner sep=0] (image) at (0,0) {\fbox{\includegraphics[width=\textwidth,trim=0cm 10cm 0cm 10cm,clip]{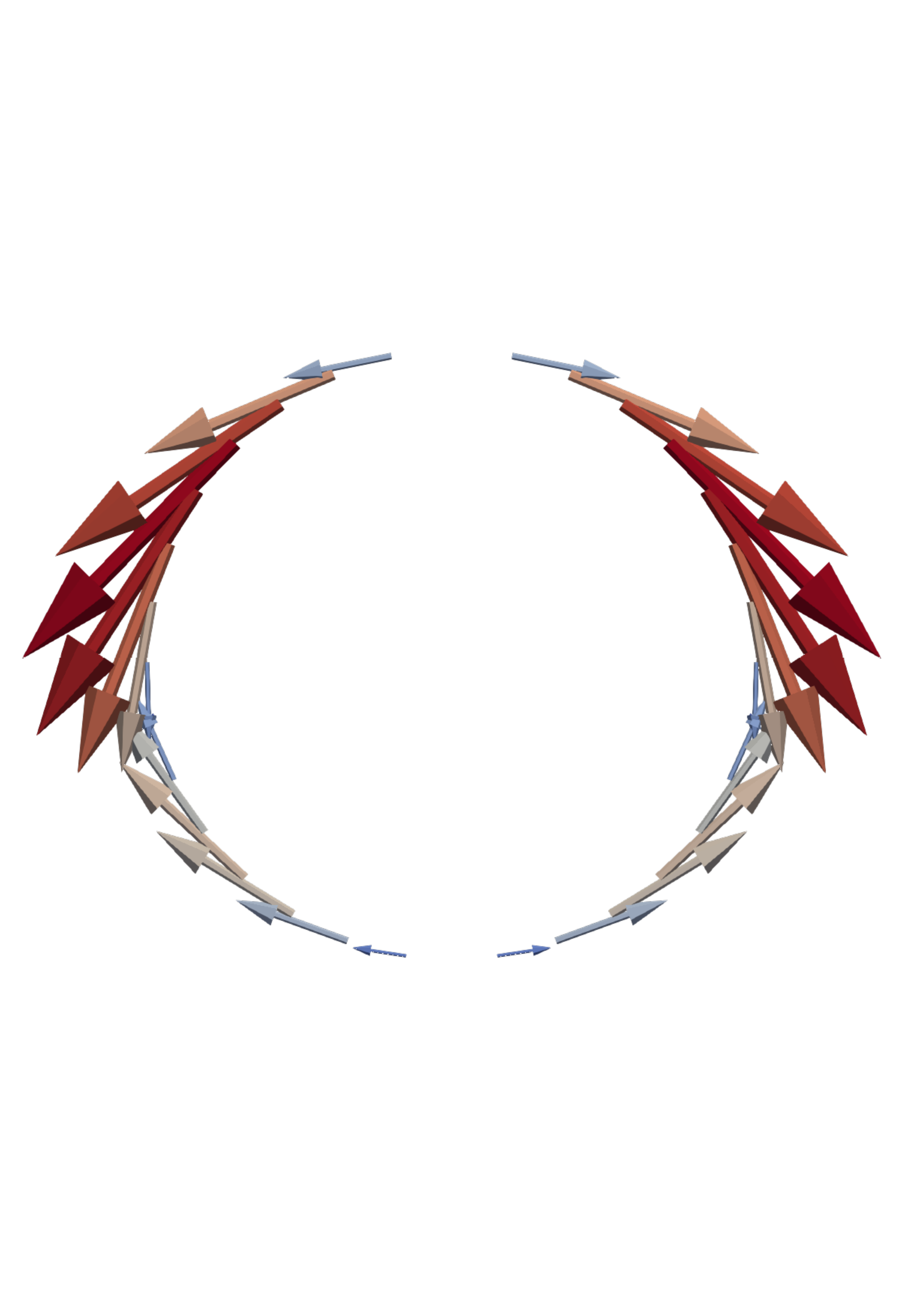}}};
    \begin{scope}[x={(image.south east)},y={(image.north west)}]
    \draw[->,>=stealth,thick] (0.5,0.5) -- (0.5,0.8);
    \node at (0.40,0.60) {\normalsize{$\mathbf{v}_{c}$}};
    \node at (0.15,0.85) {\normalsize{$\mathbf{u}_{s}$}};
    \node at (0.56,0.67) {\normalsize{$\vartheta$}};
    \draw [->] (0.5,0.6) arc (90:65:15pt);
    \node at (0.85,0.85) {\normalsize{$r$}};
    \draw[->,>=stealth,thick] (0.5,0.5) -- (0.80,0.80);
    \end{scope}
    \end{tikzpicture}\par
    \begin{tikzpicture}
    \node[anchor=south west,inner sep=0] (image) at (0,0) {\fbox{\includegraphics[width=\textwidth,trim=0cm 10cm 0cm 10cm,clip]{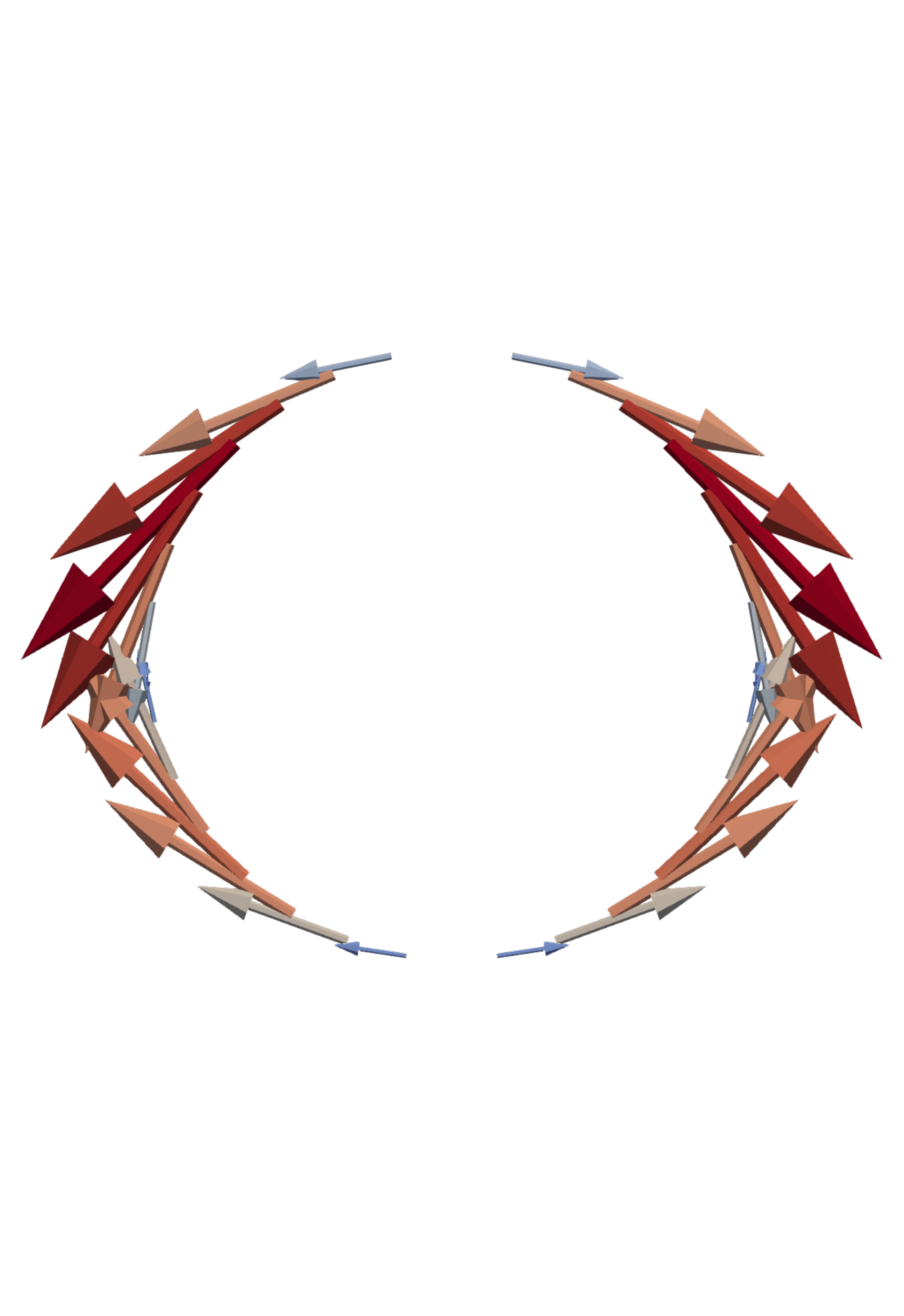}}};
    \begin{scope}[x={(image.south east)},y={(image.north west)}]
    \draw[->,>=stealth,thick] (0.5,0.5) -- (0.5,0.8);
    \node at (0.40,0.60) {\normalsize{$\mathbf{v}_{c}$}};
    \node at (0.15,0.85) {\normalsize{$\mathbf{f}_{s}$}};
    \end{scope}
    \end{tikzpicture}
    \end{center}
    \end{subfigure}
    \begin{subfigure}{\fw\textwidth}%%%%%%%%%%%%%%%%%%%%
    \begin{center}
    \foreach \y in {eul,lag}{
    \begin{tikzpicture}
    \node[anchor=south west,inner sep=0] (image) at (0,0) {\fbox{\includegraphics[width=.485\textwidth,trim=0cm 0cm 0cm 0cm,clip]{betap5\y.pdf}}};
    \begin{scope}[x={(image.south east)},y={(image.north west)}]
    \draw[->,>=stealth,thick] (0.5,0.5) -- (0.5,0.6);
    \end{scope}
    \end{tikzpicture}
    }
    \end{center}
    \end{subfigure}
\end{center}
\caption{{\em Pullers:} Same as Figure \ref{fig:sqrmodelneutral} for $\beta=+0.5$ (top) and $\beta=+5$ (bottom).}\label{fig:sqrmodelpuller}
\end{figure}
%%%%%%%%%%%%%%%%%%%%%%%%%%%%%%%%%%%%%%%%

\begin{figure}[!ht]
\centering
\includegraphics[width = .6\textwidth]{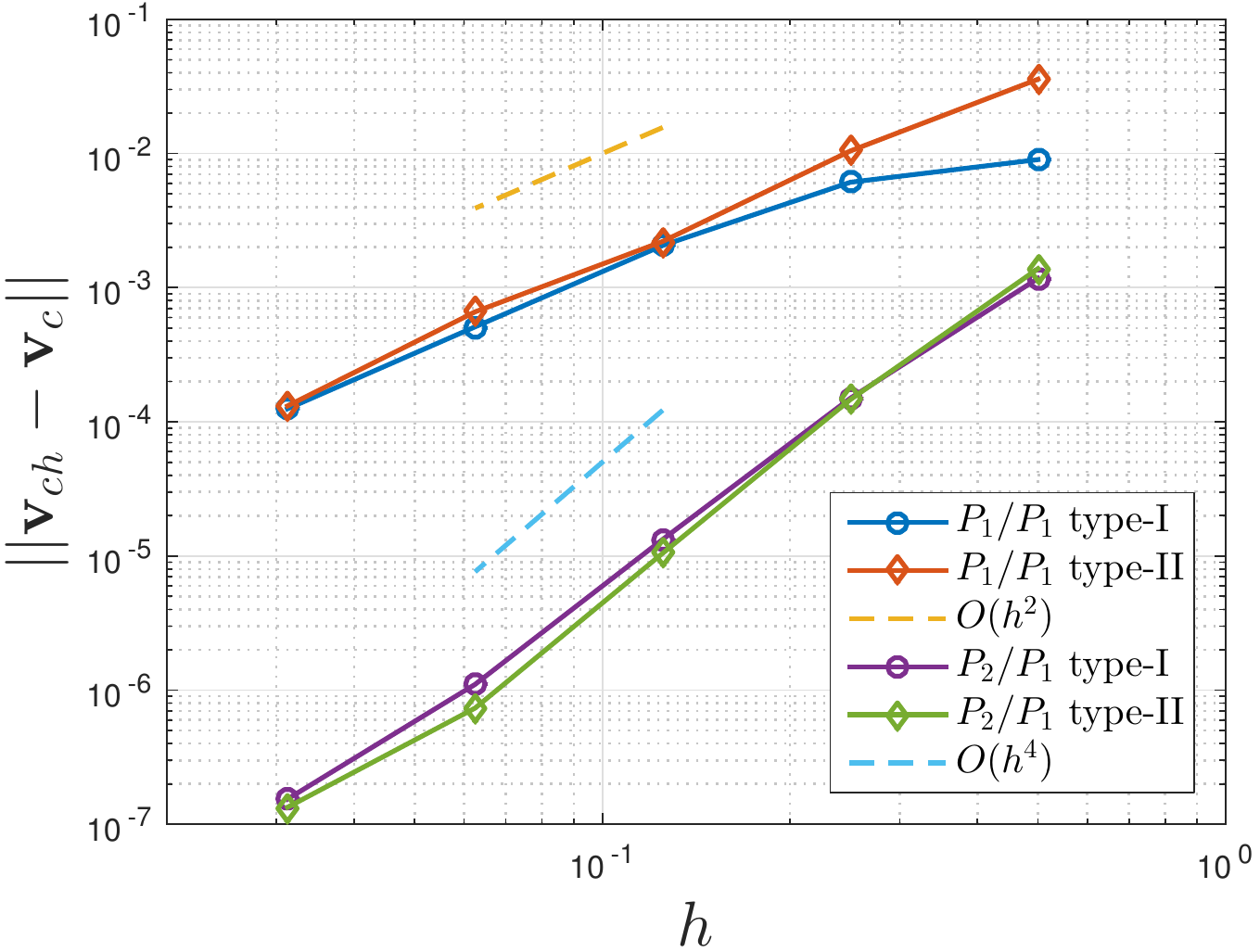}
\caption{Mesh convergence of translational velocity of an axisymmetric squirmer.} \label{fig:bnp_cmp}
\end{figure}

The problem was solved for a squirmer of radius $R=1$ inside a domain of size 300. The coarsest mesh used, corresponding to refinement $k=0$ in the tables, had element size of $h_0=0.5$ close to the squirmer and much larger away from it, totalling 388 elements. The size of the elements close to the squirmer in finer meshes is approximately $h = h_0 2^{-k}$, where $k = 1,\dots,5$ is the index of refinement. The finest mesh ($k=5$) contained 337792 triangles. The case reported here corresponds to $B_2=0$, but the results are similar for pullers and pushers. Figure \ref{fig:bnp_cmp} shows the convergence of the translational velocity, which exhibits second order for the GLS (Galerkin Least Squares \cite{hughes1986new}) stabilized $P_{1}/P_{1}$ elements and fourth order for $P_{2}/P_{1}$ elements.

Focusing on the stabilized $P_{1}/P_{1}$ case, which is the only one considered hereafter, second and first order of convergence, in the $L^{2}$-norm, are observed for the fluid velocity and fluid pressure, res\-pectively. These orders hold both for the type-I squirmer (Table \ref{tab:conv_dtp}) and for the type-II one (Table \ref{tab:conv_ntp}). Also shown are the errors in the $L^{\infty}$-norm, which exhibits roughly similar, though more erratic, behavior.

\begin{table*}[!h]
\caption{Convergence study for type-I squirmer ($P_1/P_1$ stabilized element).}\label{tab:conv_dtp}
\centering
\begin{tabular}{c|cc|cc|cc|cc}
\toprule
$k$ & $\norm{\mathbf{u}-\mathbf{u}_{h}}_{L^{2}}$ & Order & $\norm{p-p_{h}}_{L^{2}}$ & Order & $\norm{\mathbf{u}-\mathbf{u}_{h}}_{L^{\infty}}$ & Order & $\norm{p-p_{h}}_{L^{\infty}}$ & Order \\
\midrule
0 & 1.1321e-01 &        & 2.9923e-01 &        & 1.1889e-01 &        & 2.6134e-01 &        \\
1 & 5.4979e-02 & 1.0421 & 2.1324e-01 & 0.4888 & 5.9124e-02 & 1.0079 & 3.5854e-01 & -0.4562\\
2 & 1.8638e-02 & 1.5606 & 1.2364e-01 & 0.7863 & 2.7356e-02 & 1.1119 & 2.8143e-01 & 0.3493 \\
3 & 4.9802e-03 & 1.9040 & 5.6942e-02 & 1.1186 & 1.0345e-02 & 1.4028 & 2.1074e-01 & 0.4173 \\
4 & 1.3595e-03 & 1.8731 & 2.0263e-02 & 1.4906 & 2.5764e-03 & 2.0056 & 9.8449e-02 & 1.0981 \\
5 & 4.9656e-04 & 1.4531 & 8.1991e-03 & 1.3053 & 1.0501e-03 & 1.2948 & 5.6718e-02 & 0.7956 \\
\bottomrule
\end{tabular}
\end{table*}

\begin{table*}[!h]
\caption{Convergence study for type-II squirmer ($P_1/P_1$ stabilized element).}\label{tab:conv_ntp}
\centering
\begin{tabular}{c|cc|cc|cc|cc}
\toprule
$k$ & $\norm{\mathbf{u}-\mathbf{u}_{h}}_{L^{2}}$ & Order & $\norm{p-p_{h}}_{L^{2}}$ & Order & $\norm{\mathbf{u}-\mathbf{u}_{h}}_{L^{\infty}}$ & Order & $\norm{p-p_{h}}_{L^{\infty}}$ & Order \\
\midrule
0 & 1.4659e-01 &        & 2.2759e-01 &        & 1.1958e-01 &        &  1.9774e-01 &        \\
1 & 8.5926e-02 & 0.7707 & 1.7848e-01 & 0.3507 & 6.4272e-02 & 0.8958 &  3.3336e-01 & -0.7535\\
2 & 2.9136e-02 & 1.5603 & 1.1530e-01 & 0.6303 & 2.5722e-02 & 1.3212 &  2.8250e-01 & 0.2388 \\
3 & 7.8808e-03 & 1.8864 & 5.4555e-01 & 1.0797 & 8.1751e-03 & 1.6737 &  1.8417e-01 & 0.6172 \\
4 & 2.2595e-03 & 1.8023 & 1.9556e-02 & 1.4800 & 2.3110e-03 & 1.8227 &  1.1186e-01 & 0.7194 \\
5 & 6.4311e-04 & 1.8129 & 8.0467e-03 & 1.2812 & 8.2458e-04 & 1.4868 &  6.1905e-02 & 0.8536 \\
\bottomrule
\end{tabular}
\end{table*}
%%%%%%%%%%%%%%%%%%%%%%%%%%%%%%%%%%%%%%%%
\subsection{Squirmer at finite Reynolds number}

For spherical steady squirmers of type-I there also exist results at finite Reynolds number ($\Rey= \frac{\rho v_{c} R}{\mu}$) that are used here for verification. Chisholm et. al. \cite{chisholm2016squirmer} conducted a careful numerical study, which showed that inertial effects monotonically increase the speed $v_{c}$ of a pusher, while for a puller the speed decreases at first and then, for sufficiently high $\beta$ and $\Rey$, it starts increasing again with $\Rey$. Also available for comparison are the asymptotic expansions of order $O\paren{\Rey}$ and $O\paren{\Rey^{2}}$ found by Wang and Ardekani \cite{wang2012inertial} and by Khair et. al. \cite{khair2014expansions}, respectively.

\begin{figure}[!ht]%[pt]
\begin{center}
%\begin{subfigure}[t]{0.5\textwidth}
\includegraphics[width=.75\textwidth]{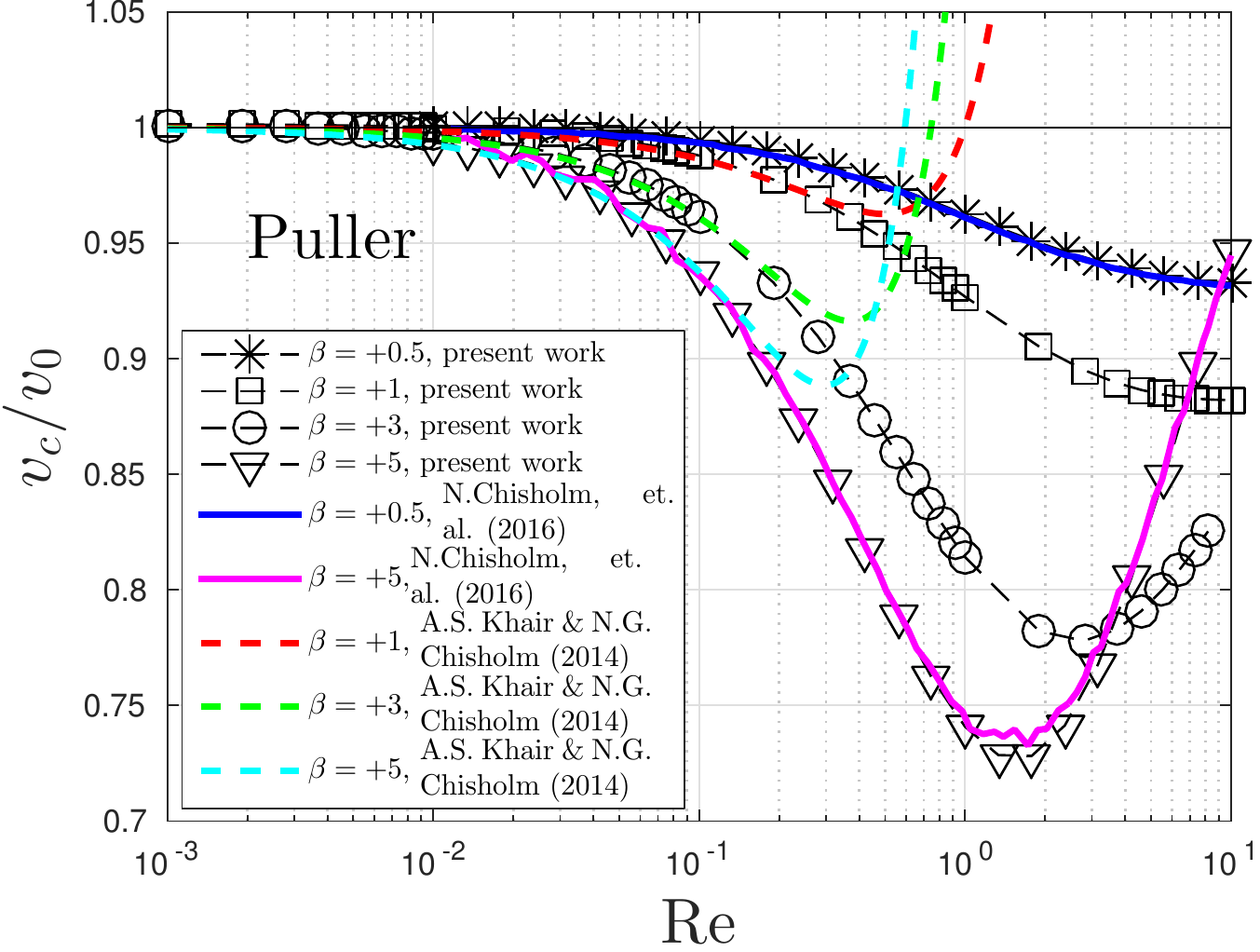}
%\end{subfigure}
%\begin{subfigure}[t]{0.5\textwidth}
\includegraphics[width=.75\textwidth]{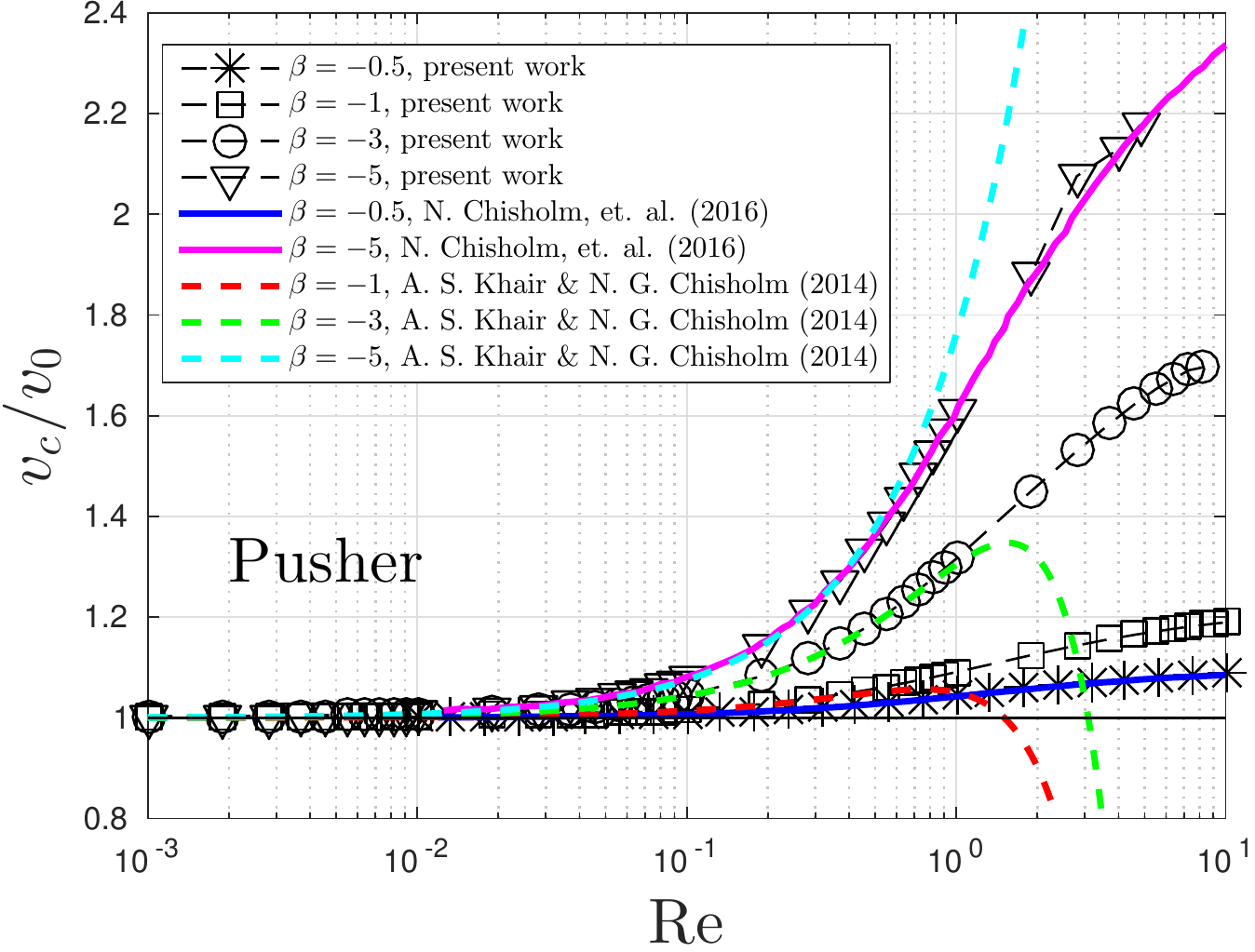}
%\end{subfigure}
\end{center}
\caption{Translational speed of the spherical squirmer, normalized with that of a neutral one ($v_{0} = \frac{2}{3}B_{1}$, for all $\Rey$), for different Reynolds numbers and different values of $\beta$. Also shown are the results of Chisholm et. al. \cite{chisholm2016squirmer} (solid lines) and of the asymptotic expansion of Khair et. al. \cite{khair2014expansions} (dashed lines).} \label{fig:Re_cmp}
\end{figure}

We simulated squirmers corresponding to $\beta = \pm 0.5, \pm 1, \pm 3, \pm 5$ with $\Rey$ ranging from $10^{-3}$ to $10^{1}$. The mesh was selected fine enough to obtain mesh-independent results. Figure \ref{fig:Re_cmp} shows the obtained behavior of $v_{c}$, normalized with the speed of a neutral squirmer $v_{0} = \frac{2}{3}B_{1}$, as a function of $\Rey$ compared to previous results in the literature. The excellent agreement serves as verification of the model and code.

For completeness, Figures \ref{fig:sqrmodelbp05}-\ref{fig:sqrmodelbn5} illustrate the fluid variables for $\beta = \pm 0.5$, $\pm 5$ and $\Rey = 10^{-2}, 1, 10^{2}$. Shown are the pressure field and the streamlines in the laboratory frame and in a frame moving with the particle. It was impossible to attain convergence of the nonlinear solver for $\beta=-5$ beyond $\Rey \simeq 8$, so Figure \ref{fig:sqrmodelbn5} lacks the last plots. Good agreement is again observed with the results of Chisholm et. al., to which the reader is referred for further discussions on squirmer hydrodynamics at finite $\Rey$.

%%%%%%%%%%%%%%%%%%%%%%%%%%%%%%%%%%%%%%%%
\newcommand\ch{.25}
\begin{figure}[!t]
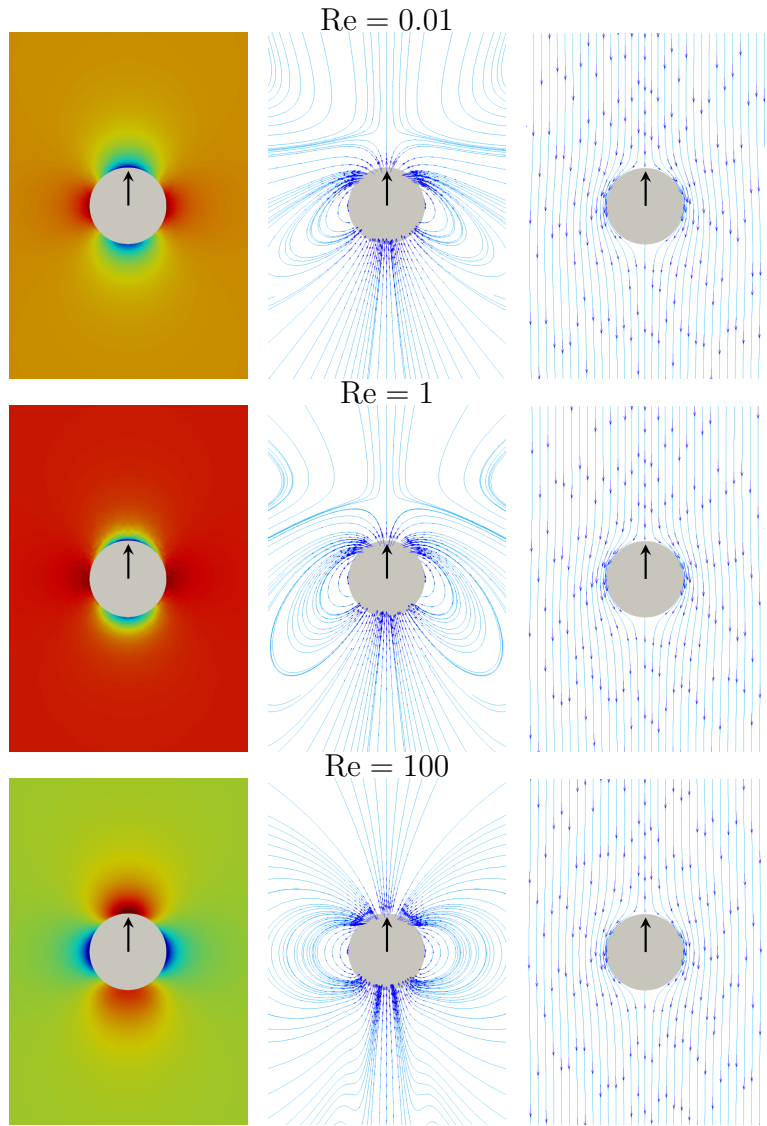

\begin{center}
    \begin{subfigure}{\textwidth}%%%%%%%%%%%%%%%%%%%%
    \begin{center}
    $\Rey = 0.01$\par
    \foreach \y in {re001,eulre001,lagre001}{
    \begin{tikzpicture}
    \node[anchor=south west,inner sep=0] (image) at (0,0) {\fbox{\includegraphics[width=\ch\textwidth,trim=0cm 0cm 0cm 0cm,clip]{betap05\y.pdf}}};
    \begin{scope}[x={(image.south east)},y={(image.north west)}]
    \draw[->,>=stealth,thick] (0.5,0.5) -- (0.5,0.6);
    \end{scope}
    \end{tikzpicture}
    }
    \end{center}
    \end{subfigure}
    \begin{subfigure}{\textwidth}%%%%%%%%%%%%%%%%%%%%
    \begin{center}
    $\Rey = 1$\par
    \foreach \y in {re1,eulre1,lagre1}{
    \begin{tikzpicture}
    \node[anchor=south west,inner sep=0] (image) at (0,0) {\fbox{\includegraphics[width=\ch\textwidth,trim=0cm 0cm 0cm 0cm,clip]{betap05\y.pdf}}};
    \begin{scope}[x={(image.south east)},y={(image.north west)}]
    \draw[->,>=stealth,thick] (0.5,0.5) -- (0.5,0.6);
    \end{scope}
    \end{tikzpicture}
    }
    \end{center}
    \end{subfigure}
    \begin{subfigure}{\textwidth}%%%%%%%%%%%%%%%%%%%%
    \begin{center}
    $\Rey = 100$\par
    \foreach \y in {re100,eulre100,lagre100}{
    \begin{tikzpicture}
    \node[anchor=south west,inner sep=0] (image) at (0,0) {\fbox{\includegraphics[width=\ch\textwidth,trim=0cm 0cm 0cm 0cm,clip]{betap05\y.pdf}}};
    \begin{scope}[x={(image.south east)},y={(image.north west)}]
    \draw[->,>=stealth,thick] (0.5,0.5) -- (0.5,0.6);
    \end{scope}
    \end{tikzpicture}
    }
    \end{center}
    \end{subfigure}
\end{center}
\caption{Pressure field (left, colors going from red to blue indicating maximum to minimum pressure values), streamlines in the laboratory frame (center) and streamlines in a frame moving with the squirmer (right) for $\beta = +0.5$ at different Reynolds numbers.} \label{fig:sqrmodelbp05}
\end{figure}
%%%%%%%%%%%%%%%%%%%%%%%%%%%%%%%%%%%%%%%%
\begin{figure}[!ht]
\begin{center}
    \begin{subfigure}{\textwidth}%%%%%%%%%%%%%%%%%%%%
    \begin{center}
    $\Rey = 0.01$\par
    \foreach \y in {re001,eulre001,lagre001}{
    \begin{tikzpicture}
    \node[anchor=south west,inner sep=0] (image) at (0,0) {\fbox{\includegraphics[width=\ch\textwidth,trim=0cm 0cm 0cm 0cm,clip]{betan05\y.pdf}}};
    \begin{scope}[x={(image.south east)},y={(image.north west)}]
    \draw[->,>=stealth,thick] (0.5,0.5) -- (0.5,0.6);
    \end{scope}
    \end{tikzpicture}
    }
    \end{center}
    \end{subfigure}
    \begin{subfigure}{\textwidth}%%%%%%%%%%%%%%%%%%%%
    \begin{center}
    $\Rey = 1$\par
    \foreach \y in {re1,eulre1,lagre1}{
    \begin{tikzpicture}
    \node[anchor=south west,inner sep=0] (image) at (0,0) {\fbox{\includegraphics[width=\ch\textwidth,trim=0cm 0cm 0cm 0cm,clip]{betan05\y.pdf}}};
    \begin{scope}[x={(image.south east)},y={(image.north west)}]
    \draw[->,>=stealth,thick] (0.5,0.5) -- (0.5,0.6);
    \end{scope}
    \end{tikzpicture}
    }
    \end{center}
    \end{subfigure}
    \begin{subfigure}{\textwidth}%%%%%%%%%%%%%%%%%%%%
    \begin{center}
    $\Rey = 100$\par
    \foreach \y in {re100,eulre100,lagre100}{
    \begin{tikzpicture}
    \node[anchor=south west,inner sep=0] (image) at (0,0) {\fbox{\includegraphics[width=\ch\textwidth,trim=0cm 0cm 0cm 0cm,clip]{betan05\y.pdf}}};
    \begin{scope}[x={(image.south east)},y={(image.north west)}]
    \draw[->,>=stealth,thick] (0.5,0.5) -- (0.5,0.6);
    \end{scope}
    \end{tikzpicture}
    }
    \end{center}
    \end{subfigure}
\end{center}
\caption{Same as Figure \ref{fig:sqrmodelbp05} for $\beta = -0.5$.} \label{fig:sqrmodelbn05}
\end{figure}
%%%%%%%%%%%%%%%%%%%%%%%%%%%%%%%%%%%%%%%%
\begin{figure}[!ht]
\begin{center}
    \begin{subfigure}{\textwidth}%%%%%%%%%%%%%%%%%%%%
    \begin{center}
    $\Rey = 0.01$\par
    \foreach \y in {re001,eulre001,lagre001}{
    \begin{tikzpicture}
    \node[anchor=south west,inner sep=0] (image) at (0,0) {\fbox{\includegraphics[width=\ch\textwidth,trim=0cm 0cm 0cm 0cm,clip]{betap5\y.pdf}}};
    \begin{scope}[x={(image.south east)},y={(image.north west)}]
    \draw[->,>=stealth,thick] (0.5,0.5) -- (0.5,0.6);
    \end{scope}
    \end{tikzpicture}
    }
    \end{center}
    \end{subfigure}
    \begin{subfigure}{\textwidth}%%%%%%%%%%%%%%%%%%%%
    \begin{center}
    $\Rey = 1$\par
    \foreach \y in {re1,eulre1,lagre1}{
    \begin{tikzpicture}
    \node[anchor=south west,inner sep=0] (image) at (0,0) {\fbox{\includegraphics[width=\ch\textwidth,trim=0cm 0cm 0cm 0cm,clip]{betap5\y.pdf}}};
    \begin{scope}[x={(image.south east)},y={(image.north west)}]
    \draw[->,>=stealth,thick] (0.5,0.5) -- (0.5,0.6);
    \end{scope}
    \end{tikzpicture}
    }
    \end{center}
    \end{subfigure}
    \begin{subfigure}{\textwidth}%%%%%%%%%%%%%%%%%%%%
    \begin{center}
    $\Rey = 100$\par
    \foreach \y in {re100,eulre100,lagre100}{
    \begin{tikzpicture}
    \node[anchor=south west,inner sep=0] (image) at (0,0) {\fbox{\includegraphics[width=\ch\textwidth,trim=0cm 0cm 0cm 0cm,clip]{betap5\y.pdf}}};
    \begin{scope}[x={(image.south east)},y={(image.north west)}]
    \draw[->,>=stealth,thick] (0.5,0.5) -- (0.5,0.6);
    \end{scope}
    \end{tikzpicture}
    }
    \end{center}
    \end{subfigure}
\end{center}
\caption{Same as Figure \ref{fig:sqrmodelbp05} for $\beta = +5$.} \label{fig:sqrmodelbp5}
\end{figure}
%%%%%%%%%%%%%%%%%%%%%%%%%%%%%%%%%%%%%%%%
\begin{figure}[!ht]
\begin{center}
    \begin{subfigure}{\textwidth}%%%%%%%%%%%%%%%%%%%%
    \begin{center}
    $\Rey = 0.01$\par
    \foreach \y in {re001,eulre001,lagre001}{
    \begin{tikzpicture}
    \node[anchor=south west,inner sep=0] (image) at (0,0) {\fbox{\includegraphics[width=\ch\textwidth,trim=0cm 0cm 0cm 0cm,clip]{betan5\y.pdf}}};
    \begin{scope}[x={(image.south east)},y={(image.north west)}]
    \draw[->,>=stealth,thick] (0.5,0.5) -- (0.5,0.6);
    \end{scope}
    \end{tikzpicture}
    }
    \end{center}
    \end{subfigure}
    \begin{subfigure}{\textwidth}%%%%%%%%%%%%%%%%%%%%
    \begin{center}
    $\Rey = 1$\par
    \foreach \y in {re1,eulre1,lagre1}{
    \begin{tikzpicture}
    \node[anchor=south west,inner sep=0] (image) at (0,0) {\fbox{\includegraphics[width=\ch\textwidth,trim=0cm 0cm 0cm 0cm,clip]{betan5\y.pdf}}};
    \begin{scope}[x={(image.south east)},y={(image.north west)}]
    \draw[->,>=stealth,thick] (0.5,0.5) -- (0.5,0.6);
    \end{scope}
    \end{tikzpicture}
    }
    \end{center}
    \end{subfigure}
\end{center}
\caption{Same as Figure \ref{fig:sqrmodelbp05} for $\beta = -5$.} \label{fig:sqrmodelbn5}
\end{figure}
%%%%%%%%%%%%%%%%%%%%%%%%%%%%%%%%%%%%%%%%

%%%%%%%%%%%%%%%%%%%%%%%%%%%%%%%%%%%%%%%%%%%%%%%%%%
%%%%%%%%%%%%%%%%%%%%%%%%%%%%%%%%%%%%%%%%%%%%%%%%%%

\section{Simulation of metachronal waves}
The cilia on the surface of a ciliated microorganism beat in a strongly organized fashion so as to break time-reversal symmetry and achieve propulsion \cite{childress1981mechanics,elgeti2013emergence}. This organization often takes the form of {\em metachronal waves}, which have been observed in \textit{Paramecium} \cite{machemer1972ciliary,tamm1972ciliary,funfak2015paramecium,ishikawa2006interaction}, \textit{Opalina ranarum} %\cite{sleigh1960form,tamm1970relation} \textit{Clamydomonas} 
\cite{goldstein2009noise,short2006flows,tamm1970relation,sleigh1960form} and flagellated \textit{Volvox algae} \cite{drescher2010fidelity,brumley2012hydrodynamic,brumley2015metachronal}.
This section discusses the implementation of metachronal waves in the ciliary envelope model discussed in this article, limiting the movements to tangential as before. Further, to simplify the exposition and render the problem two-dimensional, the tips of the cilia and the waves are assumed to move along meridian lines.

Let $s$ be the arc-length coordinate along a given meridian line. Notice that $s$ identifies a unique point $\mathbf{X}(s)$ in the reference configuration and also a unique material point on the surface of the organism's body. The tip of the cilium with its attachment point at $s$ is assumed to occupy, at time $t$, a position with arc-length coordinate denoted by $w$. We consider metachronal waves in which $s$, $w$ and $t$ are related by
\[%\begin{equation}\label{eqn:metpos}
w = s + A(s)\,\cos (ks-\omega t)~,\qquad \qquad (A\geq 0),
\]%\end{equation}
where $A$ is the amplitude of the displacement of the tip, $k = \frac{2\pi}{\lambda}$ is the spatial frequency of the wave (or wave number), $\lambda$ is the wavelength, $\omega = \frac{2\pi}{T}$ is the angular frequency and $T$ the period. For each $t$, the function $s\to w\paren{s,t}$ should be a non-decreasing function, otherwise the cilia experience tangential overlapping.

The velocity of the tip of the cilium attached at $s$ is
\[%\begin{equation}
    v(s,t) = \omega A(s)\,\sin (ks-\omega t),
\]%\end{equation}
but notice that this velocity does not take place at the point $s$ of the ciliary envelope, but rather at the point $w$. The boundary condition to be imposed at a given point of $\partial\mathcal{B}$ depends on the velocity of the ciliary envelope {\em at that point}. 
For this reason, to compute the velocity $u_{\mbox{\scriptsize{env}}}^i$ of the ciliary envelope at a mesh node $i\in\eta^U_\partial$ with arc-length coordinate $w_i$, at time $t_n$, one proceeds as follows:
\begin{enumerate}
\item Solve $F(s_i)=s_i+A(s_i)\,\cos(ks_i-\omega t_n)-w_i=0$ for $s_i$.
\item Compute $u_{\mbox{\scriptsize{env}}}^i=\omega A(s_i) \sin (ks_i-\omega t_n)$.
\end{enumerate}

For a type-I squirmer one assumes that the fluid has the same velocity as the ciliary envelope, and thus the interface condition is $\mathbf{u}_s=u_{\mbox{\scriptsize{env}}} \boldsymbol{\tau}_b$.
This type-I behavior, however, is only realistic when the cilia are very densely distributed over the body surface. In general one would expect a drag law between the cilia and the adjacent fluid, which can be modeled by a type-II squirmer with the law
\[%\begin{equation}
    \mathbf{f}_s = C_D \frac{\mu}{L} \paren{u_{\mbox{\scriptsize{env}}}-\mathbf{u}_s\cdot \boldsymbol{\tau}_b} \boldsymbol{\tau}_b,
\]%\end{equation}
where for large values of the non-dimensional drag coefficient $C_D$ one recovers the type-I behavior.

It should be noted that the cilia are typically much smaller than the body length $L$, and $A(s)$ is of the order of the cilium length $\ell$. As a consequence, to first order in $\ell/L$, one has $u_{\mbox{\scriptsize{env}}}^i(t)\simeq v(w_i,t)=\omega A(w_i)\,\sin (kw_i-\omega t)$. This is easier to implement, since it avoids the nonlinear problem in step (1) above. However, {\em this first-order approximation makes the translational velocity of the squirmer to be zero}. The propulsion by metachronal waves is a second-order effect.

As illustrative example we report here simulations of two-dimensional bodies inspired in \textit{Opalina ranarum}. We used the prolate spheroidal shape proposed by Zhang et. al. \cite{zhang2015paramecia} for ciliated organisms, which in the $x-y$ plane reads
\[
y\paren{x} = b\paren{1-\frac{x^{2}}{a^{2}}}^{\frac{1}{2}} - \varepsilon\sin\paren{\pi\frac{x}{a}},
\]
where $\varepsilon \geq 0$ is a parameter of asymmetry perturbing an ellipse with $a$ and $b$ semi-major and semi-minor axes, respectively. The adopted values are $a = 110 \mu$m,
$b = 36.3 \mu$m and $\varepsilon = 0.09 b$. The metachronal wavelength was taken with $\lambda=50 \mu$m ($k=2\pi/\lambda)$ and frequency of $5$ beats per second ($\omega=10 \pi$ rad/s). 
Finally, the amplitude function for the tangential displacement of the ciliary envelope was taken as 
\[%\begin{equation}
A(s)=K\,\tanh (\eta\,\sin(\pi\frac{s}{L}))
\]%\end{equation}
where $\eta = 5$ and $L=324 \mu$m (the semi-perimeter of the body's boundary). Essentially, $A(s)=K$ except for $s<0.1 L$ and $s>0.9 L$, where it smoothly tends to zero. The constant $K$ was taken as $K=0.02 L=6.5 \mu$m, which is consistent with a typical cilium length of $10 \mu$m. The fluid's properties were taken as $\rho=0$, $\mu=10^{-3}$ Pa-s.

\begin{figure}[t]
\centering
\includegraphics[width = .7\textwidth]{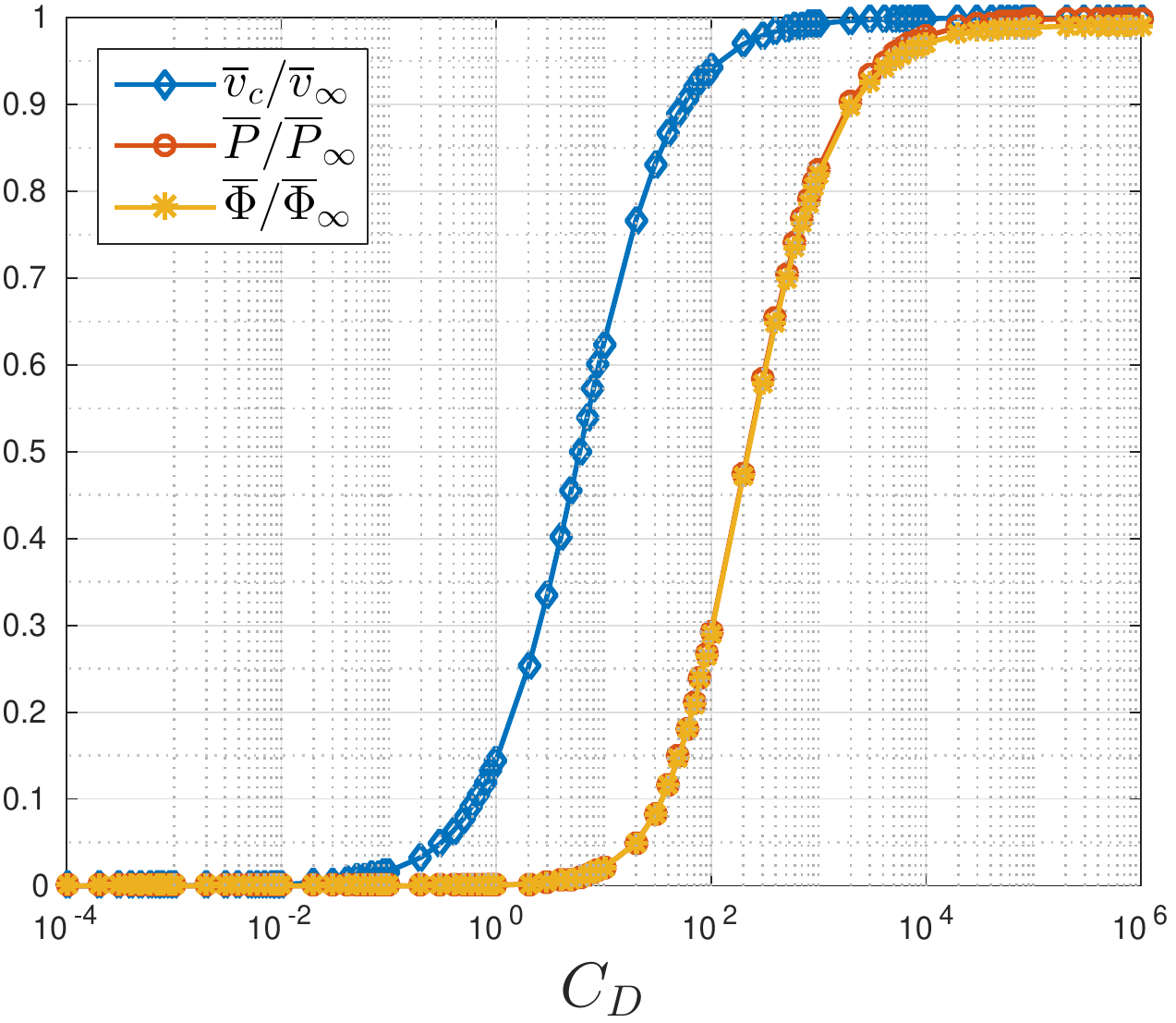}
\caption{Dimensionless time-averaged speed $\overline{v}_{c}$ and power expenditure $\overline{P}$ of the \textit{Opalina ranarum} model, and viscous dissipation of the fluid $\overline{\Phi}$ as functions of the drag coefficient $C_{D}$. The variables are scaled with those corresponding to the type-I squirmer (equivalent to $C_D\to +\infty$), denoted by $\overline{v}_{\infty}$, $\overline{P}_{\infty}$ and $\overline{\Phi}_{\infty}$.} \label{fig:vmed_cmp}
\end{figure}

In Figure \ref{fig:vmed_cmp} we show the velocity of the model as a function of the drag coefficient $C_D$. The maximum value is attained to the type-I squirmer, giving an average velocity $\overline{v}_{\infty} = 48.8 \mu$m/s. As $C_D\to +\infty$ (above, say, $10^3$) the type-II squirmer tends to this value, while for very little drag ($C_D < 10^{-1}$) the model practically does not move. Velocities of about 50 $\mu$m/s are within the observed value in {\em Opalina}. 

Also shown in the figure is the average power consumption $\overline{P}$, defined as the time average of
\[%\begin{equation}
    {P}= \int_{\partial\mathcal{B}_h} \mathbf{f}_s\cdot \mathbf{u}_s,
\]%\end{equation}
and the fluid's average viscous dissipation $\overline{\Phi}$, defined as the time average of
\[%\begin{equation}
    \Phi=\int_{\Omega_{fh}}2\mu\nabla^S\mathbf{u}:\nabla^S\mathbf{u}~d\mathbf{x}.
\]%\end{equation}
These two quantities should only differ by virtue of numerical dissipation, which is very low because the mesh is highly refined, as evidenced in Figure \ref{fig:vmed_cmp}. %, but increased ill-conditioning of the system makes them to differ for $C_D>10^3$. 
They are shown scaled by the values that correspond to the type-I squirmer.
%%, which is $\overline{P}_{\infty}= \overline{\Phi}_{\infty} = 2393.2$.

Remarkably, in the range $10<C_D<100$ the type-II squirmer attains a velocity comparable to that of the type-I squirmer with much less power expenditure. This is probably a beneficial effect of some slippage between the ciliary envelope and the fluid, reducing velocity gradients in the latter without modifying the overall hydrodynamic pattern that generates self-propulsion. This is an intriguing topic that adds interest to the simulation of type-II squirmers, but is outside the scope of this paper.

\begin{figure}[t]
\centering
\includegraphics[width = .8\textwidth]{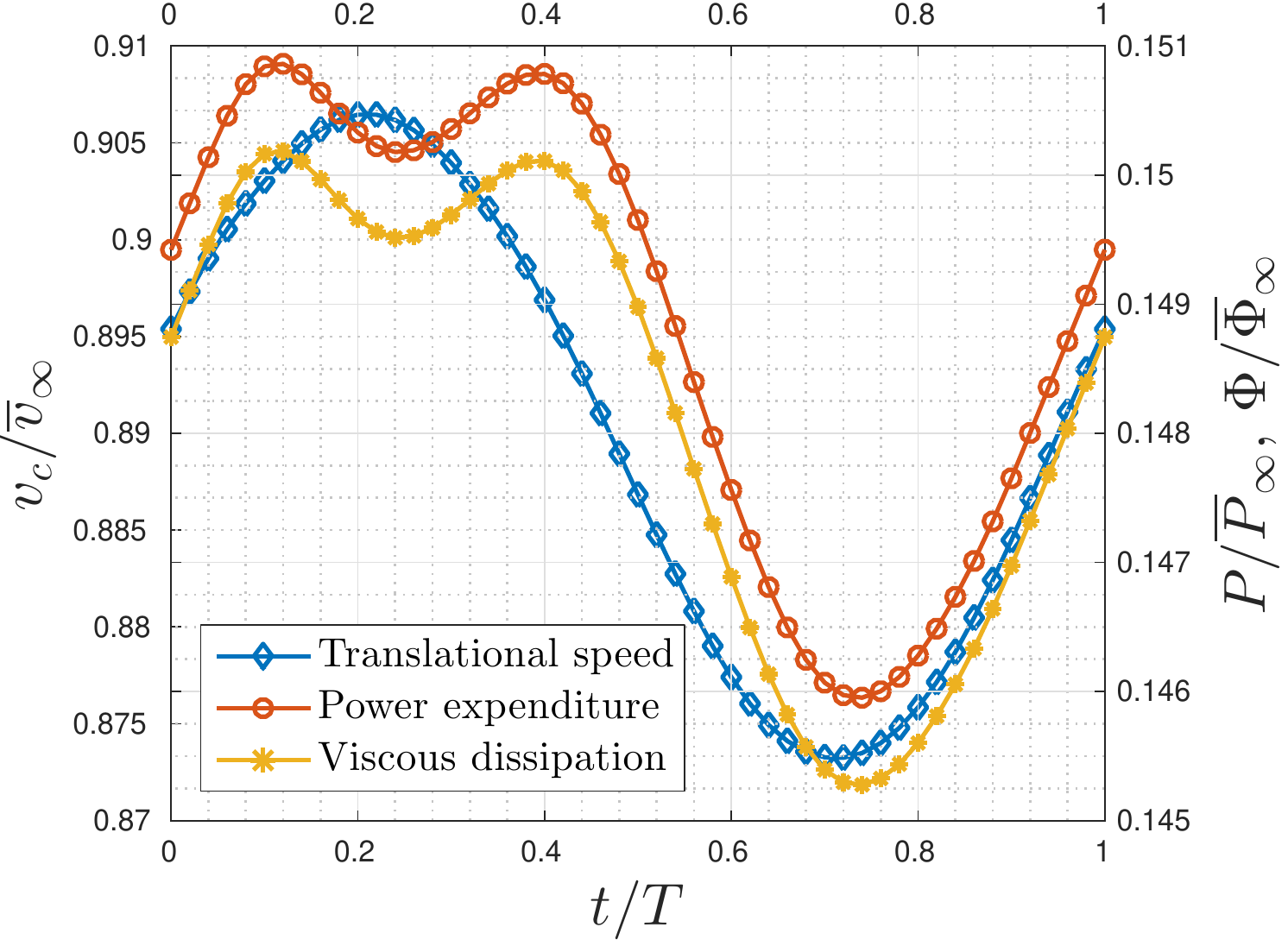}
\caption{Dimensionless speed ${v}_{c}(t)$, power expenditure ${P}(t)$ and viscous dissipation $\Phi(t)$ of the \textit{Opalina ranarum} model for the case $C_{D} = 50$. The nondimensionalization is as in Figure \ref{fig:vmed_cmp}.} \label{fig:vmed50_cmp}
\end{figure}

For completeness, some results of the model of {\em Opalina ranarum} swimming by itself are also included, now concentrating in the intermediate value $C_D=50$. Figure \ref{fig:vmed50_cmp} shows the instantaneous values of $v_c$, $P$ and $\Phi$ along one metachronal period. As advanced, the difference $P-\Phi$ is less than $0.5$\%, reflecting that numerical dissipation is indeed very low. The max/min ratio in velocity is about 1.04, from which one concludes that though the metachronal wavelength is rather large (about 1/6 of the semi-perimeter), it is short enough to produce an essentially constant $v_c$. Figure \ref{fig:opalmodel50} shows the streamlines at times $t/T=0.2$ and 0.72, roughly corresponding to the instants of maximum and minimum $v_c$. The streamlines in the laboratory frame show a stagnation point in the front of the squirmer, thus characterizing it as a {\em puller}. Close images of the corresponding pressure and velocity fields are shown in  Figure \ref{fig:opalmodel50zoom}.

%%%%%%%%%%%%%%%%%%%%%%%%%%%%%%%%%%%%%%%%
\newcommand\cl{.38}
\begin{figure}[!ht]
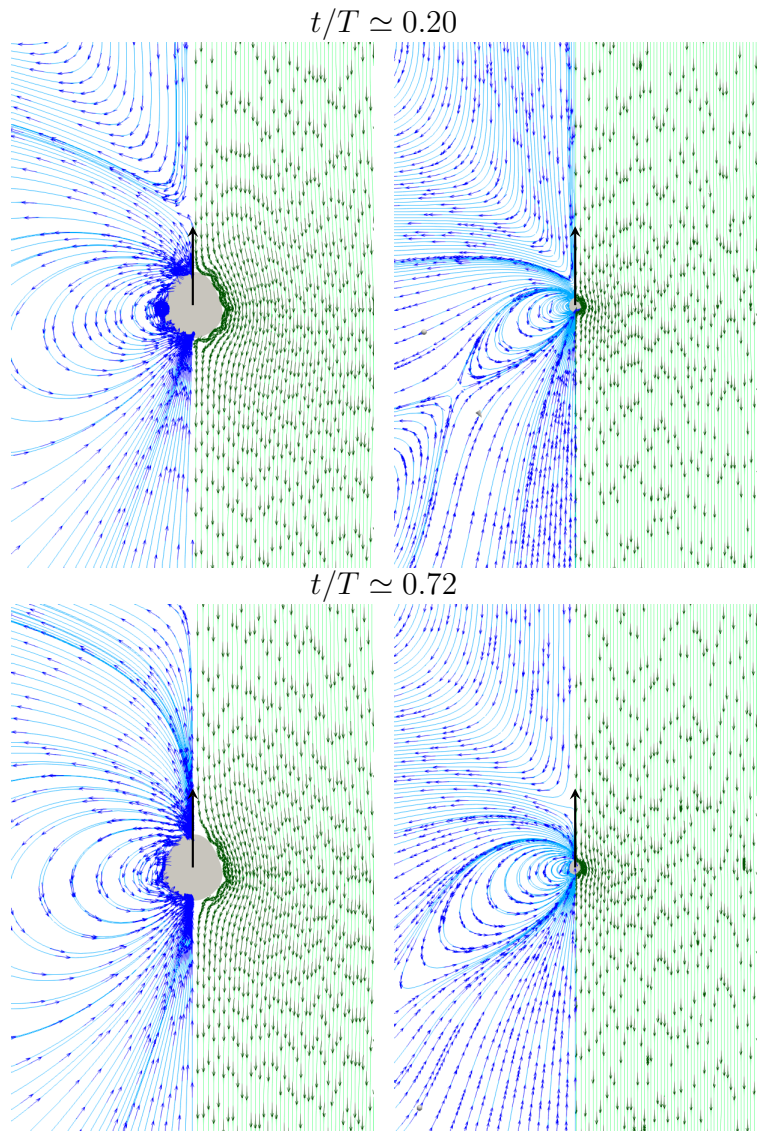

\begin{center}
    \begin{subfigure}{\textwidth}%%%%%%%%%%%%%%%%%%%%
    \begin{center}
    $t/T \simeq 0.20$\par
    \foreach \y in {vel,vell}{
    \begin{tikzpicture}
    \node[anchor=south west,inner sep=0] (image) at (0,0) {\fbox{\includegraphics[width=\cl\textwidth,trim=0cm 0cm 0cm 0cm,clip]{opal_\y_tmax04.pdf}}};
    \begin{scope}[x={(image.south east)},y={(image.north west)}]
    \draw[->,>=stealth,thick] (0.5,0.5) -- (0.5,0.65);
    \end{scope}
    \end{tikzpicture}
    }
    \end{center}
    \end{subfigure}
    \begin{subfigure}{\textwidth}%%%%%%%%%%%%%%%%%%%%
    \begin{center}
    $t/T \simeq 0.72$\par
    \foreach \y in {vel,vell}{
    \begin{tikzpicture}
    \node[anchor=south west,inner sep=0] (image) at (0,0) {\fbox{\includegraphics[width=\cl\textwidth,trim=0cm 0cm 0cm 0cm,clip]{opal_\y_tmin14.pdf}}};
    \begin{scope}[x={(image.south east)},y={(image.north west)}]
    \draw[->,>=stealth,thick] (0.5,0.5) -- (0.5,0.65);
    \end{scope}
    \end{tikzpicture}
    }
    \end{center}
    \end{subfigure}
\end{center}
\caption{Streamlines in laboratory frame (blue) and streamlines in a frame moving with the squirmer (green) for the \textit{Opalina ranarum} model with $C_D=50$. Two zooms are shown to illustrate the near and far fields at two different times ($t/T \simeq 0.20$ and $t/T\simeq 0.72$). The black arrow shows the swimming direction.} \label{fig:opalmodel50}
\end{figure}
%%%%%%%%%%%%%%%%%%%%%%%%%%%%%%%%%%%%%%%%

\begin{figure}[!ht]
\centering
\includegraphics[width = .7\textwidth]{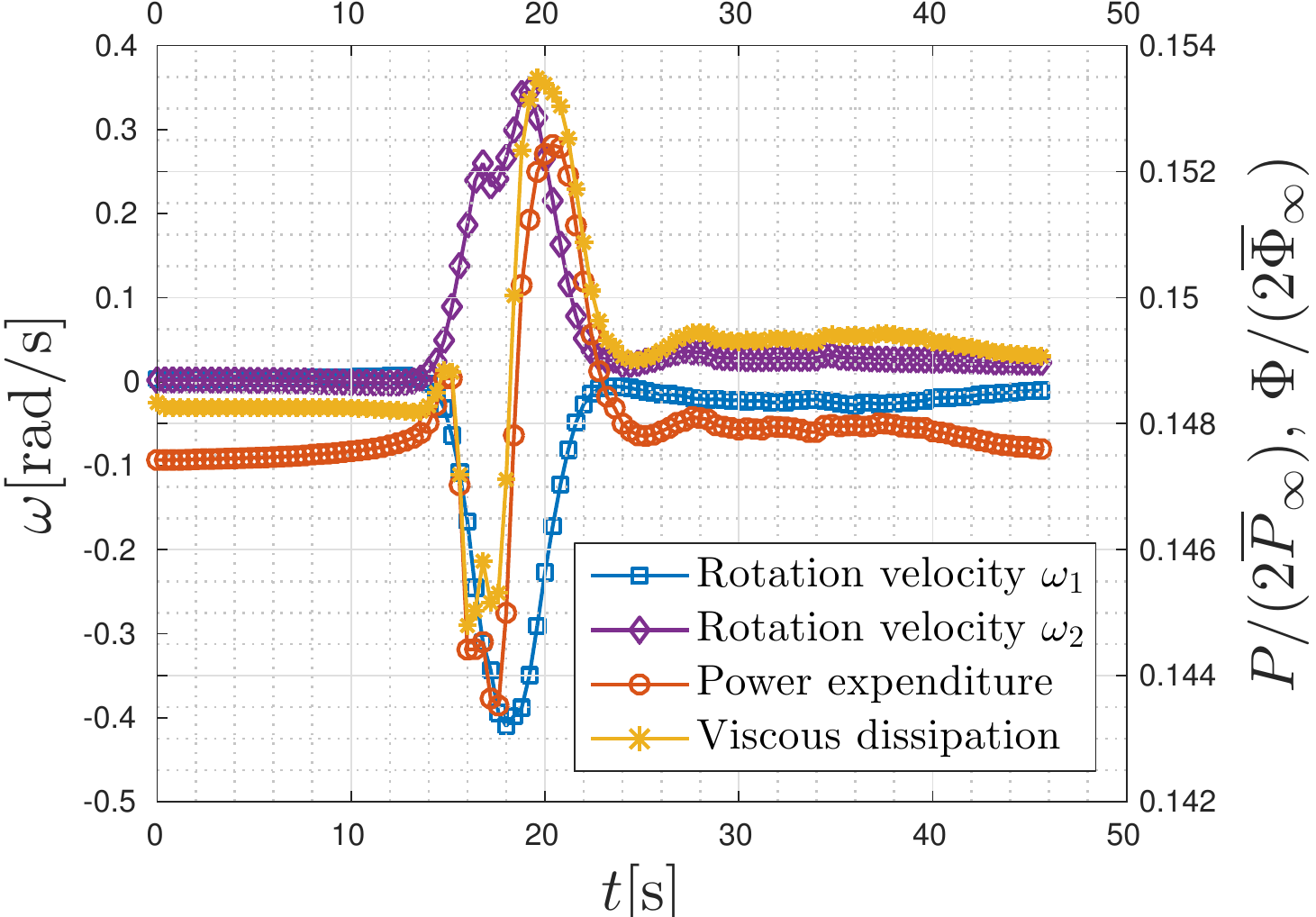}
\caption{Period-averaged rotational velocities $\omega_{1}(t)$ and $\omega_{2}(t)$ (subscripts representing each squirmer), dimensionless power expenditure ${P}(t)$ and viscous dissipation $\Phi(t)$ of the \textit{Opalina ranarum} interaction model for the case $C_{D} = 50$. The nondimensionalization is as in Figure \ref{fig:vmed_cmp} considering that there are now two squirmers.} \label{fig:temp_cmp}
\end{figure}

\begin{figure}[!ht]
\begin{center}
    \includegraphics[width=.75\textwidth,trim=0cm 0cm 0cm 0cm,clip]{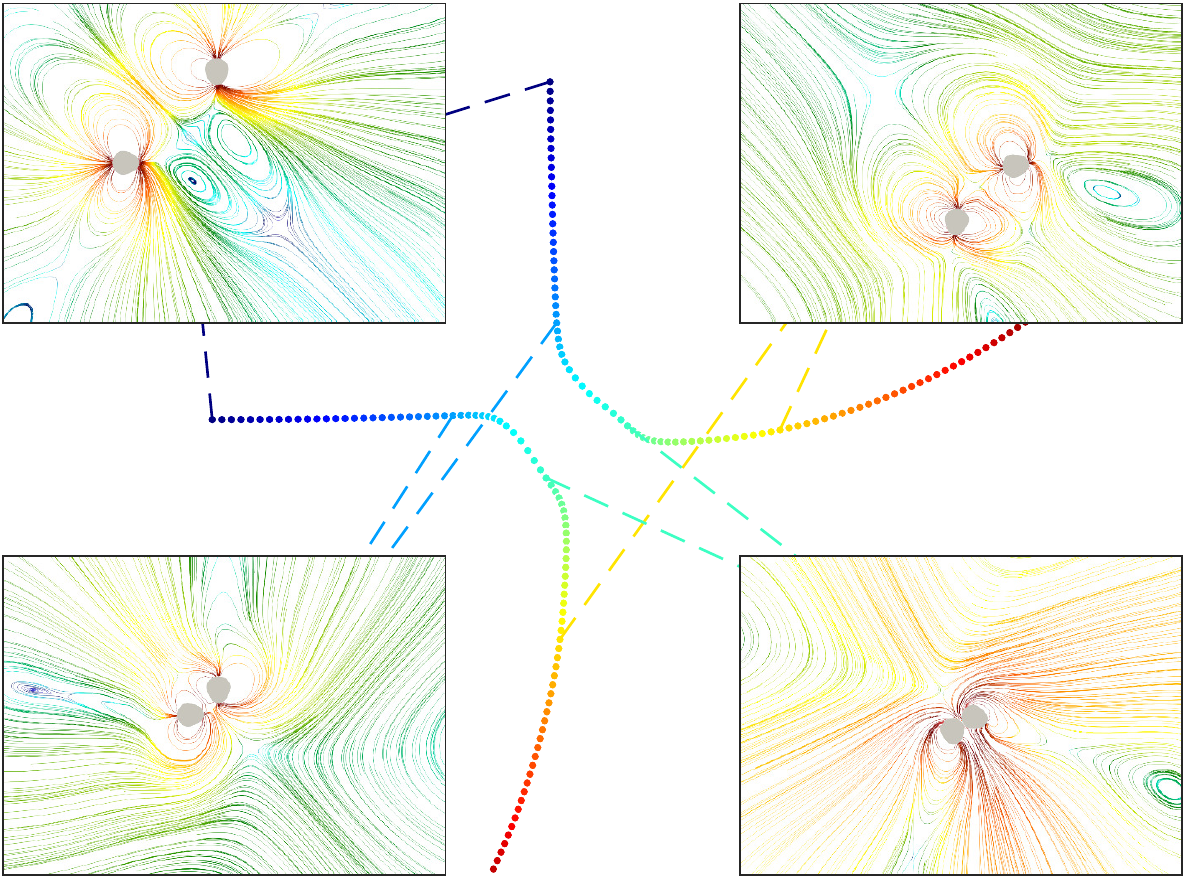}
    \hspace{.5cm}
    \includegraphics[width=.1\textwidth,trim=0cm 0cm 0cm 0cm,clip]{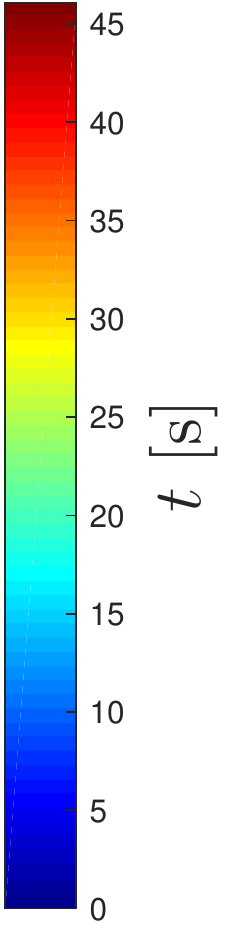}
\end{center}
\caption{Trajectories of two \textit{Opalina ranarum} in hydrodynamic interaction. The colored sequence of points indicates the position of the centroid of each \textit{Opalina ranarum} from $t = 0$ s (blue) to $t = 46$ s (red). Streamlines are shown for specific times evidencing the changes in the fluid velocity field throughout the interaction process. The streamlines are colored with the norm of the velocity, in logarithmic scale.}\label{fig:2opaltrj}
\end{figure}

Finally, a two-dimensional simulation of two \textit{Opalina ranarum} individuals interacting is presented, again for the intermediate case $C_{D} = 50$. Their unperturbed initial directions of locomotion are orthogonal, in such a way that they eventually meet and interact hydrodynamically (no contact or repulsion force has been added). The detailed velocity and pressure fields during the interaction process are shown in Figure \ref{fig:sqrmodel50int} for instants of time between $t = 13$ s and $t = 23$ s (approximately the time that takes the strong interaction), evidencing the mutually induced change of orientation. A global view of the interaction is presented in Figure \ref{fig:2opaltrj} in which the first swimmer, going from the left to the right, encounters the second one, going from top to bottom. Streamlines are also shown at four specific times depicting the emergence of recirculating regions and stagnation points during the interaction process. In a similar way as in Figure \ref{fig:vmed50_cmp}, period-averaged $P$ and $\Phi$ are plotted in Figure \ref{fig:temp_cmp} along the whole time dependent collision progress, together with the rotational velocities $\omega_{1}$ and $\omega_{2}$ of the first and second swimmer, respectively. The numerical dissipation is low (less than 1\%) throughout the simulation. In the course of the strong interaction, the first and second swimmers attain rotational velocities of $-0.45$rad/s and $0.35$rad/s, respectively.

%%%%%%%%%%%%%%%%%%%%%%%%%%%%%%%%%%%%%%%%
\begin{figure}[t]
\begin{center}
    \begin{subfigure}{\textwidth}%%%%%%%%%%%%%%%%%%%%
    \begin{center}
    $t/T \simeq 0.20 \hspace{.32\textwidth} t/T \simeq 0.72$\par
    \foreach \y in {tmax04,tmin14}{
    \begin{tikzpicture}
    \node[anchor=south west,inner sep=0] (image) at (0,0) {\fbox{\includegraphics[width=.49\textwidth,trim=0cm 0cm 0cm 0cm,clip]{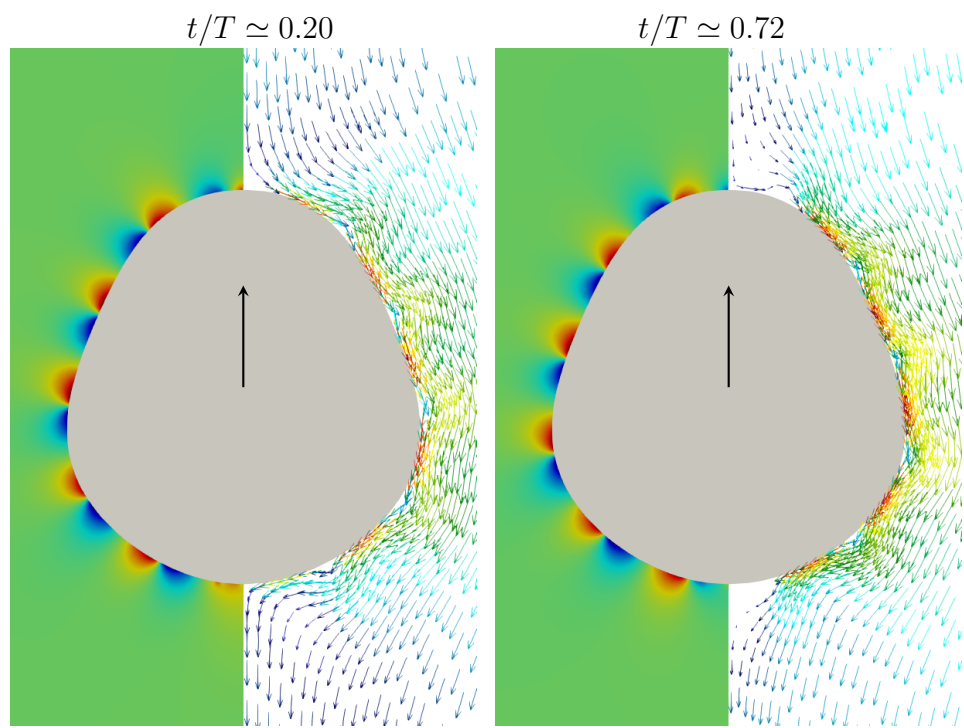}}};
    \begin{scope}[x={(image.south east)},y={(image.north west)}]
    \draw[->,>=stealth,thick] (0.5,0.5) -- (0.5,0.65);
    \end{scope}
    \end{tikzpicture}
    }
    \end{center}
    \end{subfigure}
\end{center}
\caption{Detailed view of pressure (left on each figure) and velocity fields (right on each figure) for the same model in Figure \ref{fig:opalmodel50} (in both cases, colors going from red to blue indicate maximum to minimum pressure and speed values).} \label{fig:opalmodel50zoom}% Streamlines in laboratory frame (blue) and streamlines in a frame moving with the squirmer (green) for the \textit{Opalina ranarum} model with $C_D=50$. Two zooms are shown to illustrate the near and far fields at two different times ($t/T \simeq 0.20$ and $t/T\simeq 0.72$). The black arrow shows the swimming direction.} 
\end{figure}
%%%%%%%%%%%%%%%%%%%%%%%%%%%%%%%%%%%%%%%%

%%%%%%%%%%%%%%%%%%%%%%%%%%%%%%%%%%%%%%%%
\newcommand{\ci}{.3}
\begin{figure}[!ht]
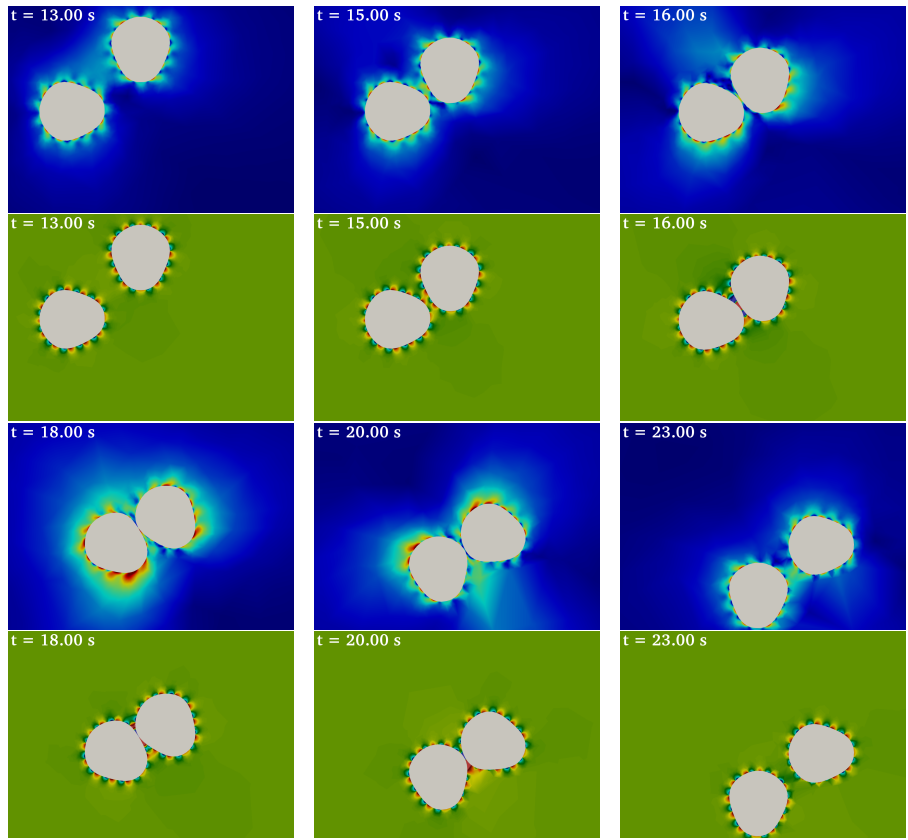

\begin{center}
    \foreach \x in {vel,pre}{
    \begin{subfigure}{\textwidth}%%%%%%%%%%%%%%%%%%%%
    \begin{center}
    \foreach \y in {1300,1500,1600}{
    \begin{tikzpicture}
    \node[anchor=south west,inner sep=0] (image) at (0,0) {\fbox{\includegraphics[width=\ci\textwidth,trim=0cm 0cm 0cm 0cm,clip]{2opal_\x_t\y.pdf}}};
    \end{tikzpicture}
    }
    \end{center}
    \end{subfigure}
    }
    \foreach \x in {vel,pre}{
    \begin{subfigure}{\textwidth}%%%%%%%%%%%%%%%%%%%%
    \begin{center}
    \foreach \y in {1800,2000,2300}{
    \begin{tikzpicture}
    \node[anchor=south west,inner sep=0] (image) at (0,0) {\fbox{\includegraphics[width=\ci\textwidth,trim=0cm 0cm 0cm 0cm,clip]{2opal_\x_t\y.pdf}}};
    \end{tikzpicture}
    }
    \end{center}
    \end{subfigure}
    }
\end{center}
\caption{Detailed view of fluid velocity (first and third rows) and pressure fields (second and fourth rows) during the interaction process (for each cases, colors going from red to blue indicate maximum to minimum speed and pressure values).} \label{fig:sqrmodel50int}
\end{figure}
%%%%%%%%%%%%%%%%%%%%%%%%%%%%%%%%%%%%%%%%

%%%%%%%%%%%%%%%%%%%%%%%%%%%%%%%%%%%%%%%%

%%%%%%%%%%%%%%%%%%%%%%%%%%%%%%%%%%%%%%%%%%%%%%%%%%

%%%%%%%%%%%%%%%%%%%%%%%%%%%%%%%%%%%%%%%%%%%%%%%%%%
%%%%%%%%%%%%%%%%%%%%%%%%%%%%%%%%%%%%%%%%%%%%%%%%%%
\clearpage
\section{Conclusions}
In this article the mathematical setting, numerical approximation and implementation details needed for successful simulation of organisms and active particles known as squirmers have been presented. The formulation describes two types (I and II) of interface conditions between the body and the surrounding fluid, of which the second is a versatile coupling-force condition that had been scarcely treated in the literature. The differential and weak formulations have been introduced for both types of squirmers in such a way that they are readily associated with formulations familiar to practitioners of computational fluid dynamics (CFD) that use finite element or finite volume codes. Special care was taken to provide the procedures that turn a generic CFD solver into a squirmer simulator. The techniques apply to squirmers of any shape, contemplate inertial or rheological nonlinearities, and can handle interactions of any number of simultaneous squirmers in domains of arbitrary geometry. Hopefully, this will encourage other researchers to implement the proposed techniques into open-source libraries and commercial codes. In this way, investigations of the fascinating individual and collective behavior of phoretic particles and living micro-organisms will become more accessible to students and specialists of other areas.

%%%%%%%%%%%%%%%%%%%%%%%%%%%%%%%%%%%%%%%%%%%%%%%%%%
%%%%%%%%%%%%%%%%%%%%%%%%%%%%%%%%%%%%%%%%%%%%%%%%%%

\section*{Acknowledgments}
The authors gratefully acknowledge the financial support received from S\~ao Paulo Research Foundation (FAPESP) (Grants \#2014/14720-8, \#2013/07375-0 (CEPID), \#2018/08752-5) and from the Brazilian National Council for Scientific and Technological Development (CNPq) (Grants \#305599/2017-8)
%(Grants \#308728/2013-0, \#447904/2014-0).

%%%%%%%%%%%%%%%%%%%%%%%%%%%%%%%%%%%%%%%%%%%%%%%%%%
%%%%%%%%%%%%%%%%%%%%%%%%%%%%%%%%%%%%%%%%%%%%%%%%%%

%\section*{References}
\bibliographystyle{amsplain}
\bibliography{/home/stevens/Documentos/USP/Decimo/paper2v/sqm_vol/refs.bib}

\end{document}